\documentclass[matprg]{svjour3}

\usepackage[top=2.5cm,bottom=2.5cm,right=2.5cm,left=2.5cm]{geometry}

\smartqed  
\usepackage{graphicx}
\usepackage{amsmath,amssymb,amstext}
\usepackage{mathrsfs}
\usepackage{multirow}
\usepackage{algorithm}
\usepackage{algorithmic}
\usepackage{threeparttable}
\usepackage{subfigure}
\usepackage[multiple]{footmisc}
\usepackage{color}
\makeatletter
\def\@thmcountersep{.}
\makeatother
\spnewtheorem{defi}{Definition}[section]{\bfseries}{\normalfont}
\spnewtheorem{alg}{Algorithm}[section]{\bfseries}{\normalfont}
\spnewtheorem{ass}{Assumption}[section]{\bfseries}{\normalfont}
\spnewtheorem{lem}{Lemma}[section]{\bfseries}{\normalfont}
\spnewtheorem{thm}{Theorem}[section]{\bfseries}{\normalfont}
\newcommand{\tsum}{\textstyle\sum}

\numberwithin{equation}{section}

\begin{document}
	
	\title{A Unified Single-loop Alternating Gradient Projection Algorithm for Nonconvex-Concave and Convex-Nonconcave Minimax Problems\thanks{Z. Xu was supported by National Natural Science Foundation of China under the grant 12071279 and by General Project of Shanghai Natural Science Foundation
			(No. 20ZR1420600). G. Lan's research was partly supported by National Science Foundation (NSF) grant CCF 1909298.}
	}
	
	\titlerunning{A Unified Single-loop AGP Algorithm for Minimax Problems}        
	
	\author{Zi Xu       \and
		Huiling Zhang \and
		Yang Xu \and Guanghui Lan
	}
	
	
	\institute{Zi Xu\at
		Department of Mathematics, College of Sciences, Shanghai University, Shanghai 200444, P.R.China. \\
		\email{xuzi@shu.edu.cn}           
		\and Huiling Zhang \at
		Department of Mathematics,  College of Sciences, Shanghai University, Shanghai 200444, P.R.China.\\ \email{zhl18720009@i.shu.edu.cn}
		\and Yang Xu\at
		Department of Mathematics, College of Sciences, Shanghai University, Shanghai 200444, P.R.China. \\ \email{fgjfhg3013@i.shu.edu.cn}
		\and Guanghui Lan \at
		Industrial and Systems Engineering, Georgia Institute of Technology, Atlanta, GA, 30332, USA.\\
		Corresponding author. \email{george.lan@isye.gatech.edu}
	}
	\date{Received: date / Accepted: date}

	\maketitle
	
	\begin{abstract}
		Much recent research effort has been directed to the development of efficient algorithms for solving minimax problems with theoretical convergence guarantees	due to the relevance of these problems to a few emergent applications.
		In this paper, we propose a unified single-loop alternating gradient projection  (AGP) algorithm for solving smooth nonconvex-(strongly) concave and (strongly) convex-nonconcave minimax problems. AGP employs simple gradient projection steps for updating the primal and dual variables alternatively at each iteration. We show that it can find an $\varepsilon$-stationary point of the	objective function in $\mathcal{O}\left( \varepsilon ^{-2} \right)$ (resp. $\mathcal{O}\left( \varepsilon ^{-4} \right)$) iterations under nonconvex-strongly concave (resp. nonconvex-concave) setting. Moreover, its gradient complexity to obtain an $\varepsilon$-stationary point of the objective function is bounded by $\mathcal{O}\left( \varepsilon ^{-2} \right)$ (resp., $\mathcal{O}\left( \varepsilon ^{-4} \right)$) under the strongly convex-nonconcave (resp.,  convex-nonconcave) setting. To the best of our knowledge, this is the first time that a simple and unified single-loop algorithm is developed for solving both nonconvex-(strongly) concave and (strongly) convex-nonconcave minimax problems. Moreover, the complexity results for solving the latter (strongly) convex-nonconcave minimax problems have never been obtained before in the literature. Numerical results show the efficiency of the proposed AGP algorithm.
		
		Furthermore, we extend the AGP algorithm by presenting a block alternating proximal gradient (BAPG) algorithm for solving more general multi-block nonsmooth nonconvex-(strongly) concave and (strongly) convex-nonconcave minimax problems. We can similarly establish the gradient complexity of the proposed algorithm under these four different settings.
	\end{abstract}
	
	\section{Introduction}
	\label{intro}
	We consider the following minimax optimization problem:
	\begin{equation}\label{problem:1}
		\min \limits_{x\in \mathcal{X}}\max \limits_{y\in \mathcal{Y}}f(x,y),\tag{P}
	\end{equation}
	where $\mathcal{X}\subseteq \mathbb{R}^n$ and $\mathcal{Y}\subseteq \mathbb{R}^m$ are nonempty closed and bounded convex sets, and $f: \mathcal{X} \times \mathcal{Y} \rightarrow \mathbb{R}$
	is a smooth function. This problem has attracted more  attention due to its wide applications in machine learning, signal processing, and many other research fields in recent years.
	Many practical problems can be formulated as in~\eqref{problem:1}, such as the power control and transceiver design problem in signal processing \cite{Lu},  distributed
	nonconvex optimization \cite{Mateos,Liao,Giannakis}, robust learning over multiple domains \cite{Qian}, statistical learning \cite{Abadeh,Giordano} and many others.
	
	Minimax optimization problems have been studied for many years, but most previous works focused on convex-concave minimax
	problems, i.e., $f(x,y)$ is convex with respect to $x$ and concave with respect to $y$~\cite{Nedic,Boyd,Chen}. Under this setting, Nemirovski~\cite{Nemi2004} proposed a mirror-prox algorithm which returns an $\varepsilon$-saddle point within the complexity of $\mathcal{O} (1/\varepsilon)$ when $\mathcal{X}$ and  $\mathcal{Y}$ are bounded sets. Nesterov~\cite{Nes2007} developed a dual extrapolation algorithm which owns the same complexity bound as in~\cite{Nemi2004}.  Monteiro and Svaiter~\cite{Mon2010,Mon2011} extended the complexity result to unbounded sets and composite objectives by using the hybrid proximal extragradient algorithm with a different termination criterion. Tseng~\cite{Tseng2008} proved the same result using a refined convergence analysis. Abernethy et al.~\cite{Abernethy} presented a Hamiltonian gradient descent algorithm with last-iterate convergence under a ``sufficient bilinear" condition. A few other papers have studied special cases in the convex-concave setting, for more details, we refer to~\cite{Chen,Chen2017,Dang,He2016,Lan2016,Lin2020,Ouyang2015,Ouyang2019} and the references therein.
	
	However, recent applications of problem~\eqref{problem:1}
	in machine learning and signal processing urge the necessity of moving beyond this classical setting.
	For example, in a typical generative adversary network (GANs) problem formulation \cite{Maziar}, the objective function $f(x,y)$ is nonconvex with respect to (w.r.t.) $x$ and concave w.r.t. $y$.
	As another example, in a robust support vector machine problem, we have
	$y \equiv (y_u, y_v) \in \mathbb{R}^m \times \mathbb{R}$, $x \equiv (x_w, x_b)\in \mathbb{R}^m \times \mathbb{R}$
	and the objective function $f(x, y)= y_v \left[\langle x_w, y_u\rangle + x_b\right]$. This problem is
	convex  w.r.t. $x$, but not necessarily concave w.r.t. $y$.
	
	Most of recent studies focus on nonconvex-(strongly) concave minimax problems. For nonconvex-strongly concave minimax problem, several recent works~\cite{Jin,Rafique,Lin2019,Lu} have studied various algorithms,
	and all of them can achieve the gradient complexity of $\tilde{\mathcal{O}}\left( \kappa_y^2\varepsilon ^{-2} \right)$ in terms of stationary point of $\Phi (\cdot) = \max_{y\in \mathcal{Y}} f(\cdot, y)$ (when $\mathcal{X}=\mathbb{R}^n$), or stationary point of $f$,  where $\kappa_y$ is the condition number for $f(x,\cdot)$. Lin et al.\cite{Lin2020} propose an accelerated algorithm which can improve the gradient complexity bound to $\tilde{\mathcal{O}}\left(\sqrt{\kappa_y}\varepsilon ^{-2} \right)$. 
	
	For general nonconvex-concave minimax problem, there are two types of algorithms, i.e., nested-loop algorithms and single-loop algorithms. One intensively studied type is nested-loop algorithms. Rafique et al. \cite{Rafique} propose a proximally guided stochastic mirror descent method (PG-SMD/PGSVRG), which updates $x$ and $y$ simultaneously,
	and provably converges to an approximate stationary point of  $\Phi (\cdot) = \max_{y\in \mathcal{Y}} f(\cdot, y)$. However, only partial convergence results were established for nonconvex-linear minimax problem. Nouriehed et al. \cite{Nouriehed} propose an alternative multi-step framework that finds an $\varepsilon$-first order Nash equilibrium of $f$ with $\tilde{\mathcal{O}}\left( \varepsilon ^{-3.5} \right)$ gradient evaluations. Very recently, Thekumparampil et al.~\cite{Thek} propose a proximal dual implicit accelerated gradient algorithm
	and proved that the algorithm finds an approximate stationary point of $\Phi (\cdot) = \max_{y\in \mathcal{Y}} f(\cdot, y)$ with the rate of $\tilde{\mathcal{O}}\left( \varepsilon ^{-3} \right)$.  Under an equivalent notion of stationary point, Kong and Monteiro~\cite{Kong} propose an accelerated inexact proximal point smoothing method to achieve the same rate, however, at each outer iteration of their algorithm, a perturbed smooth approximation of the inner maximization subproblem  needs to be solved, and the complexity of solving the inner problem has not been considered in \cite{Kong}. Lin et al.~\cite{Lin2020} propose a class of accelerated algorithms for smooth nonconvex-concave minimax problems, which achieves a gradient complexity bound of $\tilde{\mathcal{O}}\left( \varepsilon ^{-2.5} \right)$ in terms of stationary point of $f$. All these nested-loop algorithms either employ multiple gradient ascent steps for $y$'s update to solve the inner subproblem exactly or inexactly, or further do similar acceleration for $x$'s update by adding regularization terms to the inner objective function, and thus are relatively complicated to be implemented.
	
	On the other hand, fewer studies have been directed to single-loop algorithms for nonconvex-concave minimax problems,
	even though these methods are more popular in practice due to their simplicity. One such method is the gradient descent-ascent (GDA) method,
	which performs a gradient descent step on $x$ and a gradient ascent step on $y$ simultaneously at each iteration.
	However, this algorithm fails to converge even for simple bilinear zero-sum games \cite{Letcher}. Many improved GDA algorithms have been proposed in \cite{Chambolle,Daskalakis17,Daskalakis18,Gidel18-1,Gidel18-2,Ho2016}. However, the theoretical understanding of GDA is fairly limited. Very recently, by setting the stepsize of updating $x$ in the order of $\varepsilon^4$, Lin et al.~\cite{Lin2019} proved that the iteration complexity of GDA  to find an $\varepsilon$-stationary point of $\Phi (\cdot) = \max_{y\in \mathcal{Y}} f(\cdot, y)$ is bounded by $\tilde{\mathcal{O}} (\varepsilon^{-6})$ for nonconvex-concave minimax problems when $\mathcal{X}=\mathbb{R}^n$ and $\mathcal{Y}$ is a convex compact set. If both $x$ and $y$ are constrained, the complexity of GDA still remains unknown.  Jin et al. \cite{Jin} propose a GDmax algorithm with iteration complexity $\tilde{\mathcal{O}}\left( \varepsilon ^{-6} \right) $, which corresponds to the number of times the inner maximization problem is solved. Moreover, Lu et al.~\cite{Lu} proposed another single-loop algorithm for nonconvex minimax problems, namely
	the Hybrid Block Successive Approximation (HiBSA) algorithm, which can obtain an $\varepsilon$-stationary point of $f(x,y)$ in $\tilde{\mathcal{O}}\left( \varepsilon ^{-4} \right)$ iterations
	when $\mathcal{X}$ and $\mathcal{Y}$ are  convex compact sets. At each iteration, for updating $x$, one has to solve the subproblem of minimizing a strongly convex function majorizing the
	original function $f$ for a fixed $y$,
	while for updating $y$ one has to solve the subproblem of maximizing the original function $f(x, y)$ plus some regularization terms for a fixed $x$. Under the nonconvex-concave (not strongly concave) setting, both complexity results in \cite{Lu} and \cite{Jin} count the number of times the inner maximization problem is solved, without taking into account
	the complexity of solving the inner maximization problem. Moreover, different stopping rules have been adopted in these existing works, e.g., \cite{Lin2020,Lu,Nouriehed}.
	
	As mentioned earlier, another interesting class of minimax problems is the (strongly) convex-nonconcave setting of \eqref{problem:1}, i.e., $f(x,y)$ is (strongly) convex w.r.t. $x$ and nonconcave w.r.t. $y$.
	However, for any given $x$, to solve the inner maximization subproblem, i.e., $\max_{y\in \mathcal{Y}} f(x, y)$, is already NP-hard.
	Due to this reason, almost all the existing nested-loop algorithms will lose their theoretical guarantees since they need to solve the inner subproblem exactly, or approximately with an error proportional to the accuracy $\varepsilon$. Most existing single-loop algorithms, e.g., HiBSA or GDmax, will also get stuck, since they require the solution of the inner maximization problem.
	One possible alternative approach would be to switch the order of the ``$\min$" and ``$\max$" operators. However, in general, $\min \limits_{x\in \mathcal{X}}\max \limits_{y\in \mathcal{Y}}\ f(x,y) \not= \max \limits_{y\in \mathcal{Y}}	\min \limits_{x\in \mathcal{X}}\ f(x,y)$  when $f(x,y)$ is nonconvex w.r.t. $x$ or nonconcave w.r.t. $y$. The set of stationary points
	for these problems obtained by switching the order of ``$\min$" and ``$\max$" could also be different under some criterions, e.g., $\|\nabla \Phi (\cdot)\|\leq \varepsilon $ with $\Phi (\cdot) = \max_{y\in \mathcal{Y}} f(\cdot, y)$ which is meaningful when $f(\cdot, y)$ is strongly concave with respect to $y$ as defined in \cite{Lin2020}. Whereas the set of stationary points for the above two problems might the same (e.g., in terms of the stationarity of $f$),
	the algorithms applied to these problems may have drastically different trajectories and would converge to quite different solutions.

	It is worth mentioning the most general nonconvex-nonconcave minimax problems, i.e., $f(x,y)$ is nonconvex w.r.t. $x$ and nonconcave w.r.t $y$.
	It should be noted that it is unclear whether a stationary point exists or not for solving these general minimax problems.
	Most recent works aimed at defining a notion of goodness or developing new practical algorithms for reducing oscillations
	and speeding up the convergence of gradient dynamics \cite{Adolphs,Daskalakis18,Heusel,Hsieh2018,Mazumdar}.
	The convergence results for all these algorithms hold either in local region or asymptotically and hence can not imply the global convergence rate. Recently, Flokas et al.~\cite{Flokas} analyze the GDA algorithm for nonconvex-nonconcave Zero-Sum Games, i.e., $\mathcal{X}=\mathbb{R}^n$ and $\mathcal{Y}=\mathbb{R}^m$.
	Lin et al.~\cite{Lin2018} propose a proximal algorithm to solve a special case of the nonconvex-nonconcave minimax problem, where the function $f(x,y)$ satisfies a
	generalized monotone variational inequality condition (see \cite{DangLan12-1}), and show its convergence to stationary points. Nouriehed et al.~\cite{Nouriehed} study a special class of nonconvex-nonconcave minimax problems, in which $f(x,\cdot)$ satisfies the Polyak-{\L}ojasiewic(PL) condition. More recently, Yang et al.~\cite{Yang2020} study the two-sided PL minimax problems and proposed a variance reduced strategy for solving them.

	

	{\bfseries Contributions}. In this paper, we focus on single-loop algorithms for solving nonconvex-concave and convex-nonconcave minimax problems. Our main contributions are as follows.
	
	We propose a simple and unified single-loop Alternating Gradient Projection (AGP) algorithm for solving both  nonconvex-(strongly) concave and (strongly) convex-nonconcave smooth minimax problems. At each iteration, only simple gradient projection steps are employed for updating $x$ and $y$ alternatively.
	
	We analyze the gradient complexity of the proposed unified AGP algorithm under four different settings.
	For the nonconvex-concave setting, we show that an $\varepsilon$-stationary point of $f$ can be obtained in $\mathcal{O}\left( \varepsilon ^{-2} \right)$ (resp. $\mathcal{O}\left( \varepsilon ^{-4} \right)$) iterations for nonconvex-strongly concave (resp. nonconvex concave) minimax problems.
	To the best of our knowledge, these represent the state-of-the-art single loop algorithms under nonconvex-concave setting. Secondly, we show that the gradient complexity to obtain an $\varepsilon$-stationary point of $f$ is $\mathcal{O}\left( \varepsilon ^{-2} \right)$ (resp., $\mathcal{O}\left( \varepsilon ^{-4} \right)$) under the strongly convex-nonconcave (resp.,  convex-nonconcave ) setting. To the best of our knowledge, these
	are the first two theoretically guaranteed convergence results reported in the literature under these settings. Existing single-loop algorithms under both nonconvex-concave and convex-nonconcave settings are summarized in Table \ref{complexity}.

	As shown in Table~\ref{complexity}, AGP matches the best-known $O(\epsilon^{-2})$ complexity for  the nonconvex-strongly concave setting, whereas for the general nonconvex-concave setting, it improves the best-known complexity for single-loop algorithms from $O(\epsilon^{-6})$ to $O(\epsilon^{-4})$. A key step in our development is to construct a suitable potential function which involves function value plus some distance between two adjacent iterates for the nonconvex-strongly concave and strongly convex-nonconcave case, whereas an additional quadratic regularization term of $x$ or $y$ is needed for the
	general nonconvex-concave and convex-nonconcave case (see Lemma 3.4 and Lemma 4.4). To the best of our knowledge, this is the first time that such potential functions have been constructed for solving minimax problems. Moreover, the functional descent results in Lemmas 3.1, 3.3, 4.1, and 4.3 have not been reported before and appear to be novel from our point of view.
	It should be noted that the complexity of AGP does not match the best-known  $O(\epsilon^{-2.5})$ complexity possessed by nested-loop algorithms for the general nonconvex-concave case.
	One possible reason is that only one gradient projection step is employed for solving the inner problem in single loop algorithms, and hence
the error associated with the gradient for the outer problem will be larger than that for nested loop algorithms and such errors will accumulate as the algorithms proceeds.
This also explains why their iteration complexity analysis might be more difficult than that of nested-looped algorithms. It is not evident to us
whether or not this complexity bound obtained for single-loop algorithms can be further improved. Nevertheless, due to its simplicity and the fact that it does not require
the input and fine-tuning of too many algorithmic parameters,
 the proposed AGP algorithm can numerically outperform the state-of-the-art nested loop algorithms as shown in Section 6.

	Furthermore, AGP provides a flexible algorithmic framework that can be easily generalized to solve more complicated minimax problems.
	Based on the basic idea of the AGP algorithm, we propose a block alternating proximal gradient (BAPG) algorithm for solving more general nonsmooth multi-block nonconvex-(strongly) concave and (strongly) convex-nonconcave minimax problems. Each BAPG iteration requires only simple proximal gradient steps to update each block of the multi-block variables alternatively. We prove the gradient complexity of the proposed BAPG algorithm under these four different settings. To the best of our knowledge, existing algorithms (especially the ones with nested loop) can hardly be extended to these more complicated multi-block settings and these complexity results have not been obtained before.
	
	\begin{table}[!ht]
		\centering
		\caption{Single-loop Algorithms for solving \eqref{problem:1} in one block smooth nonconvex-concave and convex-nonconcave settings}\label{complexity}
		\begin{threeparttable}
			\begin{tabular}{l|c|l|c|c}
				\hline
				\multirow{2}{*}{Algorithms} & \multicolumn{2}{c|}{Nonconvex-concave setting} & \multicolumn{2}{c}{Convex-Nonconcave setting}\\ \cline{2-5}
				& strongly concave & general concave & strongly convex & general convex\\ \hline
				GDA\cite{Lin2019} & $\mathcal{O}\left( \varepsilon ^{-2} \right)$  &  $\tilde{\mathcal{O}}\left( \varepsilon ^{-6} \right)$\tnote{1} & Unknown & Unknown\\
				GDmax\cite{Jin} &  $\mathcal{O}\left( \varepsilon ^{-2} \right)$   & $\tilde{\mathcal{O}}\left( \varepsilon ^{-6} \right)$\tnote{1} $^,$\tnote{3}  &  Unknown &  Unknown\\
				HiBSA\cite{Lu} &  $\mathcal{O}\left( \varepsilon ^{-2} \right)$ & $\tilde{\mathcal{O}}\left( \varepsilon ^{-4} \right) $\tnote{2} $^,$\tnote{3} &  Unknown  &  Unknown\\
				{\bfseries AGP} &  $\mathcal{O}\left( \varepsilon ^{-2} \right)$ &  $\mathcal{O}\left( \varepsilon ^{-4} \right)\tnote{2}$ & $\mathcal{O}\left( \varepsilon ^{-2} \right)$ & $\mathcal{O}\left( \varepsilon ^{-4} \right)$\\
				\hline
			\end{tabular}
			\begin{tablenotes}
				\footnotesize
				\item[1] This complexity is to obtain an $\varepsilon$-stationary point of $\Phi (\cdot) = \max_{y\in \mathcal{Y}} f(\cdot, y)$, i.e., $\|\nabla \Phi (\cdot)\|\leq \varepsilon $ when  $\mathcal{X}=\mathbb{R}^n$ and $\mathcal{Y}$ is a convex compact set.
				\item[2] This complexity is to obtain an $\varepsilon$-stationary point of $f$ which is defined as in Section \ref{sec:2}, when $\mathcal{X}$ and $\mathcal{Y}$ are both convex compact sets.
				\item[3] This complexity corresponds to the number of times the inner maximization problem is solved. Thus it does not consider the complexity of solving the inner problem.
			\end{tablenotes}
		\end{threeparttable}
	\end{table}

	The rest of this paper is organized as follows. In Section 2, we propose a unified alternating gradient projection (AGP) algorithm for nonconvex-(strongly) concave and  (strongly) convex-concave minimax problems, and we then analyze the corresponding gradient complexity for four different settings in Section 3 and Section 4. We propose a block alternating proximal gradient (BAPG) algorithm for solving more general nonsmooth multi-block nonconvex-(strongly) concave and (strongly) convex-nonconcave minimax problems, and also establish the corresponding gradient complexity for four different settings in Section 5. We report some numerical results in Section 6 and make some concluding remarks in the last section.

	{\bfseries Notation}. For vectors, we use $\|\cdot\|$ to denote  the $l_2$-norm. For a function $f(x,y):\mathbb{R}^n\times\mathbb{R}^m\rightarrow \mathbb{R}$, we use $\nabla_{x} f(x,y)$ (or $\nabla_{y} f(x, y)$) to denote the partial gradient of $f$ with respect to the first variable (or the second variable) at point $(x, y)$.
	Let $\mathcal{P}_{\mathcal{X}}$ and $\mathcal{P}_{\mathcal{Y}}$ denote projections onto the sets $\mathcal{X}$ and $\mathcal{Y}$. Finally, we use the notation $\mathcal{O} (\cdot)$ to hide only absolute constants which do not depend on any problem parameter, and $\tilde{\mathcal{O}}(\cdot)$ notation to hide only absolute constants and log factors.

A continuously differentiable function $f(\cdot)$ is called $\theta$-strongly convex if there exists a constant $\theta> 0$ such
that for any $x,y \in \mathcal{X}$,
\begin{equation}
  f(y)\ge f(x)+\langle \nabla f(x),y-x \rangle+\frac{\theta}{2}\|y-x\|^2.\label{scon}
\end{equation}
If $-f$ satisfies \eqref{scon}, $f(\cdot)$ is called $\theta$-strongly concave. A pair $(x^*, y^*)$ is a Nash equilibrium (or equivalently a saddle point) of function $f$, if $\forall x\in \mathcal{X}, $ $\forall y\in \mathcal{Y}$,
\begin{equation}
f(x^*, y)\leq f(x^*, y^*)\leq f(x, y^*).\label{saddlepoint}
\end{equation}
A pair $(x^*, y^*)$ is a local Nash equilibrium (or equivalently a local saddle point) of $f$, if there exists $\delta>0$ such that for any $(x,y)$ in $\mathcal{X}\times \mathcal{Y}$ and $\|x-x^*\|\le \delta, \|y-y^*\|\le \delta$, \eqref{saddlepoint} is satisfied. A pair $(x^*, y^*)$ is an $\varepsilon$-first-order Nash equilibrium of function $f$, if $\mathcal{X}(x^*,y^*)\le \varepsilon$ and $\mathcal{Y}(x^*,y^*)\le \varepsilon$, where
\begin{align*}
\mathcal{X}(x^*,y^*)&:=-\min_x \langle \nabla_xf(x^*,y^*),x-x^* \rangle \quad s.t.~ x\in \mathcal{X},\|x-x^*\|\le 1,\\
\mathcal{Y}(x^*,y^*)&:=\max_y \langle \nabla_yf(x^*,y^*),y-y^* \rangle \quad s.t. ~y\in \mathcal{Y},\|y-y^*\|\le 1.
\end{align*}

	\section{An Alternating Gradient Projection Algorithm for \eqref{problem:1}}\label{sec:2}
	In this section, we propose a unified alternating gradient projection (AGP) algorithm that will
	be used later for solving a few different classes of \eqref{problem:1}.
	Each iteration of the proposed AGP algorithm consists of two gradient projection steps for updating both $x$ and $y$.
	Instead of the original function $f(x, y)$, at the $k$-th iteration, AGP uses the gradient of a regularized version of the original function, i.e.,
	\begin{equation}
		f_k (x, y) = f(x, y) + \tfrac{b_{k}}{2} \|x\|^2 -\tfrac{c_{k}}{2}\| y\|^2,\label{sec:2:1}
	\end{equation}
	where $b_{k} \ge 0$ and $c_{k}\ge 0$ are two regularization parameters.
	More specifically, for a given pair $(x_k, y_k) \in \mathcal{X} \times \mathcal{Y}$, AGP
	minimizes a linearized approximation of $f_k(x, y_k)$ plus some regularized terms to update  $x_k$ as follows:
	\begin{align}
		x_{k+1}&={\arg\min}_{x\in \mathcal{X}}\langle  \nabla _x f_k\left( x_k,y_k \right) ,x-x_k \rangle  +\tfrac{\beta_k}{2}\| x-x_k \| ^2\nonumber\\
		&= \mathcal{P}_\mathcal{X} \left( x_k - \tfrac{1}{\beta_k}  \nabla _xf ( x_k,y_k ) - \tfrac{1}{\beta_k} b_{k} x_k \right), \label{sec:2:2}
	\end{align}
	where  $\mathcal{P}_\mathcal{X}$ is the projection operator onto $\mathcal{X}$ and $\beta_k > 0$ denotes
	a stepsize parameter. Similarly, it updates $y_k$ by maximizing a linearized approximation of $f_k(x_{k+1}, y)$ minus some regularized terms, i.e.,
	\begin{align}
		y_{k+1}&={\arg\max}_{y\in \mathcal{Y}}\langle  \nabla _yf_k\left( x_{k+1},y_k \right),y-y_k \rangle  -\tfrac{\gamma_k}{2}\| y-y_k \| ^2\nonumber\\
		&= \mathcal{P}_\mathcal{Y}\left( y_k+\tfrac{1}{\gamma_k} \nabla _y {f}( x_{k+1},y_k ) - \tfrac{1}{\gamma_k} c_k y_k \right),\label{sec:2:3}
	\end{align}
	where  $\mathcal{P}_\mathcal{Y}$ is the projection operator onto  $\mathcal{Y}$ and $\gamma_k >0$ is another stepsize parameter.
	The proposed AGP method is formally stated in Algorithm \ref{alg:2.1}, where sequences $\{\beta_k\}$, $\{b_k\}$, $\{\gamma_k\}$, $\{c_k\}$ and the stopping rule in Step 4 will be specified later in each of the different problem settings to be studied.
	%
	%
	%
	%
	%
	%

	\begin{algorithm}
		\caption{(An Alternating Gradient Projection (AGP) Algorithm)}
		\label{alg:2.1}
		\begin{algorithmic}
			\STATE{\textbf{Step 1}:Input $x_1,y_1,\beta_1,\gamma_1, b_1, c_1$; Set $k=1$.}
			\STATE{\textbf{Step 2}:Calculate $\beta_k$ and $b_k$, and perform the following update for $x_k$:  	
				\qquad \begin{equation}\label{sec:2:4:agp:upx}
					x_{k+1}=\mathcal{P}_\mathcal{X} \left( x_k - \tfrac{1}{\beta_k}  \nabla _xf ( x_k,y_k ) - \tfrac{1}{\beta_k} b_k x_k \right).
			\end{equation}}		
			\STATE{\textbf{Step 3}:Calculate $\gamma_k$ and $c_k$, and perform the following update for $y_k$:
				\qquad	\begin{equation}\label{sec:2:5:agp:upy}
					y_{k+1}=\mathcal{P}_\mathcal{Y}\left( y_k+\tfrac{1}{\gamma_k} \nabla _y {f}( x_{k+1},y_k ) - \tfrac{1}{\gamma_k} c_k y_k \right).
			\end{equation}}
			\STATE{\textbf{Step 4}:If some stationary condition is satisfied, stop; otherwise, set $k=k+1, $ go to Step 2.}
		\end{algorithmic}
	\end{algorithm}
	
	Observe that if we set $b_k=0$ and $c_k=0$, the AGP algorithm is exactly the alternating version of the GDA algorithm, which is rather natural and has been widely used by practitioners for solving minimax problems, e.g., in generative adversarial networks. 	However, to the best of our knowledge, even for the alternating GDA algorithm, the convergence for solving \eqref{problem:1} has never been established before in the literature. Moreover, if $b_k$ or $c_k$ is not equal to $0$, AGP algorithm is completely new. It turns out that
	$b_k$ or $c_k$ plays a very crucial role to guarantee the convergence of AGP when solving general nonconvex-concave or convex-nonconcave minimax problems.

	Before we prove the iteration complexity of AGP algorithm for solving \eqref{problem:1}, we define the stationarity gap as the termination criterion as follows.
	
	\begin{defi}\label{sec:2:def1:G}
		At each iteration of Algorithm \ref{alg:2.1}, the stationarity gap for problem \eqref{problem:1} w.r.t. $f(x,y)$ is defined as:
		\begin{equation*}
			\nabla G\left( x_k,y_k \right) :=\left[\begin{array}{c}
				{\beta_k\left(x_k-\operatorname{\mathcal{P}_\mathcal{X}}(x_k-\frac{1}{\beta_k} \nabla_{x} f(x_k ,y_k))\right)} \\
				
				{\gamma_k\left(y_k-\operatorname{\mathcal{P}_\mathcal{Y}}(y_k+\frac{1}{\gamma_k} \nabla_{y} f(x_k ,y_k))\right)}
			\end{array}\right].
		\end{equation*}
    	We denote $\nabla {G}_k:= \nabla {G}(x_k, y_k)$,
		$(\nabla G_k)_{x}:= \beta_k( x_k-\operatorname{\mathcal{P}_\mathcal{X}}( x_k-\frac{1}{\beta_k}\nabla _{x}f( x_k,y_k)))$,
		and $(\nabla G_k)_{y}:=\gamma_k( y_k-\operatorname{\mathcal{P}_\mathcal{Y}}( y_k+ \frac{1}{\gamma_k} \nabla _{y}f( x_k,y_k )) )$.
	\end{defi}
	
	Note that the widely used metrics for convex-concave minimax problems such as the min-max value or the distance to the set of optimal solutions of the min-max problem are not applicable in nonconvex cases.
	For the latter cases we call $(x_k,y_k)$ an $\varepsilon$-stationary point of $f(x,y)$ if $ \|\nabla G(x_k,y_k)\|\leq \varepsilon$ with $\nabla G(x_k,y_k)$ being defined as in Definition \ref{sec:2:def1:G}, which has been widely used as the optimality measure, e.g., \cite{Lin2020,Lu}. In the absence of constraints, $ \| \nabla G(x_k,y_k)\| \leq \varepsilon$ reduces to the standard condition $ \|\nabla_{x} f(x_k ,y_k)\|\leq \varepsilon $ and $\|\nabla_{y} f(x_k ,y_k)\|\leq \varepsilon$ for unconstrained problems. The vector  $\nabla G(x_k,y_k)$ also refers to gradient mapping at $(x_k,y_k)$, see \cite{Nes2013} for the details. 
	
	\begin{defi}\label{sec:2:def2:tildeG}
		At each iteration of Algorithm \ref{alg:2.1}, the stationarity gap for problem \eqref{problem:1} w.r.t. $f_{k}(x,y)$ is defined as:
		\begin{equation*}
			\nabla \tilde{G}\left( x_k,y_k \right) :=\left[\begin{array}{c}
				{\beta_k\left(x_k-\operatorname{\mathcal{P}_\mathcal{X}}(x_k-\frac{1}{\beta_k} \nabla_{x}f_{k}(x_k ,y_k))\right)} \\
				
				{\gamma_k\left(y_k-\operatorname{\mathcal{P}_\mathcal{Y}}(y_k+\frac{1}{\gamma_k} \nabla_{y} f_{k}(x_k ,y_k))\right)}
			\end{array}\right].
		\end{equation*}
		We denote $\nabla {\tilde{G}}_k:=\nabla {\tilde{G}}(x_k, y_k)$,
		$(\nabla \tilde{G}_k)_{x}:= \beta_k( x_k-\operatorname{\mathcal{P}_\mathcal{X}}( x_k- \frac{1}{\beta_k} \nabla _{x}f_{k}\left( x_k,y_k \right)))$,
		 and $(\nabla \tilde{G}_k)_{y}:=\gamma_k( y_k-\operatorname{\mathcal{P}_\mathcal{Y}}( y_k+ \frac{1}{\gamma_k} \nabla _{y}f_{k}\left( x_k,y_k \right)) )$.
	\end{defi}
	
	We also need to make the following assumption about the smoothness of $f(x,y)$.
	
	\begin{ass}\label{sec:2:ass1:lip}
		$f\left(x,y\right)$ has Lipschitz continuous gradients, i.e., there exist positive scalars $L_{22}$, $L_{12}$, $L_{11}$, $L_{21}$ such that for any $x,\bar{x}\in \mathcal{X}$, $y, \bar{y}\in \mathcal{Y}$,
		\begin{align*}
			\| \nabla_{x} f\left(x, y\right)-\nabla_{x} f\left(\bar{x}, y\right)\| &\leq L_{11}\|x-\bar{x}\|,\\
			\|\nabla_{x} f\left(x, y\right)-\nabla_{x} f\left(x, \bar{y}\right)\| &\leq L_{21}\|y-\bar{y}\|,\\
			\|\nabla_{y} f\left(x, \bar{y}\right)-\nabla_{y} f\left(x, y\right)\| &\leq L_{22}\|\bar{y}-y\|,\\
			\|\nabla_{y} f\left(x, y\right)-\nabla_{y} f\left(\bar{x}, y\right)\| &\leq L_{12}\|x-\bar{x}\|.
		\end{align*}
		We denote $L:=\max\{L_{11}, L_{22}, L_{12}, L_{21}\}$.
	\end{ass}

Under this assumption, we can first prove the following lemma for estimating bounds on changes in the function value when $x_k$ or $y_k$ is updated at each iteration in Algorithm \ref{alg:2.1}.
\begin{lem}\label{lem:2.1}
Suppose that Assumption \ref{sec:2:ass1:lip} holds. Let $\{(x_k,y_k)\}$ be a sequence generated by Algorithm \ref{alg:2.1}. Then we have
\begin{align}
  f\left( x_{k+1},y_k \right) -f\left( x_k,y_k \right) &\le-\left( \beta_k -\tfrac{L_{11}}{2} \right) \lVert x_{k+1}-x_k \rVert ^2-\frac{b_k}{2}\left[\|x_{k+1}\|^2-\|x_k\|^2-\|x_{k+1}-x_k\|^2\right], \label{lem:2.1:1}\\
  f(x_{k+1},y_{k+1})-f( x_{k+1},y_{k})
			&\ge\left( \gamma_k -\tfrac{L_{22}}{2} \right) \lVert y_{k+1}-y_k \rVert ^2+\frac{c_k}{2}\left[\|y_{k+1}\|^2-\|y_k\|^2-\|y_{k+1}-y_k\|^2\right].\label{lem:2.1:2}
\end{align}
\end{lem}

\begin{proof}
  Firstly, by the optimality condition for \eqref{sec:2:4:agp:upx}, we have
		\begin{equation}\label{lem:2.1:3}
			\langle  \nabla _xf\left( x_k,y_k \right) +\beta_k \left( x_{k+1}-x_k \right)+b_kx_k ,x_k-x_{k+1} \rangle  \ge 0.
		\end{equation}
		By Assumption \ref{sec:2:ass1:lip}, the gradient of $f$ is Lipschitz continuous, implying that
		\begin{equation}
			f\left( x_{k+1},y_k \right) -f\left( x_k,y_k \right)\leq \langle  \nabla _xf\left( x_k,y_k \right) ,x_{k+1}-x_k \rangle  +\tfrac{L_{11}}{2}\lVert x_{k+1}-x_k \rVert ^2. \label{lem:2.1:4}
		\end{equation}
		Adding \eqref{lem:2.1:3} and \eqref{lem:2.1:4}, we can easily show that
		\begin{align*}
			f\left( x_{k+1},y_k \right) -f\left( x_k,y_k \right) &\le -\left( \beta_k -\tfrac{L_{11}}{2} \right) \lVert x_{k+1}-x_k \rVert ^2-b_k\langle x_k, x_{k+1}-x_k\rangle\nonumber\\
            &= -\left( \beta_k -\tfrac{L_{11}}{2} \right) \lVert x_{k+1}-x_k \rVert ^2-\frac{b_k}{2}\left[\|x_{k+1}\|^2-\|x_k\|^2-\|x_{k+1}-x_k\|^2\right].
		\end{align*}
This completes the proof of \eqref{lem:2.1:1}. By the optimality condition for \eqref{sec:2:5:agp:upy}, we have
		\begin{equation}\label{lem:2.1:6}
			\langle  \nabla _{y}f( x_{k+1},y_{k})-\gamma_k( y_{k+1}-y_k)-c_ky_k ,y_k-y_{k+1} \rangle  \leq 0.
		\end{equation}
		By  Assumption \ref{sec:2:ass1:lip}, the gradient of $f$ is Lipschitz continuous,  implying that
		\begin{align}\label{lem:2.1:7}
			f(x_{k+1},y_{k+1})-f( x_{k+1},y_{k})
			\ge  \langle \nabla _{y}f(x_{k+1},y_{k}),y_{k+1}
    -y_{k} \rangle -\frac{L_{22}}{2}\| y_{k+1}-y_{k} \|^2.
		\end{align}
			The result \eqref{lem:2.1:2} then follows by adding \eqref{lem:2.1:6} and \eqref{lem:2.1:7}.

\end{proof}

		In the following two sections, we will establish the iteration complexity of the AGP algorithm under four different problem settings. Although there are some different technical details under different problem settings, the main process used in these proofs is similar. Firstly, we estimate a bound on the change of function values between two adjacent iterates, shown as in Lemmas 3.1, 3.3, 4.1 and 4.3 respectively.  Then, we estimate an upper bound on the weighted distance between two adjacent pairs of iterates by constructing a suitable potential function according to different properties of the objective function, shown as in Lemmas 3.2, 3.4, 4.2 and 4.4 respectively. Finally, by using those upper bounds, we prove the complexity of the algorithm through some careful parameter selection, shown as in Theorems 3.1, 3.2, 4.1 and 4.2.
	\section{Complexity Analysis for Nonconvex-Concave Minimax Problems}\label{sec:3}	
	
	%

	\subsection{Nonconvex-Strongly Concave Setting}
	In this subsection, we analyze the iteration complexity of Algorithm \ref{alg:2.1} for solving nonconvex-strongly concave minimax optimization problems, i.e., $f(x,y)$ is nonconvex w.r.t. $x$ for any fixed $y\in \mathcal{Y}$,
	and $\mu$-strongly concave  w.r.t. $y$ for any given $x\in \mathcal{X}$. Under this setting,  $\forall k\ge 1$, we set
	\begin{equation}\label{subsec:3.1:1}
		\beta_k=\eta, \gamma_k=\tfrac{1}{\rho}, \ \mbox{and} \ b_k=c_k=0
	\end{equation}
	in Algorithm \ref{alg:2.1}, and simplify the update for $x_k$ and $y_k$ as follows:
	\begin{align}
		x_{k+1}&=\mathcal{P}_\mathcal{X} \left( x_k - \tfrac{1}{\eta}  \nabla _xf ( x_k,y_k )  \right),\label{subsec:3.1:2:x}\\
		y_{k+1}&=\mathcal{P}_\mathcal{Y}\left( y_k+\rho \nabla _y {f}( x_{k+1},y_k ) \right),\label{subsec:3.1:3:y}
	\end{align}
	which is exactly the alternating version of GDA algorithm. Our goal in the remaining part of this subsection is to establish the iteration complexity of Algorithm \ref{alg:2.1} under the nonconvex-strongly concave setting.

	We now establish an important recursion for the AGP algorithm.
	\begin{lem}\label{lem:3.2:blostr:y}
		Suppose that Assumption \ref{sec:2:ass1:lip} holds. Let $\{\left(x_k,y_k\right)\}$ be a sequence generated by Algorithm \ref{alg:2.1}
		with parameter settings in \eqref{subsec:3.1:1}. If $\eta > L_{11}$, then we have
		\begin{align}\label{lem:3.2:1}
			&f( x_{k+1},y_{k+1}) -f\left( x_k,y_k \right)  \nonumber \\
			\leq& -\left(\tfrac{\eta}{2}-\tfrac{L_{12}^{2}\rho}{2}\right)\| x_{k+1}-x_{k} \| ^2-\left(\tfrac{\mu}{2}-\tfrac{1}{\rho}\right)\| y_{k+1}-y_k \|^2\nonumber \\
			&-\left(\mu-\tfrac{1}{2\rho}-\tfrac{\rho L_{22}^2}{2}\right)\|y_k-y_{k-1} \|^2.
		\end{align}
	\end{lem}
	
	\begin{proof}
		The optimality condition for $y_k$ in \eqref{subsec:3.1:3:y} implies that $\forall y\in \mathcal{Y}$ and $\forall k\geq 1$,
		\begin{equation}\label{lem:3.2:2}
			\langle \nabla _yf(x_{k+1},y_k)-\frac{1}{\rho}(y_{k+1}-y_k),y-y_{k+1} \rangle \le 0.
		\end{equation}
		By choosing $y=y_k$ in \eqref{lem:3.2:2}, we have
		\begin{equation}\label{lem:3.2:3}
			\langle \nabla _yf(x_{k+1},y_k)-\frac{1}{\rho}(y_{k+1}-y_k),y_k-y_{k+1} \rangle \le 0.
		\end{equation}
		On the other hand, by replacing $k$ with $k-1$ and choosing $y=y_{k+1}$ in \eqref{lem:3.2:2}, we obtain
		\begin{equation}\label{lem:3.2:4}
			\langle \nabla _yf(x_{k},y_{k-1})-\frac{1}{\rho}(y_k-y_{k-1}),y_{k+1}-y_k \rangle \le 0,
		\end{equation}
		which, in view of the fact that $f\left( x,y  \right)$ is $\mu $-strongly concave w.r.t. $y$ for any given $x\in \mathcal{X}$,
		then implies that
		\begin{align}\label{lem:3.2:5}
			&f( x_{k+1},y_{k+1}) -f( x_{k+1},y_k ) \nonumber\\
			\leq &  \langle  \nabla _yf( x_{k+1},y_k ),y_{k+1}-y_k\rangle  -\frac{\mu}{2}\| y_{k+1}-y_k \|^2\  \nonumber\\
			\leq &\langle  \nabla _yf( x_{k+1},y_k ) -\nabla _yf( x_k,y_{k-1}) ,y_{k+1}-y_k \rangle   \nonumber\\
			& +\frac{1}{\rho}\langle  y_k-y_{k-1},y_{k+1}-y_k\rangle  -\frac{\mu}{2}\| y_{k+1}-y_k \|^2.
		\end{align}
		Denoting $v_{k+1}:=\left(y_{k+1}-y_k\right)-\left(y_k-y_{k-1}\right)$, we can write the first inner product term in the r.h.s.
		of \eqref{lem:3.2:5} as
		\begin{align}\label{lem:3.2:6}
			&\langle  \nabla _y f ( x_{k+1},y_k ) -\nabla _y f ( x_k,y_{k-1}) ,y_{k+1}-y_k \rangle  \nonumber\\
			=&\langle  \nabla _y f ( x_{k+1},y_k ) -\nabla _y f \left( x_k,y_{k} \right) ,y_{k+1}-y_k \rangle  \nonumber\\
			& + \langle  \nabla _y f \left( x_k,y_{k} \right) -\nabla _y f ( x_k,y_{k-1}) ,v_{k+1} \rangle  \nonumber\\
			&+\langle  \nabla _y f \left( x_k,y_{k}\right) -\nabla _y f ( x_k,y_{k-1}) ,y_k-y_{k-1} \rangle .
		\end{align}
		Next, we estimate the three terms in the right hand side of \eqref{lem:3.2:6} respectively.
		By Assumption \ref{sec:2:ass1:lip} and the Cauchy-Schwarz inequality, we can bound the first two terms according to
		\begin{align}\label{lem:3.2:7}
			&\langle  \nabla _yf( x_{k+1},y_k ) -\nabla _yf\left( x_k,y_{k} \right) ,y_{k+1}-y_k \rangle  \nonumber\\
			\leq &\frac{L_{12}^{2}\rho}{2}\| x_{k+1}-x_k \|^2+\frac{1}{2\rho}\| y_{k+1}-y_k \|^2,
		\end{align}
		and
		\begin{align}\label{lem:3.2:8}
			&\langle  \nabla _yf\left(  x_k,y_{k} \right) -\nabla _yf( x_k,y_{k-1}) ,v_{k+1} \rangle  \leq \frac{\rho L_{22}^{2}}{2}\|y_k-y_{k-1} \|^2+\frac{1}{2\rho}\|v_{k+1} \|^2.
		\end{align}
		For the third term, by $\mu $-strongly-concavity  of $f$ with respect to $y$,
		\begin{align}\label{lem:3.2:9}
			\quad \langle  \nabla _yf\left( x_k,y_{k} \right) -\nabla _yf( x_k,y_{k-1}) ,y_k-y_{k-1} \rangle  \le -\mu \| y_k-y_{k-1} \|^2 .
		\end{align}
		Moreover, it can be easily checked that
		\begin{align}\label{lem:3.2:10:3points}
			&\langle  y_k-y_{k-1},y_{k+1}-y_k \rangle =\frac{1}{2}\| y_k-y_{k-1} \|^2+\frac{1}{2}\| y_{k+1}-y_k \|^2-\frac{1}{2}\| v_{k+1} \|^2.
		\end{align}
		Plugging \eqref{lem:3.2:6}-\eqref{lem:3.2:10:3points} into \eqref{lem:3.2:5} and rearranging the terms, we conclude that
		\begin{align}\label{lem:3.2:11}
			&f( x_{k+1},y_{k+1}) -f( x_{k+1},y_{k})  \nonumber \\
			\leq& \frac{L_{12}^{2}\rho}{2}\| x_{k+1}-x_k \|^2-(\mu-\frac{1}{2\rho}-\frac{\rho L_{22}^2}{2})\| y_k-y_{k-1}\|^2\nonumber\\
			&-(\frac{\mu}{2}-\frac{1}{\rho})\|y_{k+1}-y_k \|^2.
		\end{align}
By setting $\beta_k=\eta$, $b_k=0$ in \eqref{lem:2.1:1} of Lemma \ref{lem:2.1}  and the assumption $\eta > L_{11}$, we have
		\begin{equation}\label{lem:3.1:1}
			f(x_{k+1},y_{k})-f(x_{k},y_k)\le -\frac{\eta}{2}\| x_{k+1}-x_k \|^2.
		\end{equation}
The proof is completed by combining \eqref{lem:3.2:11} with \eqref{lem:3.1:1}.
	\end{proof}
	
	One may want to take the telescoping sum of \eqref{lem:3.2:1} in order to provide
	a bound on $\sum_k (\| x_{k+1}-x_k \|^2+\| y_{k+1}-y_k \|^2)$.
	However, since $(1/\rho + \rho L_{22}^2)/2 \ge L_{22} \ge \mu$, the coefficient
	of the third term in the r.h.s. of  \eqref{lem:3.2:1} is always positive.
	As a result, we need to further refine this relation as shown below.	
	
	\begin{lem}\label{lem:3.3:blostr:F}
		Suppose that Assumption \ref{sec:2:ass1:lip} holds. Let $\{\left(x_k,y_k\right)\}$ be a sequence generated by Algorithm \ref{alg:2.1}
		with parameter settings in \eqref{subsec:3.1:1}.
		Also let us denote
		\begin{align*}
			f_{k+1} := f( x_{k+1},y_{k+1}), \ \
			S_{k+1} :=\frac{2}{\rho^2\mu}\|y_{k+1}-y_k \|^2, \\
			F_{k+1} :=f_{k+1}+S_{k+1}- (\mu+\frac{7}{2\rho}-\frac{\rho L_{22}^2}{2}-\frac{2L_{22}^2}{\mu})\|y_{k+1}-y_k\|^2.
		\end{align*}
		If $\eta>L_{11}$,
		then
		\begin{align}\label{lem:3.3:2}
			F_{k+1}-F_k\leq& -\left(\tfrac{\eta}{2}-\tfrac{\rho L_{12}^2}{2}-\tfrac{2L_{12}^2}{\rho \mu^2}\right)\| x_{k+1}-x_k \|^2 \nonumber\\
			& -\left(\tfrac{3\mu-\rho L_{22}^2}{2}+\tfrac{\mu-4\rho L_{22}^2}{2\rho\mu}\right)\| y_{k+1}-y_k \|^2.
		\end{align}
	\end{lem}
	
	\begin{proof}
		First by \eqref{lem:3.2:3} and \eqref{lem:3.2:4}, we have
		\begin{align}\label{lem:3.3:3}
			\frac{1}{\rho}\langle v_{k+1},y_{k+1} -y_k \rangle
			\le \langle \nabla _yf( x_{k+1},y_k ) -\nabla _yf\left( x_{k},y_{k-1} \right),y_{k+1}-y_k \rangle,
		\end{align}
		which together with \eqref{lem:3.2:6} then imply that
		\begin{align}\label{lem:3.3:4}
			\frac{1}{\rho}\langle v_{k+1},y_{k+1} -y_k  \rangle  \leq&\langle  \nabla _yf( x_{k+1},y_k ) -\nabla _yf\left( x_{k},y_k \right) ,y_{k+1}-y_k \rangle  \nonumber\\
			& + \langle  \nabla _yf\left( x_{k},y_k \right) -\nabla _yf\left( x_{k},y_{k-1} \right) ,v_{k+1} \rangle  \nonumber\\
			&+\langle  \nabla _yf\left( x_{k},y_k \right) -\nabla _yf\left( x_{k},y_{k-1} \right) ,y_{k}-y_{k-1} \rangle .
		\end{align}
		Similar to \eqref{lem:3.2:7}, we can easily see that
		\begin{align}\label{lem:3.3:5}
			\left\langle  \nabla _yf( x_{k+1},y_k ) -\nabla _yf\left( x_{k},y_k  \right) ,y_{k+1}-y_k \right\rangle
			\leq &\frac{L_{12}^{2}}{2 \mu}\|x_{k+1}-x_{k} \|^2+\frac{\mu}{2}\|y_{k+1}-y_k \|^2.
		\end{align}
		By plugging \eqref{lem:3.2:8},\eqref{lem:3.2:9},\eqref{lem:3.3:5} into \eqref{lem:3.3:4},
		and using the identity
		$\frac{1}{\rho}\langle v_{k+1}, y_{k+1}-y_k \rangle  = \frac{1}{2\rho}\|y_{k+1}-y_k\|^2+\frac{1}{2\rho}\| v_{k+1}\| ^2-\frac{1}{2\rho}\|y_{k}-y_{k-1}\|^2$,
		we conclude that
		\begin{align}\label{lem:3.3:6}
			&\frac{1}{2\rho}\|y_{k+1}-y_k\|^2+\frac{1}{2\rho}\| v_{k+1}\|^2-\frac{1}{2\rho}\|y_{k}-y_{k-1}\|^2\nonumber\\
			\leq& \frac{L_{12}^{2}}{2\mu}\| x_{k+1}-x_{k}\|^2+\frac{\mu}{2}\| y_{k+1}-y_k \|^2+\frac{\rho L_{22}^2}{2}\|y_{k}-y_{k-1} \|^2\nonumber\\
			&+\frac{1}{2\rho}\| v_{k+1}\|^2 -\mu\|y_{k}-y_{k-1} \|^2.
		\end{align}
		Rearranging the terms of \eqref{lem:3.3:6}, we have
		\begin{align}\label{lem:3.3:7}
			&\frac{1}{2\rho}\|y_{k+1}-y_k\|^2-\frac{1}{2\rho}\|y_{k}-y_{k-1}\|^2\nonumber\\
			\leq& \frac{L_{12}^{2}}{2\mu}\| x_{k+1}-x_{k}\|^2+\frac{\mu}{2}\| y_{k+1}-y_k \|^2-\left(\mu-\frac{\rho L_{22}^2}{2}\right)\|y_{k}-y_{k-1} \|^2.
		\end{align}
		Multiplying $\frac{4}{\rho\mu}$ on both sides of \eqref{lem:3.3:7} and using the definition of $S_{k+1}$, we obtain
		\begin{align*}
			S_{k+1}-S_k \leq & \frac{2L_{12}^{2}}{\mu^2\rho}\| x_{k+1}-x_{k}\|^2+\frac{2}{\rho}\| y_{k+1}-y_k \|^2-\left(\frac{4}{\rho}-\frac{2 L_{22}^2}{\mu}\right)\|y_{k}-y_{k-1} \|^2.
		\end{align*}
		It then follows from \eqref{lem:3.2:1} in Lemma \ref{lem:3.2:blostr:y} and the definition of $F_k$ that
		\begin{align*}
			F_{k+1}-F_k\leq& -\left(\frac{\eta}{2}-\frac{\rho L_{12}^2}{2}-\frac{2L_{12}^2}{\rho \mu^2}\right)\| x_{k+1}-x_{k}\|^2 \nonumber\\
			& -\left(\frac{3\mu-\rho L_{22}^2}{2}+\frac{\mu-4\rho L_{22}^2}{2\rho\mu}\right)\| y_{k+1}-y_k \|^2.
		\end{align*}
	\end{proof}
	
	We are now ready to establish the iteration complexity for the AGP algorithm in the nonconvex-strongly concave setting.
	In particular,
	letting $\nabla G\left( x_k,y_k \right)$ be defined as  in Definition \ref{sec:2:def1:G} and $\varepsilon>0$ be
	a given target accuracy, we provide
	a bound on $T(\varepsilon)$, the first iteration index to achieve an $\varepsilon$-stationary point, i.e., $ \| \nabla G(x_k,y_k)\| \leq \varepsilon$, which is equivalent to
	\begin{equation}\label{subsec:3.1:4:defT}
		T(\varepsilon):=\min\{k \mid \|\nabla G(x_k,y_k)\|\leq \varepsilon \}.
	\end{equation}

	\begin{thm}\label{thm:3.1:blockstr}
		Suppose that Assumption \ref{sec:2:ass1:lip} holds. Let $\{\left(x_k,y_k\right)\}$ be a sequence generated by Algorithm \ref{alg:2.1}
		with parameter settings in \eqref{subsec:3.1:1}.
		If the relations $\eta>L_{11}, \eta>L_{12}^2\rho+\frac{4L_{12}^2}{\rho \mu^2}, \ \mbox{and} \ \rho\leq \frac{\mu}{4L_{22}^2}$ are satisfied, then
		it holds that $$ T\left( \varepsilon \right) \le \frac{F_1-\underline{F}}{d_1\varepsilon ^2},$$
		where $	d_1 :=\min \left\{\frac{\eta}{2}-\frac{\rho L_{12}^2}{2}-\frac{2L_{12}^2}{\rho \mu^2},\frac{3\mu-\rho L_{22}^2}{2}+\frac{\mu-4\rho L_{22}^2}{2\rho\mu} \right\}/\max \left\{ \eta^2+2L_{12}^2,\frac{2}{\rho^2} \right\}$ and  $ 			\underline{F}:=\underline{f}-(\mu+\frac{7}{2\rho}-\frac{\rho L_{22}^2}{2}-\frac{2L_{22}^2}{\mu})\sigma_y^2$ with $\underline{f} := \min_{(x,y) \in \mathcal{X} \times \mathcal{Y}} f(x,y)$ and $\sigma _y:= \max\{\|y_1-y_2\| \mid y_1,y_2 \in  \mathcal{Y}\}$.
	\end{thm}
	
	\begin{proof}
		It follows immediately from \eqref{subsec:3.1:1} and \eqref{subsec:3.1:2:x} that
		\begin{align}\label{thm:3.1:1:gx}
			\|(\nabla G_k)_{x} \|= \eta\| x_{k+1}-x_k \|.
		\end{align}
		On the other hand,  by \eqref{subsec:3.1:3:y} and the triangle inequality, we obtain that
		\begin{align}\label{thm:3.1:2:gy}
			&\| (\nabla G_k)_y \|\nonumber\\
			\le&\frac{1}{\rho}\|\operatorname{ \mathcal{P}_\mathcal{Y}}( y_k+\rho \nabla _{y}f( x_{k+1} ,y_{k}) )-\operatorname{ \mathcal{P}_\mathcal{Y}}\left( y_k+\rho \nabla _{y}f\left( x_{k} ,y_{k} \right) \right)\|+\frac{1}{\rho} \| y_{k+1}-y_k \|\nonumber\\
			\le&L_{12}\|x_{k+1}-x_k\|+\frac{1}{\rho} \| y_{k+1}-y_k \|.
		\end{align}
		By combining \eqref{thm:3.1:1:gx} and \eqref{thm:3.1:2:gy}, and using the Cauchy-Schwarz inequality, we obtain
		\begin{align}\label{thm:3.1:3:g}
			\| \nabla G_k \|^2
			\leq& \left(\eta^2+2L_{12}^2\right)\| x_{k+1}-x_k \|^2+\frac{2}{\rho ^2}\|y_{k+1}-y_k \|^2.
		\end{align}
		Observing that $d_1 > 0$. Multiplying both sides of \eqref{thm:3.1:3:g} by $d_1$, and using \eqref{lem:3.3:2} in Lemma \ref{lem:3.3:blostr:F},
		we have
		\begin{align}\label{thm:3.1:4:d1}
			d_1\|\nabla G_k \|^2\le F_k -F_{k+1}.
		\end{align}
		Summing up the above inequalities from $k=1$ to $k=T(\varepsilon)$, we obtain
		\begin{align}\label{thm:3.1:5}
			\sum_{k=1}^{T\left( \varepsilon \right)}{d_1\|\nabla G_k \| ^2}\le F_1 -F_{T\left( \varepsilon \right)+1}.
		\end{align}
		Note that by the definition of $F_{k+1}$ in Lemma \ref{lem:3.3:blostr:F}, we have
		\begin{align*}
			F_{T\left( \varepsilon \right)+1}&=
			f_{T\left( \varepsilon \right)+1}+S_{T\left( \varepsilon \right)+1}-(\mu+\frac{7}{2\rho}-\frac{\rho L_{22}^2}{2}-\frac{2L_{22}^2}{\mu})\|y_{T\left( \varepsilon \right)+1}-y_{T\left( \varepsilon \right)} \|^2\\
			&\ge \underline{f}-(\mu+\frac{7}{2\rho}-\frac{\rho L_{22}^2}{2}-\frac{2L_{22}^2}{\mu})\sigma_y^2 = \underline{F},
		\end{align*}
		where the inequality follows from the definitions of $\underline{f}$ and  $\sigma _y$, and the
		facts that $S_k\ge 0$ ($\forall k\ge 1$) and $\mu+\frac{7}{2\rho}-\frac{\rho L_{22}^2}{2}-\frac{2L_{22}^2}{\mu} \ge 0 $ due to
		the selection of $\rho$.
		We then conclude from \eqref{thm:3.1:5} that $\sum_{k=1}^{T\left( \varepsilon \right)}{d_1\| \nabla G_k \|^2}\le F_1 -F_{T\left( \varepsilon \right)+1}\le F_1-\underline{F}$ which, in view of the definition of $T(\varepsilon)$, implies that $\varepsilon ^2\le (F_1-\underline{F})/(T( \varepsilon )\cdot d_1)$
		or equivalently, $T\left( \varepsilon \right) \le (F_1-\underline{F})/(d_1\varepsilon ^2)$.
	\end{proof}

	A few remarks are in place for the results obtained in Theorem \ref{thm:3.1:blockstr}.
	First, in view of	Theorem \ref{thm:3.1:blockstr}, the gradient complexity of Algorithm \ref{alg:2.1} to obtain a stationary
		point that satisfies $\|\nabla G(x_k,y_k)\|\leq \varepsilon$ is given by $\mathcal{O} (L^2\varepsilon^{-2})$ under
		the nonconvex-strongly concave setting.
	Second, under this setting, very few single-loop algorithms have been investigated
	although a few other existing algorithms can achieve similar complexity bound.
	 In particular, it seems that even
	the complexity bound for the GDA algorithm remains unknown under this setting when $\mathcal{X}$ and $\mathcal{Y}$ are both convex compact sets. Compared to the state-of-the-art algorithm in \cite{Lin2020}, we improve the complexity for the nonconvex-strongly concave setting by a logarithmic factor, since we do not need to solve the inner problem at each iteration.
	Third, as mentioned earlier, Algorithm \ref{alg:2.1} is a single-loop alternating gradient projection method with constant stepsizes,
	which is very easy to implement in practice.
	
	\subsection{Complexity Analysis for General Nonconvex-Concave Setting}
	In this subsection, we analyze the iteration complexity of Algorithm \ref{alg:2.1} for solving \eqref{problem:1} under the general nonconvex-concave setting. Under this setting, $\forall k\ge 1$, we set
	\begin{equation}\label{subsec:3.2:1}
		b_k=0,\quad \beta_k=\bar{\eta}+\bar{\beta}_k,\quad \gamma_k=\frac{1}{\bar{\rho}},
	\end{equation}
	where $\bar{\eta}>0$, $\bar{\rho}>0$ are two constants, and $\bar{\beta}_k$ are stepsize parameters to be defined later. We need to make the following assumption on the parameters $c_k$.
	
	\begin{ass}\label{subsec:3.2:ass1:ck}
		$\{c_k\}$ is a nonnegative monotonically decreasing sequence.
	\end{ass}

By Assumption \ref{sec:2:ass1:lip} and $\nabla _yf_{k-1}\left( x,y \right) =\nabla _yf\left( x,y \right) -c_{k-1} y$, we have
	\begin{align}\label{subsec:3.2:2:newLy}	
		&\|\nabla_{y} f_{k-1}\left(x_{k},y_{k}\right)-\nabla_{y}f_{k-1}(x_{k},y_{k-1})\|\nonumber\\
		=&\|\nabla_{y}f\left(x_{k},y_{k}\right)-\nabla_{y}f\left(x_{k},y_{k-1}\right)-c_{k-1}\left(y_{k}-y_{k-1}\right)\|\nonumber\\
		\leq& \left(L_{22}+c_{k-1}\right)\| y_{k}-y_{k-1}\|.
	\end{align}
	Denoting $L_{22}^{'}=L_{22}+c_1$, by Assumption \ref{subsec:3.2:ass1:ck} and  \eqref{subsec:3.2:2:newLy}, we have
	$$\|\nabla_{y} f_{k-1}\left(x_{k},y_{k}\right)-\nabla_{y}f_{k-1}\left(x_{k},y_{k-1}\right)\|\leq L_{22}^{'}\| y_{k}-y_{k-1}\|.$$
	It then follows from the above inequality and the strong concavity of $f_{k-1}\left(x_{k},y \right)$ w.r.t. $y$ (Theorem 2.1.12 in \cite{Nestrov}) that
	\begin{align}\label{subsec:3.2:3:key}
		&\quad\langle  \nabla _yf_{k-1}\left( x_{k},y_{k} \right) -\nabla _yf_{k-1}\left( x_{k},y_{k-1} \right) ,y_{k}-y_{k-1} \rangle \nonumber\\
		&\leq-\frac{1}{L_{22}^{'}+c_{k-1}}\|\nabla _yf_{k-1}\left(x_{k},y_{k} \right)-\nabla _yf_{k-1}\left( x_{k},y_{k-1} \right) \|^2-\frac{c_{k-1}L_{22}^{'}}{L_{22}^{'}+c_{k-1}}\|y_{k}-y_{k-1} \|^2.
	\end{align}	
	This is a key inequality that we will use to establish some important recursions for the AGP method under the nonconvex-concave setting in the following two results. This is also one of the key differences between the proof for the nonconvex-strongly concave setting and the one for the nonconvex-concave setting.
	\begin{lem}\label{lem:3.5}
		Suppose that Assumptions \ref{sec:2:ass1:lip} and \ref{subsec:3.2:ass1:ck} hold. Let $\{(x_k,y_k)\}$ be a sequence generated by Algorithm \ref{alg:2.1} with parameter settings in \eqref{subsec:3.2:1}. If $\forall k,\bar{\beta}_k > L_{11}$ and $\bar{\rho} \le \frac{2}{L_{22}^{'}+c_1}$, then
		\begin{align}\label{lem:3.5:1:conclu}
			&f(x_{k+1},y_{k+1})-f(x_k,y_k)\nonumber \\
			\leq& -\left( \bar{\eta} +\frac{\bar{\beta}_k}{2}-\frac{\bar{\rho} L_{12}^{2}}{2} \right) \| x_{k+1}-x_k \| ^2+\frac{1}{\bar{\rho}} \|y_{k+1}-y_k \| ^2+\frac{1}{2\bar{\rho}}\| y_k-y_{k-1} \|^2\nonumber\\
			&+\frac{c_{k-1}}{2}(\| y_{k+1}\| ^2-\| y_k\|^2).
		\end{align}
	\end{lem}
	\begin{proof}
		The optimality condition for $y_k$ in \eqref{sec:2:5:agp:upy} implies that, $\forall y\in \mathcal{Y}$, $\forall k\geq 1$,
		\begin{equation}\label{lem:3.5:2:opty:k+1}
			\langle \nabla _yf_k(x_{k+1},y_k)-\frac{1}{\bar{\rho}}(y_{k+1}-y_k) ,y-y_{k+1} \rangle \le 0.
		\end{equation}
		By choosing $y=y_k$ in \eqref{lem:3.5:2:opty:k+1}, we have
		\begin{equation}\label{lem:3.5:3:newopty:k+1}
			\langle \nabla _yf_k(x_{k+1},y_k)-\frac{1}{\bar{\rho}}(y_{k+1}-y_k),y_k-y_{k+1} \rangle \le 0.
		\end{equation}
		On the other hand, by replacing $k$ with $k-1$ and choosing $y=y_{k+1}$ in \eqref{lem:3.5:2:opty:k+1}, we obtain
		\begin{equation}\label{lem:3.5:4:opty:k}
			\langle \nabla _yf_{k-1}(x_{k},y_{k-1})-\frac{1}{\bar{\rho}}(y_k-y_{k-1}) ,y_{k+1}-y_k \rangle \le 0.
		\end{equation}
		The concavity of  $f_k(x_{k+1},y)$ w.r.t. $y$ together with  \eqref{lem:3.5:4:opty:k} then imply that
		\begin{align}\label{lem:3.5:5}
			&f_k(x_{k+1},y_{k+1})-f_k(x_{k+1},y_k)\nonumber\\
			\le& \langle \nabla _yf_k(x_{k+1},y_k)-\nabla _xf_{k-1}(x_{k},y_{k-1}),y_{k+1}-y_k \rangle+\frac{1}{\bar{\rho}}\langle y_k-y_{k-1},y_{k+1}-y_k \rangle\nonumber\\
			=&\langle \nabla _yf_k(x_{k+1},y_k)-\nabla _yf_{k-1}(x_k,y_{k}),y_{k+1}-y_k \rangle+\langle \nabla _yf_{k-1}(x_k,y_{k})-\nabla _yf_{k-1}(x_{k},y_{k-1}),v_{k+1} \rangle\nonumber\\
			&+\langle \nabla _yf_{k-1}(x_k,y_{k})-\nabla _yf_{k-1}(x_{k},y_{k-1}),y_k-y_{k-1} \rangle+\frac{1}{\bar{\rho}}\langle y_k-y_{k-1},y_{k+1}-y_k \rangle,
		\end{align}
		where $v_{k+1}=y_{k+1}-y_k-(y_k-y_{k-1})$. We now provide bounds on the inner product terms of \eqref{lem:3.5:5}. Firstly, by the definition of $f_k(x_{k+1},y_k)$ and $f_{k-1}(x_{k},y_k)$, Assumptions \ref{sec:2:ass1:lip} and \ref{subsec:3.2:ass1:ck}, and the Cauchy-Schwarz inequality, we have
		\begin{align}\label{lem:3.5:6}
			&\langle  \nabla _yf_k( x_{k+1},y_k) -\nabla _yf_{k-1}( x_{k},y_{k}) ,y_{k+1}-y_k \rangle  \nonumber\\
			=&\langle \nabla_{y} f(x_{k+1},y_k)-\nabla_{y} f(x_k,y_{k}), y_{k+1}-y_k\rangle-\left(c_k-c_{k-1}\right)\langle y_k, y_{k+1}-y_k\rangle \nonumber\\
			\le&\frac{\bar{\rho} L_{12}^{2}}{2}\| x_{k+1}-x_{k} \|^2+\frac{1}{2\bar{\rho}}\| y_{k+1}-y_k \|^2-\frac{c_k-c_{k-1}}{2}(\| y_{k+1}\| ^2-\| y_k\|^2)\nonumber\\
			&+\frac{c_k-c_{k-1}}{2}\| y_{k+1}-y_k \|^2\nonumber\\
			\leq& \frac{\bar{\rho} L_{12}^{2}}{2}\| x_{k+1}-x_{k} \|^2+\frac{1}{2\bar{\rho}}\| y_{k+1}-y_k \|^2-\frac{c_k-c_{k-1}}{2}(\| y_{k+1}\| ^2-\| y_k\|^2).	
		\end{align}
		Secondly, by the Cauchy-Schwarz inequality,
		\begin{align}\label{lem:3.5:7}
			\langle  \nabla _yf_{k-1}( x_k,y_{k}) -\nabla _yf_{k-1}( x_{k},y_{k-1}) ,v_{k+1} \rangle\leq& \frac{\bar{\rho}}{2}\|\nabla _yf_{k-1}( x_k,y_{k}) -\nabla _yf_{k-1}( x_{k},y_{k-1}) \|^2 \nonumber\\ &+\frac{1}{2\bar{\rho}}\|v_{k+1} \|^2.
		\end{align}
		Thirdly, it follows from  \eqref{subsec:3.2:3:key} that
		\begin{align}\label{lem:3.5:8}
			&\langle \nabla _yf_{k-1}\left( x_k,y_{k} \right) -\nabla _yf_{k-1}\left( x_{k},y_{k-1} \right) ,y_k-y_{k-1} \rangle\nonumber\\
			\le&-\frac{1}{L_{22}^{'}+c_{k-1}}\| \nabla _yf_{k-1}\left( x_k,y_{k} \right)-\nabla _yf_{k-1}\left( x_{k},y_{k-1} \right) \|^2-\frac{c_{k-1}L_{22}^{'}}{L_{22}^{'}+c_{k-1}}\|y_k-y_{k-1} \|^2\nonumber\\
			\leq& -\frac{1}{L_{22}^{'}+c_{k-1}}\| \nabla _yf_{k-1}\left( x_k,y_{k} \right)-\nabla _yf_{k-1}\left( x_{k},y_{k-1} \right)\|^2.
		\end{align}
		Also observe that
		\begin{equation}\label{lem:3.5:9}
			\frac{1}{\bar{\rho}}\langle  y_{k+1}-y_k,y_k-y_{k-1} \rangle =\frac{1}{2\bar{\rho}}\|y_{k+1}-y_k \|^2+\frac{1}{2\bar{\rho}}\|y_k-y_{k-1} \|^2-\frac{1}{2\bar{\rho}}\|v_{k+1} \|^2.
		\end{equation}
		Plugging \eqref{lem:3.5:6}-\eqref{lem:3.5:9} into \eqref{lem:3.5:5}, and using the definition of $f_k(x_{k+1},y_{k+1})$ and $f_k(x_{k+1},y_k)$ and the assumption $\frac{\bar{\rho}}{2} \le \frac{1}{L_{22}^{'}+c_1}$, we obtain
		\begin{align}\label{lem:3.5:10}
			f(x_{k+1},y_{k+1})-f(x_{k+1},y_k)
			&\leq \frac{\bar{\rho} L_{12}^{2}}{2}\| x_{k+1}-x_{k}\|^2+\frac{1}{\bar{\rho}} \|y_{k+1}-y_k \| ^2\nonumber\\
			&\quad+\frac{1}{2\bar{\rho}}\| y_k-y_{k-1} \|^2+\frac{c_{k-1}}{2}(\| y_{k+1}\| ^2-\| y_k\|^2).
		\end{align}
By setting $\beta_k=\bar{\eta}+\bar{\beta}_k$, $b_k=0$ in \eqref{lem:2.1:1} of Lemma \ref{lem:2.1}  and the assumption $\bar{\beta}_k > L_{11}$, we have
		\begin{equation}\label{lem:3.4:1conclu}
			f(x_{k+1},y_{k})-f(x_{k},y_k)\le -\left( \bar{\eta} +\frac{\bar{\beta}_k}{2} \right) \| x_{k+1}-x_k \|^2.
		\end{equation}
		The result in \eqref{lem:3.5:1:conclu} then follows by adding \eqref{lem:3.5:10} and \eqref{lem:3.4:1conclu}.
	\end{proof}
	
	It turns out that from Lemma \ref{lem:3.5} we can not obtain an upper bound on the
	positively weighted sum of $\| x_{k+1}-x_k\|^2$ and $\|y_{k+1}-y_k \|^2$
	to provide an upper bound for $\|\nabla G_k\|$. We need to further refine this result in~\eqref{lem:3.5:1:conclu} to overcome this difficulty. In particular, we obtain
	below a new inequality as in \eqref{lem:3.6:5} to further investigate the relation between $\| x_{k+1}-x_k\|^2$ and $\|y_{k+1}-y_k \|^2$. 
	Then by using this new inequality, we construct a new potential function  as shown in the following important result for Algorithm~\ref{alg:2.1}. Note that in Lemma \ref{lem:3.3:blostr:F} we have constructed a potential function which involves the function value plus some distance between two adjacent iterates for the nonconvex-strongly concave cases, whereas in the following lemma an additional quadratic regularization term, i.e., $\|y_{k+1}\|^2$, is added to construct another potential function for the general nonconvex-concave case.

	\begin{lem}\label{lem:3.6}
		Suppose that Assumptions \ref{sec:2:ass1:lip} and \ref{subsec:3.2:ass1:ck} hold. Let $\{(x_k,y_k)\}$ be a sequence generated by Algorithm \ref{alg:2.1} with parameter settings in \eqref{subsec:3.2:1}. Also let us denote
		\begin{align*}
			\mathcal{S}_{k+1}&:=\frac{8}{\bar{\rho}^2c_{k+1}}\| y_{k+1}-y_k \|^2+\frac{8}{\bar{\rho}}\left(1-\frac{c_k}{c_{k+1}}\right)\|y_{k+1}\|^2,\\
			\mathcal{F}_{k+1}&:=f(x_{k+1},y_{k+1})+\mathcal{S}_{k+1}- \frac{15}{2\bar{\rho}}\| y_{k+1}-y_k \|^2-\frac{c_k}{2}\|y_{k+1}\|^2.
		\end{align*}
		If
		\begin{equation}\label{lem:3.6:1:par}
        \bar{\beta}_k >L_{11},~\frac{1}{c_{k+1}}-\frac{1}{c_k}\leq \frac{\bar{\rho}}{5},~\bar{\rho} \leq \frac{2}{L_{22}^{'}+c_1},
		\end{equation}
		then $\forall k \geq 2$,
		\begin{align}\label{lem:3.6:2:conclu}
			\mathcal{F}_{k+1}-\mathcal{F}_{k}
			&\leq -\left( \bar{\eta}+\frac{\bar{\beta}_k}{2}-\frac{\bar{\rho} L_{12}^2}{2}-\frac{16L_{12}^2}{\bar{\rho} c_{k}^2}\right)\| x_{k+1}-x_k \| ^2+\frac{c_{k-1}-c_{k}}{2}\|y_{k+1}\|^2\nonumber\\
			& \quad-\frac{9}{10\bar{\rho}}\| y_{k+1}-y_k \|^2+\frac{8}{\bar{\rho}}\left(\frac{c_{k-1}}{c_k}-\frac{c_k}{c_{k+1}} \right)\|y_{k+1}\|^2.
		\end{align}
	\end{lem}
	
	\begin{proof}
		By \eqref{lem:3.5:3:newopty:k+1} and \eqref{lem:3.5:4:opty:k}, we have
		\begin{align}\label{lem:3.6:3}
			&\frac{1}{\bar{\rho}}\langle v_{k+1},y_{k+1} -y_k \rangle \le \langle \nabla _yf_k( x_{k+1},y_k ) -\nabla _yf_{k-1}\left( x_{k},y_{k-1} \right),y_{k+1}-y_k \rangle.
		\end{align}
	Similar to \eqref{lem:3.5:5} in Lemma \ref{lem:3.5}, \eqref{lem:3.6:3} can be rewritten as
		\begin{align*}
			&\frac{1}{\bar{\rho}}\langle v_{k+1}, y_{k+1} -y_k \rangle\nonumber\\
			\leq&\langle \nabla _yf_k(x_{k+1},y_k)-\nabla _yf_{k-1}(x_k,y_{k}),y_{k+1}-y_k \rangle+\langle \nabla _yf_{k-1}(x_k,y_{k})-\nabla _yf_{k-1}(x_{k},y_{k-1}),v_{k+1} \rangle\nonumber\\
			&+\langle \nabla _yf_{k-1}(x_k,y_{k})-\nabla _yf_{k-1}(x_{k},y_{k-1}),y_k-y_{k-1} \rangle.
		\end{align*}
		Using an argument similar to the proof of \eqref{lem:3.5:6}-\eqref{lem:3.5:9},
		the relation in \eqref{subsec:3.2:3:key}, and the Cauchy-Schwarz inequality, we conclude from the above inequality that
		\begin{align}\label{lem:3.6:4}	&\frac{1}{2\bar{\rho}}\|y_{k+1}-y_k\|^2+\frac{1}{2\bar{\rho}}\|v_{k+1}\|^2-\frac{1}{2\bar{\rho}}\|y_k-y_{k-1}\|^2\nonumber\\
			\leq& \frac{L_{12}^{2}}{2a_k}\|x_{k+1}-x_{k} \|^2+\frac{a_k}{2}\| y_{k+1}-y_k \|^2-\frac{c_k-c_{k-1}}{2}(\| y_{k+1}\|^2-\| y_k\| ^2)\nonumber\\
			&+\frac{\bar{\rho}}{2}\|  \nabla _yf_{k-1}(x_k,y_{k})-\nabla _yf_{k-1}(x_{k},y_{k-1})\|^2-\frac{1}{L_{22}^{'}+c_{k-1}}\|\nabla _yf_{k-1}(x_k,y_{k})-\nabla _yf_{k-1}(x_{k},y_{k-1}) \|^2\nonumber\\
			&+\frac{1}{2\bar{\rho}}\|v_{k+1} \|^2-\frac{c_{k-1}L_{22}^{'}}{L_{22}^{'}+c_{k-1}}\| y_k-y_{k-1} \|^2+ \frac{c_k-c_{k-1}}{2}\|y_{k+1}-y_k\|^2,
		\end{align}
		for any $a_k>0$. Observing by $c_1\leq L_{22}^{'}$, and Assumption \ref{subsec:3.2:ass1:ck}, we have
		\begin{equation*}
			-\frac{c_{k-1}L_{22}^{'}}{c_{k-1}+L_{22}^{'}}\leq-\frac{c_{k-1}L_{22}^{'}}{2L_{22}^{'}}
			=-\frac{c_{k-1}}{2}\leq-\frac{c_k}{2}.
		\end{equation*}
		Combining $\bar{\rho} \le \frac{2}{L_{22}^{'}+c_1}$ and rearranging the terms in \eqref{lem:3.6:4}, we obtain
		\begin{align*}
			& \frac{1}{2\bar{\rho}}\| y_{k+1}-y_k \|^2+\frac{c_k-c_{k-1}}{2}\|y_{k+1}\|^2 \nonumber\\
			\leq&\frac{1}{2\bar{\rho}}\|y_k-y_{k-1}\|^2+\frac{c_k-c_{k-1}}{2}\|y_k\|^2+\frac{L_{12}^{2}}{2a_k}\|x_{k+1}-x_{k} \|^2+\frac{a_k}{2}\| y_{k+1}-y_k \|^2\nonumber\\
			&-\frac{c_k}{2}\|y_k-y_{k-1} \|^2.
		\end{align*}
		By multiplying $\frac{16}{\bar{\rho} c_k}$ on both sides of the above inequality, we then obtain
		\begin{align}\label{lem:3.6:5}
			&\frac{8}{\bar{\rho}^2c_k}\| y_{k+1}-y_k \|^2+\frac{8}{\bar{\rho}}\left(1-\frac{c_{k-1}}{c_k}\right)\|y_{k+1}\|^2\nonumber\\
			\leq& \frac{8}{\bar{\rho}^2c_k}\|y_k-y_{k-1}\|^2 +\frac{8}{\bar{\rho}}(1-\frac{c_{k-1}}{c_k}
			)\|y_k\|^2+\frac{8L_{12}^2}{\bar{\rho} c_ka_k}\|x_{k+1}-x_{k} \|^2\nonumber\\
			& +\frac{8a_k}{\bar{\rho} c_k}\| y_{k+1}-y_k \|^2 -\frac{8}{\bar{\rho}}\|y_k-y_{k-1}\|^2.
		\end{align}
		Setting $a_k=\frac{c_k}{2}$ in the above inequality, and using the definition of $\mathcal{S}_{k+1}$ and \eqref{lem:3.6:1:par}, we have
		\begin{align}\label{lem:3.6:6}
			\mathcal{S}_{k+1}-\mathcal{S}_{k}
			&\leq \frac{8}{\bar{\rho}}\left(\frac{c_{k-1}}{c_k}-\frac{c_k}{c_{k+1}} \right)\|y_{k+1}\|^2+\frac{16L_{12}^2}{\bar{\rho} c_k^2}\|x_{k+1}-x_{k} \|^2+\frac{28}{5\bar{\rho}}\| y_{k+1}-y_k \|^2\nonumber\\
			&\quad-\frac{8}{\bar{\rho}}\|y_k-y_{k-1}\|^2.
		\end{align}
		Combining \eqref{lem:3.6:6} and \eqref{lem:3.5:1:conclu} in Lemma \ref{lem:3.5}, and using the definition of $\mathcal{F}_{k+1}$, we conclude
		\begin{align*}
			\mathcal{F}_{k+1}-\mathcal{F}_{k}
			&\leq -\left( \bar{\eta}+\frac{\bar{\beta}_k}{2}-\frac{\bar{\rho} L_{12}^2}{2}-\frac{16L_{12}^2}{\bar{\rho} c_{k}^2}\right)\| x_{k+1}-x_k \| ^2+\frac{c_{k-1}-c_{k}}{2}\|y_{k+1}\|^2\\
			&\quad-\frac{9}{10\bar{\rho}}\| y_{k+1}-y_k \|^2+\frac{8}{\bar{\rho}}\left(\frac{c_{k-1}}{c_k}-\frac{c_k}{c_{k+1}} \right)\|y_{k+1}\|^2.
		\end{align*}
	\end{proof}
	
	%
	
	We are now ready to establish the iteration complexity for the AGP algorithm to achieve an $\varepsilon$-stationary point in the general nonconvex-concave setting.

	\begin{thm}\label{thm3.2}
		Suppose that Assumption \ref{sec:2:ass1:lip} holds. Let $\{(x_k,y_k)\}$ be a sequence generated by Algorithm \ref{alg:2.1} with parameter settings in \eqref{subsec:3.2:1}.
	If $\bar{\rho} \le \frac{1}{10L_{22}}$, $\tau>\max\{\frac{19^2(L_{11}+2\bar{\eta}-\bar{\rho}L_{12}^2)}{20^2\cdot16\bar{\rho}L_{12}^2}, 2\}$, $c_k=\frac{19}{20\bar{\rho} k^{\text{1/}4}}, k\geq 1$, then for any given $\varepsilon>0$,
		$$T(\varepsilon)
		\leq \max \left(\left( \frac{2\cdot80^2\bar{\rho} (\tau -2) L_{12}^2d_2d_3 }{19^2\varepsilon ^2}+1 \right) ^2, \frac{19^4\hat{\sigma}_y^4}{10^4\bar{\rho} ^4\varepsilon ^4}\right),$$
		where $d_2=\mathcal{F}_2-\underline{\mathcal{F}}+ (8\cdot 2^{1/4} + \frac{19}{40} )\frac{\hat{\sigma}_y^2}{\bar{\rho}}$, $\bar{d}_1=\frac{8\tau^2}{(\tau-2)^2} + \frac{19^4\left(2\left(\bar{\rho} L_{12}^2-\bar{\eta}\right)^2 + 2L_{12}^2\right)}{64\cdot20^4 \bar{\rho}^2 (\tau-2)^2 L_{12}^4}$, 	$d_3=\max\{\bar{d}_1,\frac{19^2}{1440 \sqrt{2}(\tau -2 )\bar{\rho}^2 L_{12}^2}\}$, $\hat{\sigma}_y:= \max\{\|y\| \mid y \in  \mathcal{Y}\}$,
		$\underline{\mathcal{F}} := \underline{f}-\left(2^{13/4}+15+\frac{19}{40}\right)\frac{\hat{\sigma}_y^2}{\bar{\rho}}$ with  $\underline{f} := \min_{(x,y) \in \mathcal{X} \times \mathcal{Y}} f(x,y)$.
	\end{thm}
	
	\begin{proof}
		By $\bar{\rho}\le \frac{1}{10L_{22}}$ and $c_k=\frac{19}{20\bar{\rho} k^{\text{1/}4}},\forall k\geq 1$, let us denote $\bar{\beta}_k=\bar{\rho} L_{12}^2+\frac{16\tau L_{12}^2}{\bar{\rho} c_{k}^2}-2\bar{\eta}$, $\alpha_k=\frac{8(\tau - 2)L_{12}^2}{\bar{\rho} c_{k}^2}$, we can easily see that the relations in \eqref{lem:3.6:1:par} are satisfied. It follows from the selection of $\bar{\beta}_k$ and $\alpha_k$ that
		$$\bar{\eta}+\frac{\bar{\beta}_k}{2}-\frac{\bar{\rho} L_{12}^2}{2}-\frac{16L_{12}^2}{\bar{\rho} c_{k}^2}=\alpha_k.$$
		This observation, in view of Lemma \ref{lem:3.6}, then immediately implies that
		\begin{align}\label{thm3.2:1}
			\alpha_k\|x_{k+1}-x_k \|^2+\frac{9}{10\bar{\rho}}\| y_{k+1}-y_k \|^2
			\leq& \mathcal{F}_{k}-\mathcal{F}_{k+1}+\frac{8}{\bar{\rho}}\left(\frac{c_{k-1}}{c_k} -\frac{c_k}{c_{k+1}} \right)\|y_{k+1}\|^2\nonumber\\
			&+\frac{c_{k-1}-c_k}{2}\|y_{k+1}\|^2.
		\end{align}
		We can easily check from the definition of $f_{k} (x, y)$ that  $$\|\nabla {G}_k\| -\|\nabla \tilde{G}_k \| \leq c_{k}\| y_k \|.$$
		By replacing $f$ with $f_{k}$, $\eta$ with $\bar{\beta}_k+\bar{\eta}$, similar to \eqref{thm:3.1:1:gx} and \eqref{thm:3.1:2:gy}, we immediately obtain that
		%
		\begin{align}\label{thm3.2:2:gx}
			\|(\nabla \tilde{G}_k)_{x} \|
			= \left(\bar{\beta}_k+\bar{\eta}\right)\| x_{k+1}-x_k \|,
		\end{align}
		and
		\begin{align}\label{thm3.2:3:gy}
			\| (\nabla \tilde{G}_k)_y \|
			\le&\frac{1}{\bar{\rho}} \| y_{k+1}-y_k \|+L_{12}\|x_{k+1}-x_k\|.
		\end{align}
		Combining \eqref{thm3.2:2:gx} and \eqref{thm3.2:3:gy}, and using the Cauchy-Schwarz inequality, we have
		\begin{align}\label{thm3.2:4:g}
			\| \nabla \tilde{G}_k \|^2
			\leq& \left(\left(\bar{\beta}_k+\bar{\eta}\right)^2+2L_{12}^2\right)\| x_{k+1}-x_k \|^2+\frac{2}{\bar{\rho} ^2}\|y_{k+1}-y_k \|^2.
		\end{align}
		Since both $\alpha_k$ and $\bar{\beta}_k$ are in the same order when $k$ becomes large enough, it then follows from the definition of $\bar{d}_1$ that $\forall k\ge 1$,
		\begin{align}\label{thm3.2:5:d1}
			\frac{\left(\bar{\beta}_k+\bar{\eta}\right)^2+2L_{12}^{2}}{\alpha _k^2} &=  \frac{\left( \bar{\rho} L_{12}^2+\frac{16\tau L_{12}^2}{\bar{\rho} c_{k}^2}-\bar{\eta} \right) ^2+2L_{12}^{2}}{\alpha _k^2}\nonumber\\
			&\leq \frac{2\left(\frac{16\tau L_{12}^2}{\bar{\rho} c_{k}^2}\right)^2+ 2\left(\bar{\rho} L_{12}^2-\bar{\eta}\right)^2 + 2L_{12}^2}{\alpha _k^2}\nonumber\\
			&= \frac{8\tau^2}{(\tau-2)^2} + \frac{19^4\left(2\left(\bar{\rho} L_{12}^2-\bar{\eta}\right)^2 + 2L_{12}^2\right)}{64\cdot20^4 \bar{\rho}^2 (\tau-2)^2 L_{12}^4k} \leq \bar{d}_1.
		\end{align}
		Combining the previous two inequalities in \eqref{thm3.2:5:d1} and \eqref{thm3.2:4:g}, we obtain
		\begin{equation}\label{thm3.2:6:newG}
			\| \nabla \tilde{G}_k \|^2
			\leq \bar{d}_1(\alpha_k)^2\| x_{k+1}-x_k \|^2+ \frac{2}{\bar{\rho} ^2}\|y_{k+1}-y_k \|^2.
		\end{equation}
	Denote $d_k^{(2)}=\frac{1}{\max \left\{ \bar{d}_1\alpha_k,\frac{20}{9\bar{\rho}} \right\}}$.
		By multiplying $d_k^{(2)}$ on the both sides of \eqref{thm3.2:6:newG}, and using \eqref{thm3.2:1}, we have
		\begin{align}\label{thm3.2:7}
			d_k^{(2)}\| \nabla \tilde{G}_k \|^2
			\leq &\mathcal{F}_{k}-\mathcal{F}_{k+1}+\frac{8}{\bar{\rho}}\left(\frac{c_{k-1}}{c_k} -\frac{c_k}{c_{k+1}} \right)\|y_{k+1}\|^2+\frac{c_{k-1}-c_k}{2}\|y_{k+1}\|^2,
		\end{align}
		where the last inequality follows since $d_{k}^{(2)}$ is a decreasing sequence. Denoting
		$$\tilde{T}(\varepsilon):=\min\{k \mid \|\nabla \tilde{G}(x_k,y_k) \|\leq \frac{\varepsilon}{2}, k\geq 2\},$$
		Summing both sides of \eqref{thm3.2:7} from $k=2$ to $k=\tilde{T}(\varepsilon)$, we then obtain
		\begin{align}\label{thm3.2:8}
&\sum_{k=2}^{\tilde{T}(\varepsilon)}{d_k^{(2)}\lVert \nabla \tilde{G}_k \rVert ^2}\nonumber\\
\leq& \mathcal{F}_2-\mathcal{F}_{\tilde{T}(\varepsilon)}+\frac{8}{\bar{\rho}}\left( \frac{c_1}{c_2}-\frac{c_{\tilde{T}(\varepsilon)}}{c_{\tilde{T}(\varepsilon)+1} } \right)\hat{\sigma}_y^2+\frac{c_1-c_{\tilde{T}(\varepsilon)}}{2}\hat{\sigma}_y^2\nonumber\\
			\leq&\mathcal{F}_2-\mathcal{F}_{\tilde{T}(\varepsilon)}+\frac{8c_1}{\bar{\rho} c_2}\hat{\sigma}_y^2+\frac{c_1}{2}\hat{\sigma}_y^2\nonumber\\
			=& \mathcal{F}_2-\mathcal{F}_{\tilde{T}(\varepsilon)}+ \left(8\cdot 2^{1/4} + \frac{19}{40} \right)\frac{\hat{\sigma}_y^2}{\bar{\rho}}.
		\end{align}

Note that by the definition of $\mathcal{F}_{k+1}$ in Lemma \ref{lem:3.6}, we have
\begin{align*}
\mathcal{F}_{\tilde{T}(\varepsilon)}&\ge\underline{f}-\frac{8\hat{\sigma}_y^2}{\bar{\rho}}\left(1+\frac{1}{\tilde{T}(\varepsilon)}\right)^{1/4}
-\frac{15\hat{\sigma}_y^2}{\bar{\rho}}-\frac{c_1\hat{\sigma}_y^2}{2}\\
&=\underline{f}-\left(2^{13/4}+15+\frac{19}{40}\right)\frac{\hat{\sigma}_y^2}{\bar{\rho}}=\underline{\mathcal{F}},
\end{align*}
where $\underline{f} := \min_{(x,y) \in \mathcal{X} \times \mathcal{Y}} f(x,y)$. We then conclude from \eqref{thm3.2:8} that

\begin{align}\label{thm3.2:8n}
\sum_{k=2}^{\tilde{T}(\varepsilon)}{d_k^{(2)}\lVert \nabla \tilde{G}_k \rVert ^2}\le \mathcal{F}_2-\underline{\mathcal{F}}+ \left(8\cdot 2^{1/4} + \frac{19}{40} \right)\frac{\hat{\sigma}_y^2}{\bar{\rho}}=d_2.
\end{align}
		We can see
		from the selection of $d_3$ that $d_3=\max\{\bar{d}_1,\frac{20}{9\bar{\rho}\alpha_2}\}$. Observe that $\alpha_k$ is an increasing sequence, when $k\geq 2$, we have that $d_3\ge\max\{\bar{d}_1,\frac{20}{9\bar{\rho}\alpha_k}\}$, which implies that $d_k^{(2)}\geq\frac{1}{d_3\alpha_k}$, by multiplying $d_3$ on the both sides of \eqref{thm3.2:8n}, and combining the definition of $d_2$, we have $	\sum_{k=2}^{\tilde{T}(\varepsilon)} \frac{1}{\alpha_k} \|\nabla \tilde{G}_k \|^2\leq d_2d_3$,
		which, by the definition of $\tilde{T}(\varepsilon)$, implies that
		\begin{equation}\label{thm3.2:9}
			\frac{\varepsilon^2}{4} \leq \frac{d_2d_3}{\tsum_{k=2}^{\tilde{T}(\varepsilon)}{\frac{1}{\alpha_k}}}.
		\end{equation}
		Note that when $c_k=\frac{19}{20\bar{\rho} k^{1/4}}$, $\alpha_k = \frac{2\cdot40^2\bar{\rho} (\tau -2) L_{12}^2 \sqrt{k}}{19^2}$. By using the fact $\sum_{k=2}^{\tilde{T}(\varepsilon)}1/\sqrt{k}\geq \sqrt{\tilde{T}(\varepsilon)}-1$ and \eqref{thm3.2:9}, we conclude $\frac{\varepsilon^2}{4} \leq \frac{2\cdot40^2\bar{\rho}(\tau -2) L_{12}^2d_2d_3}{19^2\left(\sqrt{\tilde{T}(\varepsilon)}-1\right)}$ or equivalently,
		\begin{equation*}
			\tilde{T}(\varepsilon) \le \left( \frac{2\cdot80^2\bar{\rho} (\tau -2) L_{12}^2d_2d_3 }{19^2\varepsilon ^2}+1 \right) ^2.
		\end{equation*}
		On the other hand, if $k \ge \frac{19^4\hat{\sigma}_y^4}{10^4\bar{\rho} ^4\varepsilon ^4}$, then $c_{k}=\frac{19}{20\bar{\rho} k^{1/4}} \leq \frac{\varepsilon}{2\hat{\sigma}_y}$. This inequality together with the definition of $\hat{\sigma}_y$ then imply that $c_{k}\| y_k \|\leq \frac{\varepsilon}{2}$. Therefore, there exists a
		\begin{align*}
			T(\varepsilon) &\leq \max (	\tilde{T}(\varepsilon), \frac{19^4\hat{\sigma} _y^4}{10^4\bar{\rho}^4\varepsilon ^4})\\
			&\leq \max \left(	\left( \frac{2\cdot80^2\bar{\rho} (\tau -2) L_{12}^2d_2d_3 }{19^2\varepsilon ^2}+1 \right) ^2, \frac{19^4\hat{\sigma}_y^4}{10^4\bar{\rho} ^4\varepsilon ^4}\right),
		\end{align*}
		such that
		$\|\nabla G_k\|\leq \| \nabla \tilde{G}_k \|+ c_{k}\| y_k \| \leq \frac{\varepsilon}{2}+\frac{\varepsilon}{2} =\varepsilon$.
	\end{proof}

	According to Theorem \ref{thm3.2}, we can show that, by specifying $\gamma_k$ and $c_k$ as in the order of $k^{1/2}$ and $1/k^{1/4}$, the gradient complexity of Algorithm \ref{alg:2.1} to obtain an $\varepsilon$-stationarity point of $f$ in the nonconvex-concave setting can be bounded by $\mathcal{O} (L^4\varepsilon^{-4})$. In this setting the stepsize for updating $x_k$ is in the order of $k^{-1/2}$, while	the one for updating $y_k$ is a constant at iteration $k$.

	It is worth mentioning two closely related works to ours. One is the HiBSA algorithm proposed in \cite{Lu} for nonconvex-concave problems.
	It seems that our iteration complexity is slightly better than that of the HiBSA algorithm by  a logarithmic factor.
	More importantly, the two subproblems for updating $x_k$ and $y_k$ in AGP are much easier than those in the HiBSA algorithm, as
	the latter method needs to solve $\max_{y\in \mathcal{Y}} f(x, y)$ at each iteration to update $y_k$ under the general nonconvex-concave (not necessarily strongly concave) setting.
	The other related method is the GDA algorithm. By setting the stepsize to update $x_k$ in the order of $\varepsilon^4$, Lin et al.~\cite{Lin2019} proved that the iteration complexity of GDA to return an $\varepsilon$-stationary point of $\Phi (\cdot) = \max_{y\in \mathcal{Y}} f(\cdot, y)$ is bounded by $\tilde{\mathcal{O}} (\varepsilon^{-6})$ for nonconvex-concave minimax problems when $\mathcal{X}=\mathbb{R}^n$ and $\mathcal{Y}$ is a convex compact set. In this paper, both $\mathcal{X}$ and $\mathcal{Y}$ are convex compact sets, hence the outer problem is a constrained minimization problem instead of an unconstrained one. Under this setting, the stopping rule, i.e.,  $\|\nabla \Phi(\cdot)\|\leq \varepsilon$,  used in \cite{Lin2019} for unconstrained minimization setting needs to be reconsidered. Moreover, the step sizes to update $x_k$ in AGP is a decreasing sequence in the order of $k^{-1/2}$, which is significantly larger than  ${\cal O} (\varepsilon^4)$ stepsize allowed in GDA, and thus may result in better practical performance.
	
	Note that for solving unconstrained bilinear minimax problem, the alternating GDA algorithm with any fixed step size will cause recurrence~\cite{Bailey2020}.
	In the classic convex optimization literature, such a divergence issue was usually handled by incorporating averaging, smoothing, or direct acceleration techniques (see, e.g., Sections 3.5-3.8 and Sections 4.3-4.5 of~\cite{Lan2020book}). However, the proposed AGP algorithm for nonconvex minimax problems is not equivalent to the alternating GDA algorithm,
	since a regularized version of the original function is incorporated, and a variable stepsize policy has been used for the nonconvex-concave setting. 
Hence, the results of our paper do not contradict with existing ones.


\section{Complexity Analysis for Convex-Nonconcave  Minimax Problems}\label{sec4}

In this section, we establish the convergence of AGP algorithm for the cases where $f$ is convex w.r.t. $x$, but possibly nonconcave w.r.t. $y$. These are important minimax problems but the studies on their solution methods are still quite limited. Although there exists some symmetry between nonconvex-concave and convex-nonconcave minimax problems,
the complexity analysis of the same AGP algorithm for nonconvex-concave setting cannot be trivially extended to that for convex-nonconcave setting.

\subsection{Strongly Convex-Nonconcave Setting}\label{subsec4}
In this subsection, we analyze the iteration complexity of Algorithm \ref{alg:2.1} for solving strongly convex-nonconcave minimax optimization problems \eqref{problem:1}, i.e., $f(x,y)$ is $\theta$-strongly convex w.r.t. $x$ for any fixed $y\in \mathcal{Y}$,
and nonconcave  w.r.t. $y$ for any given $x\in \mathcal{X}$. Under this setting,   $\forall k\ge 1$, we set
\begin{equation}
\beta_k=\tfrac{1}{\zeta}, \gamma_k=\nu,  b_k=c_k=0,\label{subsec:4.1:1:para}
\end{equation}
in  Algorithm \ref{alg:2.1}, and simplify the update for $x_k$ and $y_k$ as follows:
\begin{align}
	x_{k+1}&=\mathcal{P}_\mathcal{X}\left( x_k- \zeta \nabla _x {f}( x_{k},y_k) \right),\label{subsec:4.1:2:x}\\
	y_{k+1}&=\mathcal{P}_\mathcal{Y} \left(y_k + \tfrac{1}{\nu}  \nabla _yf ( x_{k+1},y_k )  \right)\label{subsec:4.1:3:y}.
\end{align}
%
%
Our goal in the remaining part of this subsection is to establish the iteration complexity of Algorithm \ref{alg:2.1} under
the strongly convex-nonconcave setting.
The convergence analysis for this setting is different from that of the nonconvex-strongly concave setting in Section 3.
	\begin{lem}\label{lem:4.2}
	Suppose that Assumption \ref{sec:2:ass1:lip} holds. Let $\{\left(x_{k},y_{k}\right)\}$ be a sequence generated by Algorithm \ref{alg:2.1}
	with parameter settings in \eqref{subsec:4.1:1:para}. If $\nu > L_{22}$, then we have
	\begin{align}\label{lem:4.2:1:conclu}
		&f( x_{k+1},y_{k+1}) -f( x_{k},y_{k})  \nonumber \\
		\geq& \frac{\nu}{2}\| y_{k+1}-y_{k} \| ^2-\frac{L_{21}^{2}\zeta}{2}\| y_{k}-y_{k-1} \|^2+(\theta-\frac{1}{2\zeta}-\frac{\zeta L_{11}^2}{2})\| x_{k}-x_{k-1}\|^2\nonumber\\
		&+ (\frac{\theta}{2}-\frac{1}{\zeta})\|x_{k+1}-x_{k} \|^2.
	\end{align}
\end{lem}

\begin{proof}
	Similar to \eqref{lem:3.2:2}-\eqref{lem:3.2:4} in the proof of Lemma \ref{lem:3.2:blostr:y}, by the optimality condition for $x_k$ in \eqref{subsec:4.1:2:x} implies that $\forall x\in \mathcal{X}$ and $\forall k\geq 1$,
	\begin{align}
		\langle \nabla _xf(x_{k},y_{k})+\frac{1}{\zeta}(x_{k+1}-x_{k}),x-x_{k+1} \rangle &\ge 0,\label{lem:4.2:2:opty}\\
		\langle \nabla _xf(x_{k},y_{k})+\frac{1}{\zeta}(x_{k+1}-x_{k}) ,x_{k}-x_{k+1} \rangle &\ge 0,\label{lem:4.2:3:nopty}\\
		\langle \nabla _xf(x_{k-1},y_{k-1})+\frac{1}{\zeta}(x_{k}-x_{k-1}),x_{k+1}-x_{k} \rangle &\ge 0,\label{lem:4.2:4:opty:k}
	\end{align}
which, in view of the fact that $f\left( x,y  \right)$ is $\theta $-strongly convex w.r.t. $x$ for any given $y\in \mathcal{Y}$,
	then implies that
	\begin{align}\label{lem:4.2:5}
		&f( x_{k+1},y_{k}) -f( x_{k},y_{k} ) \nonumber\\
		\geq &  \langle  \nabla _xf( x_{k},y_{k} ),x_{k+1}-x_{k}\rangle  +\frac{\theta}{2}\| x_{k+1}-x_{k} \|^2\  \nonumber\\
		\geq &\langle  \nabla _xf( x_{k},y_{k} ) -\nabla _xf( x_{k-1},y_{k-1}) ,x_{k+1}-x_{k} \rangle  \nonumber\\
		& -\frac{1}{\zeta}\langle  x_{k}-x_{k-1},x_{k+1}-x_{k} \rangle +\frac{\theta}{2}\| x_{k+1}-x_{k} \|^2.
	\end{align}
	Denoting $m_{k+1}:=\left(x_{k+1}-x_k\right)-\left(x_k-x_{k-1}\right)$, by Assumption \ref{sec:2:ass1:lip}, the Cauchy-Schwarz inequality and
	the $\theta $-strongly convexity  of $f$ w.r.t. $x$, we can estimate the first inner product term in the r.h.s.
	of \eqref{lem:4.2:5} as
	\begin{align}
		&\langle  \nabla _x f ( x_{k},y_{k} ) -\nabla _x f ( x_{k-1},y_{k-1}) ,x_{k+1}-x_{k} \rangle  \nonumber\\
		=&\langle  \nabla _x f ( x_{k},y_{k} ) -\nabla _x f \left( x_{k},y_{k-1} \right) ,x_{k+1}-x_{k} \rangle  \nonumber\\
		& + \langle  \nabla _x f \left( x_{k},y_{k-1} \right) -\nabla _x f ( x_{k-1},y_{k-1}) ,m_{k+1} \rangle  \nonumber\\
		&+\langle  \nabla _x f \left( x_{k},y_{k-1}\right) -\nabla _x f ( x_{k-1},y_{k-1}) ,x_{k}-x_{k-1} \rangle \label{lem:4.2:6}\\
	\geq &-\frac{L_{21}^{2}\zeta}{2}\| y_{k}-y_{k-1} \|^2-\frac{1}{2\zeta}\| x_{k+1}-x_k \|^2-\frac{\zeta L_{11}^{2}}{2}\|x_k-x_{k-1} \|^2\nonumber\\
	&-\frac{1}{2\zeta}\|m_{k+1} \|^2 +	\theta \| x_k-x_{k-1} \|^2.\label{lem:4.2:7}
	\end{align}
	Moreover, it can be easily checked that
\begin{align}\label{lem:4.2:8:3points}
		&\langle  x_k-x_{k-1},x_{k+1}-x_k \rangle =\frac{1}{2}\| x_k-x_{k-1} \|^2+\frac{1}{2}\| x_{k+1}-x_k \|^2-\frac{1}{2}\| m_{k+1} \|^2.
	\end{align}
	Plugging \eqref{lem:4.2:7} and \eqref{lem:4.2:8:3points} into \eqref{lem:4.2:5} and rearranging the terms, we conclude that
	\begin{align}\label{lem:4.2:9}
	f( x_{k+1},y_{k}) -f( x_{k},y_{k})	\geq& -\frac{L_{21}^{2}\zeta}{2}\| y_{k}-y_{k-1} \|^2+(\theta-\frac{1}{2\zeta}-\frac{\zeta L_{11}^2}{2})\| x_{k}-x_{k-1}\|^2\nonumber\\
		&+ (\frac{\theta}{2}-\frac{1}{\zeta})\|x_{k+1}-x_{k} \|^2.
	\end{align}
By setting $\gamma_k=\nu$, $c_k=0$ in \eqref{lem:2.1:2} of Lemma \ref{lem:2.1}  and the assumption $\nu > L_{22}$, we have
		\begin{equation}\label{lem:4.1:1:conclu}
			f( x_{k+1},y_{k+1}) -f( x_{k+1},y_k )\ge \frac{\nu}{2}\| y_{k+1}-y_k \|^2.
		\end{equation}
	The proof is completed by combining \eqref{lem:4.2:9} with \eqref{lem:4.1:1:conclu}.
\end{proof}

Next, we further refine this relation in \eqref{lem:4.2:1:conclu} as shown below. Note that in Lemma \ref{lem:3.3:blostr:F} we have constructed a potential function for
the nonconvex-strongly concave case, whereas in the following lemma, an additional term involving the distance between two adjacent iterates, i.e., $\|x_{k+1}-x_k\|^2$, is added to construct another potential function for
the strongly convex-nonconcave case.

\begin{lem}\label{lem:4.3}
	Suppose that Assumption \ref{sec:2:ass1:lip} holds. Let $\{\left(x_k,y_k\right)\}$ be a sequence generated by Algorithm \ref{alg:2.1}
	with parameter settings in \eqref{subsec:4.1:1:para}. Denote
	$$\hat{f}_{k+1}:=f\left( x_{k+1},y_{k+1} \right),\ \hat{S}_{k+1}:=-\frac{2}{\zeta^2\theta}\lVert x_{k+1}-x_k \rVert^2,$$
	$$\hat{F}_{k+1}:=\hat{f}_{k+1}+\hat{S}_{k+1}+ (\theta+\frac{7}{2\zeta}-\frac{\zeta L_{11}^2}{2}-\frac{2L_{11}^2}{\theta})\lVert x_{k+1}-x_{k} \rVert^2-(\frac{L_{21}^2\zeta}{2}+\frac{2L_{21}^2}{\theta^2\zeta})\|y_{k+1}-y_k\|^2.$$
	If $\nu>L_{22}$,
	then $\forall k \geq 1$,
	\begin{align}\label{lem:4.3:2:conclus}
		\hat{F}_{k+1}-\hat{F}_k\geq& \left(\frac{\nu}{2}-\frac{\zeta L_{21}^2}{2}-\frac{2L_{21}^2}{\zeta \theta^2}\right)\lVert y_{k+1}-y_k \rVert ^2\nonumber\\
		& +\left(\frac{3\theta-\zeta L_{11}^2}{2}+\frac{\theta-4\zeta L_{11}^2}{2\zeta\theta}\right)\lVert x_{k+1}-x_k \rVert ^2.
	\end{align}
\end{lem}

\begin{proof}
First by \eqref{lem:4.2:3:nopty} and \eqref{lem:4.2:4:opty:k}, we have
	\begin{align}\label{lem:4.3:3}
		\frac{1}{\zeta}\langle m_{k+1}, x_k -x_{k+1} \rangle
		\ge \langle \nabla _xf( x_k,y_k ) -\nabla _xf\left( x_{k-1},y_{k-1} \right),x_{k+1}-x_k \rangle,
	\end{align}
which together with \eqref{lem:4.2:6} then imply that
	\begin{align}
		\frac{1}{\zeta}\langle m_{k+1},x_k-x_{k+1}  \rangle  \geq&\langle  \nabla _xf( x_{k},y_{k} ) -\nabla _xf\left( x_{k},y_{k-1} \right) ,x_{k+1}-x_{k} \rangle  \nonumber\\
		& + \langle  \nabla _xf\left( x_{k},y_{k-1} \right) -\nabla _xf\left( x_{k-1},y_{k-1} \right) ,m_{k+1} \rangle  \nonumber\\
		&+\langle  \nabla _xf\left( x_k,y_{k-1} \right) -\nabla _xf\left( x_{k-1},y_{k-1} \right) ,x_{k}-x_{k-1} \rangle \label{lem:4.3:4}\\
		\geq& -\frac{L_{21}^{2}}{2 \theta}\|y_{k}-y_{k-1} \|^2-\frac{\theta}{2}\|x_{k+1}-x_{k} \|^2-\frac{\zeta L_{11}^2}{2}\|x_{k}-x_{k-1} \|^2\nonumber\\
		&-\frac{1}{2\zeta}\| m_{k+1}\|^2+\theta\|x_{k}-x_{k-1} \|^2,\label{lem:4.3:5}
	\end{align}
where the second inequality is similar to the proof of \eqref{lem:4.2:7} except that for the first term in the r.h.s., we use $\langle  \nabla _xf( x_{k},y_{k} ) -\nabla _xf\left( x_{k},y_{k-1} \right) ,x_{k+1}-x_{k} \rangle \geq -\frac{L_{21}^{2}}{2 \theta}\|y_{k}-y_{k-1} \|^2-\frac{\theta}{2}\|x_{k+1}-x_k \|^2$. By using the identity $\frac{1}{\zeta}\langle m_{k+1},x_{k}-x_{k+1} \rangle  =\frac{1}{2\zeta}\|x_{k}-x_{k-1}\|^2- \frac{1}{2\zeta}\|x_{k+1}-x_{k}\|^2-\frac{1}{2\zeta}\| m_{k+1}\| ^2$, we conclude from \eqref{lem:4.3:5} that
	\begin{align}\label{lem:4.3:6}
		&\frac{1}{2\zeta}\|x_{k}-x_{k-1}\|^2- \frac{1}{2\zeta}\|x_{k+1}-x_{k}\|^2-\frac{1}{2\zeta}\| m_{k+1}\| ^2\nonumber\\
		\geq& -\frac{L_{21}^{2}}{2 \theta}\|y_{k}-y_{k-1} \|^2-\frac{\theta}{2}\|x_{k+1}-x_k \|^2-\frac{\zeta L_{11}^2}{2}\|x_{k}-x_{k-1} \|^2\nonumber\\
		&-\frac{1}{2\zeta}\| m_{k+1}\|^2+\theta\|x_{k}-x_{k-1} \|^2.
	\end{align}
	Rearranging the terms of \eqref{lem:4.3:6}, we have
	\begin{align}\label{lem:4.3:7}
		&\frac{1}{2\zeta}\|x_{k}-x_{k-1}\|^2- \frac{1}{2\zeta}\|x_{k+1}-x_{k}\|^2\nonumber\\
		\geq& -\frac{L_{21}^{2}}{2 \theta}\|y_{k}-y_{k-1} \|^2-\frac{\theta}{2}\|x_{k+1}-x_{k} \|^2+\left(\theta-\frac{\zeta L_{11}^2}{2}\right)\|x_{k}-x_{k-1} \|^2.
	\end{align}
	Multiplying $\frac{4}{\zeta\theta}$ on both sides of \eqref{lem:4.3:7} and using the definition of $\hat{S}_{k+1}$, we obtain
	\begin{align*}
		\hat{S}_{k+1}-\hat{S}_k \geq & -\frac{2L_{21}^{2}}{\theta^2\zeta}\| y_{k}-y_{k-1}\|^2-\frac{2}{\zeta}\| x_{k+1}-x_k \|^2+\left(\frac{4}{\zeta}-\frac{2 L_{11}^2}{\theta}\right)\|x_{k}-x_{k-1} \|^2.
	\end{align*}
	The proof is completed by \eqref{lem:4.2:1:conclu} in Lemma \ref{lem:4.2} and the definition of $\hat{F}_k$.
\end{proof}

	We are now ready to establish the iteration complexity for the AGP algorithm in the strongly convex-nonconcave setting.

\begin{thm}\label{thm:4.1}
	Suppose that Assumption \ref{sec:2:ass1:lip} holds. Let $\{\left(x_{k},y_{k}\right)\}$ be a sequence generated by Algorithm \ref{alg:2.1}
	with parameter settings in \eqref{subsec:4.1:1:para}.
	If the relations $\nu>L_{22}, \nu>L_{21}^2\zeta+\frac{4L_{21}^2}{\zeta \theta^2},~\zeta\leq \frac{\theta}{4L_{11}^2}$ are satisfied, then $\forall \varepsilon>0$,  it holds that $$ T\left( \varepsilon \right) \le
	\tfrac{\overline{F}-\hat{F}_1}{\hat{d}_1\varepsilon ^2},$$
	where
	$\hat{d}_1:=\tfrac{\min \left\{\tfrac{\nu}{2}-\tfrac{\zeta L_{21}^2}{2}-\tfrac{2 L_{21}^2}{\zeta \theta^2},\tfrac{3\theta-\zeta L_{11}^2}{2}+\tfrac{\theta-4\zeta L_{11}^2}{2\zeta\theta} \right\}}{\max \left\{ \tfrac{1}{\zeta ^2}+2L_{12}^{2},2\nu^2 \right\}}
	$
	and $\overline{F}=\bar{f}+(\theta+\frac{7}{2\zeta}-\frac{\zeta L_{11}^2}{2}-\frac{2L_{11}^2}{\theta})\sigma_x^2$ with $\overline{f} := \max_{(x,y) \in \mathcal{X} \times \mathcal{Y}} f(x,y)$ and $\sigma _x= \max\{\|x_1-x_2\| \mid \forall x_1, x_2\in \mathcal{X}\}$.
\end{thm}

\begin{proof}
	By \eqref{subsec:4.1:2:x}, we immediately obtain
	\begin{equation}
		\|(\nabla G_k)_{x} \|\leq\frac{1}{\zeta}\| x_{k+1}-x_{k} \|.\label{thm:4.1:1:gx}
	\end{equation}
	On the other hand, similar to the proof of \eqref{thm:3.1:2:gy},  by \eqref{subsec:4.1:3:y} and the triangle inequality, and the nonexpansiveness of the projection operator $\operatorname{ \mathcal{P}_\mathcal{Y}}$, we conclude that
	\begin{align}\label{thm:4.1:2:gy}
		\| (\nabla G_k)_{y} \|\le& \nu\| y_{k+1}-y_k \|+L_{12}\|x_{k+1}-x_{k} \|.
			\end{align}
	By combining \eqref{thm:4.1:1:gx} and \eqref{thm:4.1:2:gy}, and using the Cauchy-Schwarz inequality, we obtain
	\begin{align}\label{thm:4.1:3:g}
		\| \nabla G_k \|^2
		\leq 2\nu^2\| y_{k+1}-y_{k} \|^2+\left( \frac{1}{\zeta ^2}+2L_{12}^{2} \right)\|x_{k+1}-x_{k} \|^2.
	\end{align}
	Observing that $\hat{d}_1 > 0$. Multiplying both sides of \eqref{thm:4.1:3:g} by $\hat{d}_1$, and using \eqref{lem:4.3:2:conclus} in Lemma \ref{lem:4.3},
	we have
	\begin{align}\label{thm:4.1:4:d1}
		\hat{d}_1\|\nabla G_k \|^2\le \hat{F}_{k+1}-\hat{F}_k.
	\end{align}
	Summing up the above inequalities from $k=1$ to $k=T(\varepsilon)$, we obtain
	\begin{align}\label{thm:4.1:5}
		\sum_{k=1}^{T\left( \varepsilon \right)}{\hat{d}_1\|\nabla G_k \| ^2}\le \hat{F}_{T\left( \varepsilon \right)+1}-\hat{F}_1.
	\end{align}
	Note that by the definition of $\hat{F}_{k+1}$ in Lemma \ref{lem:4.3}, we have
	\begin{align*}
		\hat{F}_{T\left( \varepsilon \right)+1}&\le
		\hat{f}_{T\left( \varepsilon \right)+1}+(\theta+\frac{7}{2\zeta}-\frac{\zeta L_{11}^2}{2}-\frac{2L_{11}^2}{\theta})\|x_{T\left( \varepsilon \right)+1}-x_{T\left( \varepsilon \right)} \|^2\\
		&\le \bar{f}+(\theta+\frac{7}{2\zeta}-\frac{\zeta L_{11}^2}{2}-\frac{2L_{11}^2}{\theta})\sigma_x^2 = \bar{F},
	\end{align*}
	where the last inequality follows from the definitions of $\bar{f}$ and  $\sigma _x$, and $\theta+\frac{7}{2\zeta}-\frac{\zeta L_{11}^2}{2}-\frac{2L_{11}^2}{\theta} \ge 0 $ due to
	the selection of $\zeta$.
	We then conclude from \eqref{thm:4.1:5} that $\sum_{k=1}^{T\left( \varepsilon \right)}{\hat{d}_1\| \nabla G_k \|^2}\le \bar{F}-\hat{F}_1$ which, in view of the definition of $T(\varepsilon)$, implies that $\varepsilon ^2\le (\bar{F}-\hat{F}_1)/(T( \varepsilon )\cdot \hat{d}_1)$
	or equivalently, $T\left( \varepsilon \right) \le (\bar{F}-\hat{F}_1)/(\hat{d}_1\varepsilon ^2)$.
\end{proof}

Theorem \ref{thm:4.1} shows that the number of gradient evaluations performed by Algorithm \ref{alg:2.1} to obtain an $\varepsilon$-stationary point of $f$ is bounded by $\mathcal{O}\left(L^2 \varepsilon ^{-2} \right)$ under the strongly convex-nonconcave setting. To the best of our knowledge, this is the first theoretical guarantee that has been obtained in the literature for solving this class of minimax problems.

	\subsection{Complexity Analysis for Convex-Nonconcave Setting}\label{subsec4-s}
	In this subsection, we analyze the iteration complexity of Algorithm \ref{alg:2.1} applied to
	the general convex-nonconcave setting, for which $f(x,y)$ is convex w.r.t. $x$ for any fixed $y\in \mathcal{Y}$, and nonconcave w.r.t. $y$ for any given $x\in \mathcal{X}$.
	Under this setting,  $\forall k\ge 1$, we set
	\begin{equation}
		\beta_k=\tfrac{1}{\bar{\zeta}}, \gamma_k=\bar{\nu}+\bar{\gamma}_k,  b_k=q_k, c_k=0,\label{subsec:4.2:1:para}
	\end{equation}
	where $\bar{\gamma}_k$ and $q_k$ are stepsize parameters to be defined later. We need to make the following assumption on the parameters $q_k$.
	\begin{ass}\label{subsec:4.2:ass1:qk}
		$\{q_k\}$ is a nonnegative monotonically decreasing sequence.
	\end{ass}

	By Assumption \ref{sec:2:ass1:lip} and $\nabla _xf_{k-1}\left( x,y \right) =\nabla _xf\left( x,y \right) +q_{k-1}x$, we have
\begin{align}\label{subsec:4.2:2:Lx}	
	&\|\nabla_{x} f_{k-1}(x_{k},y_{k-1})-\nabla_{x}f_{k-1}(x_{k-1},y_{k-1})\|	\leq \left(L_{11}+q_{k-1}\right)\| x_{k}-x_{k-1}\|.
\end{align}
Denoting $L_{11}^{'}=L_{11}+q_1$, by Assumption \ref{subsec:4.2:ass1:qk} and  \eqref{subsec:4.2:2:Lx}, we have
$$\|\nabla_{x} f_{k-1}(x_{k},y_{k-1})-\nabla_{x}f_{k-1}(x_{k-1},y_{k-1})\|\leq L_{11}^{'}\| x_{k}-x_{k-1}\}.$$
It then follows from the above inequality and the strong convexity of $f_{k-1}\left(x,y_{k-1} \right)$ w.r.t. $x$ (Theorem 2.1.12 in \cite{Nestrov}) that
\begin{align}\label{subsec:4.2:3:key}
	&\quad\langle \nabla_{x} f_{k-1}(x_{k},y_{k-1})-\nabla_{x}f_{k-1}(x_{k-1},y_{k-1}) ,x_{k}-x_{k-1} \rangle \nonumber\\
	&\geq \frac{1}{L_{11}^{'}+q_{k-1}}\|\nabla_{x} f_{k-1}(x_{k},y_{k-1})-\nabla_{x}f_{k-1}(x_{k-1},y_{k-1})  \|^2 + \frac{q_{k-1}L_{11}^{'}}{L_{11}^{'}+q_{k-1}}\|x_{k}-x_{k-1}\|^2.
\end{align}	
    By using the strong convexity of  $f_{k-1}\left(x,y_{k-1} \right)$ w.r.t. $x$ instead of the strong concavity of $f_{k-1}\left(x_{k},y \right)$ w.r.t. $y$, this inequality provides a lower bound for the inner product instead of an upper bound shown as in \eqref{subsec:3.2:3:key}. This is a key inequality that we will use to establish some important recursions for the AGP algorithm under the convex-nonconcave setting in the following two results.
	
	Similar to Lemma \ref{lem:4.2}, we first provide an estimate on the increase of the function values from $f(x_k,y_k)$ to $f(x_{k+1},y_{k+1})$.
	\begin{lem}\label{lem:4.5}
		Suppose that Assumption \ref{sec:2:ass1:lip} and \ref{subsec:4.2:ass1:qk} hold. Let $\{\left(x_k,y_k\right)\}$ be a sequence generated by Algorithm \ref{alg:2.1}
		with parameter settings in \eqref{subsec:4.2:1:para}. If $\forall k,\bar{\gamma}_k >L_{22}$ and $\bar{\zeta} \le  \frac{2}{L_{11}^{'}+q_1}$, then
		\begin{align}\label{lem:4.5:1:conclu}
			&f(x_{k+1},y_{k+1})-f(x_k,y_k)\nonumber \\
\geq& \left( \bar{\nu} +\frac{\bar{\gamma}_k}{2} \right) \| y_{k+1}-y_k \| ^2-\frac{L_{21}^{2}\bar{\zeta}}{2}\| y_k-y_{k-1}\|^2-\frac{1}{\bar{\zeta}} \|x_{k+1}-x_k \| ^2-\frac{1}{2\bar{\zeta}}\| x_k-x_{k-1} \|^2\nonumber\\
			&-\frac{q_{k-1}}{2}(\| x_{k+1}\| ^2-\| x_k\|^2).
		\end{align}
	\end{lem}
	
	\begin{proof}
		Similar to \eqref{lem:4.2:2:opty}-\eqref{lem:4.2:6} in the proof of Lemma \ref{lem:4.2}, by replacing $f$ with $f_k$, $\zeta$ with $\bar{\zeta}$ respectively, and setting $\theta=0$, we obtain that
%
		\begin{align}\label{lem:4.5:2}
			&f_k(x_{k+1},y_k)-f_k(x_k,y_k)\nonumber\\
			\geq &\langle \nabla _xf_k(x_k,y_k)-\nabla _xf_{k-1}(x_k,y_{k-1}),x_{k+1}-x_k \rangle+\langle \nabla _xf_{k-1}(x_k,y_{k-1})-\nabla _xf_{k-1}(x_{k-1},y_{k-1}),m_{k+1} \rangle\nonumber\\
			&+\langle \nabla _xf_{k-1}(x_k,y_{k-1})-\nabla _xf_{k-1}(x_{k-1},y_{k-1}),x_k-x_{k-1} \rangle-\frac{1}{\bar{\zeta}}\langle x_k-x_{k-1},x_{k+1}-x_k \rangle.
		\end{align}
 Next, we prove  lower bound for the four  terms in the r.h.s. of \eqref{lem:4.5:2} which are different from that in Lemma \ref{lem:4.2}.
 By using the convexity of $f_{k-1}\left(x,y_{k-1} \right)$ w.r.t $x$ instead of the concavity of $f_{k-1}\left(x_{k},y \right)$ w.r.t $y$, the opposite side Cauchy-Schwarz inequality, and replacing  $\bar{\rho}, c_{k}, c_{k-1},L_{12}$, $L_{22}^{'}$ by  $\bar{\zeta},q_{k}, q_{k-1},L_{21}$, $L_{11}^{'}$ respectively, Assumptions \ref{sec:2:ass1:lip} and \ref{subsec:4.2:ass1:qk}, similar to the proof of \eqref{lem:3.5:6}-\eqref{lem:3.5:8} we conclude that
		\begin{align}\label{lem:4.5:3}
			&\langle  \nabla _xf_k( x_k,y_k) -\nabla _xf_{k-1}( x_k,y_{k-1}) ,x_{k+1}-x_k \rangle  \nonumber\\
\geq&-\frac{L_{21}^{2}\bar{\zeta}}{2}\| y_k-y_{k-1} \|^2-\frac{1}{2\bar{\zeta}}\| x_{k+1}-x_k \|^2+\frac{q_k-q_{k-1}}{2}(\| x_{k+1}\| ^2-\| x_k\|^2),
		\end{align}
		and
		\begin{align}\label{lem:4.5:4}
			\langle  \nabla _xf_{k-1}( x_k,y_{k-1}) -\nabla _xf_{k-1}( x_{k-1},y_{k-1}) ,m_{k+1} \rangle &	\geq  -\frac{\bar{\zeta}}{2}\|\nabla _xf_{k-1}( x_k,y_{k-1}) -\nabla _xf_{k-1}( x_{k-1},y_{k-1}) \|^2 \nonumber\\
			&\quad -\frac{1}{2\bar{\zeta}}\|m_{k+1} \|^2,
		\end{align}
	and
		\begin{align}\label{lem:4.5:5}
			&\langle \nabla _xf_{k-1}\left( x_k,y_{k-1} \right) -\nabla _xf_{k-1}\left( x_{k-1},y_{k-1} \right) ,x_k-x_{k-1} \rangle \nonumber\\
			& \geq  \frac{1}{L_{11}^{'}+q_{k-1}}\| \nabla _xf_{k-1}\left( x_k,y_{k-1} \right)-\nabla _xf_{k-1}\left( x_{k-1},y_{k-1} \right)\|^2.
		\end{align}
	Moreover, it can be easily checked that
		\begin{equation}\label{lem:4.5:6}
			\frac{1}{\bar{\zeta}}\langle  x_{k+1}-x_k,x_k-x_{k-1} \rangle =\frac{1}{2\bar{\zeta}}\|x_{k+1}-x_k \|^2+\frac{1}{2\bar{\zeta}}\|x_k-x_{k-1} \|^2-\frac{1}{2\bar{\zeta}}\|m_{k+1} \|^2.
		\end{equation}
		Plugging \eqref{lem:4.5:3}-\eqref{lem:4.5:6} into \eqref{lem:4.5:2}, and by the definition of $f_k(x_{k+1},y_k)$ and $ f_k(x_k,y_k)$ and  $\frac{\bar{\zeta}}{2} \le \frac{1}{L_{11}^{'}+q_1}$, we conclude that
		\begin{align}\label{lem:4.5:7}
f(x_{k+1},y_k)-f(x_k,y_k) \geq& -\frac{L_{21}^{2}\bar{\zeta}}{2}\| y_k-y_{k-1}\|^2-\frac{1}{\bar{\zeta}} \|x_{k+1}-x_k \| ^2-\frac{1}{2\bar{\zeta}}\| x_k-x_{k-1} \|^2-\frac{q_{k-1}}{2}(\| x_{k+1}\| ^2-\| x_k\|^2).
		\end{align}
By setting $\gamma_k=\bar{\nu}+\bar{\gamma}_k$, $c_k=0$ in \eqref{lem:2.1:2} of Lemma \ref{lem:2.1}  and the assumption $\bar{\gamma}_k > L_{22}$, we have
		\begin{equation}\label{lem:4.4:1:conclu}
			f(x_{k+1},y_{k+1})-f(x_{k+1},y_k)\ge \left(\bar{\nu} +\frac{\bar{\gamma}_k}{2} \right) \| y_{k+1}-y_k \|^2.
		\end{equation}
	The proof is completed by combing \eqref{lem:4.5:7} and \eqref{lem:4.4:1:conclu}.
	\end{proof}
	
We need to further refine the relation in \eqref{lem:4.5:1:conclu} in order to establish the convergence of the AGP algorithm as shown below.	Note that in Lemma \ref{lem:4.3} we have constructed a potential function which involves the function value plus some distance between two adjacent iterates for the strongly convex-nonconcave case, whereas in the following lemma an additional quadratic regularization term, i.e., $\|x_{k+1}\|^2$, is added to construct another potential function for the general convex-nonconcave case.
	\begin{lem}\label{lem:4.6}
		Suppose that Assumptions \ref{sec:2:ass1:lip} and \ref{subsec:4.2:ass1:qk} hold. Let $\{\left(x_k,y_k\right)\}$ be a sequence generated by Algorithm \ref{alg:2.1}
		with parameter settings in \eqref{subsec:4.2:1:para}. Denote
		\begin{align*}
			\hat{\mathcal{S}}_{k+1}&:=-\frac{8}{\bar{\zeta}^2q_{k+1}}\| x_{k+1}-x_k \|^2-\frac{8}{\bar{\zeta}}\left(1-\frac{q_k}{q_{k+1}}\right)\|x_{k+1}\|^2,\\
			\hat{\mathcal{F}}_{k+1}&:=f(x_{k+1},y_{k+1})+\hat{\mathcal{S}}_{k+1}+ \frac{15}{2\bar{\zeta}}\| x_{k+1}-x_k \|^2+\frac{q_k}{2}\|x_{k+1}\|^2\\
			&\quad-\left( \frac{\bar{\zeta} L_{21}^2}{2}+\frac{16L_{21}^2}{\bar{\zeta} (q_{k+1})^2} \right)\| y_{k+1}-y_k \| ^2.
		\end{align*}
		If
		\begin{equation}\label{lem:4.6:1:para}
			\bar{\gamma}_k >L_{22}, ~\frac{1}{q_{k+1}}-\frac{1}{q_k}\leq \frac{\bar{\zeta}}{5},~\bar{\zeta} \leq \frac{2}{L_{11}^{'}+q_1},
		\end{equation}
		then $\forall k \geq 1$,
		\begin{align*}
			&\hat{\mathcal{F}}_{k+1}-\hat{\mathcal{F}}_{k} \nonumber\\
		\geq 	&\left( \bar{\nu}+\frac{\bar{\gamma}_k}{2}-\frac{\bar{\zeta} L_{21}^2}{2}-\frac{16L_{21}^2}{\bar{\zeta} (q_{k+1})^2}\right)\| y_{k+1}-y_k \| ^2+\frac{q_k-q_{k-1}}{2}\|x_{k+1}\|^2\nonumber\\
			&+\frac{9}{10\bar{\zeta}}\| x_{k+1}-x_k \|^2 +\frac{8}{\bar{\zeta}}\left(\frac{q_k}{q_{k+1}}- \frac{q_{k-1}}{q_k} \right)\|x_{k+1}\|^2.
		\end{align*}
	\end{lem}
	
	\begin{proof}
		Similar to the proof of \eqref{lem:4.3:4}, by replacing $f$ with $f_k$, we have
		\begin{align*}
			\frac{1}{\bar{\zeta}}\langle m_{k+1}, x_k -x_{k+1} \rangle\geq  & \langle \nabla _xf_k(x_k,y_k)-\nabla _xf_{k-1}(x_k,y_{k-1}),x_{k+1}-x_k \rangle\nonumber\\
&+\langle \nabla _xf_{k-1}(x_k,y_{k-1})-\nabla _xf_{k-1}(x_{k-1},y_{k-1}),m_{k+1} \rangle\nonumber\\
			& +\langle \nabla _xf_{k-1}(x_k,y_{k-1})-\nabla _xf_{k-1}(x_{k-1},y_{k-1}),x_k-x_{k-1} \rangle.
		\end{align*}
		Using an argument similar to the proof of \eqref{lem:4.5:3}-\eqref{lem:4.5:6},
		and the Cauchy-Schwarz inequality, we conclude from the above inequality that
		\begin{align}\label{lem:4.6:2}	&\frac{1}{2\bar{\zeta}}\|x_k-x_{k-1}\|^2-\frac{1}{2\bar{\zeta}}\|x_{k+1}-x_k\|^2-\frac{1}{2\bar{\zeta}}\|m_{k+1}\|^2\nonumber\\
\geq& -\frac{L_{21}^{2}}{2\bar{a}_k}\|y_k-y_{k-1} \|^2-\frac{\bar{a}_k}{2}\| x_{k+1}-x_k \|^2+\frac{q_k-q_{k-1}}{2}(\| x_{k+1}\|^2-\| x_k\| ^2)\nonumber\\
			&- \frac{q_k-q_{k-1}}{2}\|x_{k+1}-x_k\|^2-\frac{\bar{\zeta}}{2}\|  \nabla _xf_{k-1}(x_k,y_{k-1})-\nabla _xf_{k-1}(x_{k-1},y_{k-1})\|^2\nonumber\\
			&+\frac{1}{L_{11}^{'}+q_{k-1}}\|\nabla _xf_{k-1}(x_k,y_{k-1})-\nabla _xf_{k-1}(x_{k-1},y_{k-1}) \|^2-\frac{1}{2\bar{\zeta}}\|m_{k+1} \|^2\nonumber\\
			&+\frac{q_{k-1}L_{11}^{'}}{L_{11}^{'}+q_{k-1}}\| x_k-x_{k-1} \|^2,
		\end{align}
		where $\bar{a}_k>0$. Observing that $q_1\leq L_{11}^{'}$ and by Assumption \ref{subsec:4.2:ass1:qk}, we have
		\begin{equation*}
			\frac{q_{k-1}L_{11}^{'}}{q_{k-1}+L_{11}^{'}}\geq\frac{q_{k-1}L_{11}^{'}}{2L_{11}^{'}}
			=\frac{q_{k-1}}{2}\geq\frac{q_k}{2}.
		\end{equation*}
	By $\bar{\zeta}\le \frac{2}{L_{11}^{'}+q_1}$ and Assumption \ref{subsec:4.2:ass1:qk}, rearranging the terms in \eqref{lem:4.6:2}, we obtain
		\begin{align*}
			& -\frac{1}{2\bar{\zeta}}\| x_{k+1}-x_k \|^2-\frac{q_k-q_{k-1}}{2}\|x_{k+1}\|^2 \nonumber\\	
			\geq &-\frac{1}{2\bar{\zeta}}\|x_k-x_{k-1}\|^2-\frac{q_k-q_{k-1}}{2}\|x_k\|^2-\frac{L_{21}^{2}}{2\bar{a}_k}\|y_k-y_{k-1} \|^2\nonumber\\
			& -\frac{\bar{a}_k}{2}\| x_{k+1}-x_k \|^2+\frac{q_k}{2}\|x_k-x_{k-1} \|^2.
		\end{align*}
		By multiplying $\frac{16}{\bar{\zeta}q_k}$ on both sides of the above inequality, we then obtain
		\begin{align}\label{lem:4.6:3}
			&-\frac{8}{\bar{\zeta}^2q_k}\| x_{k+1}-x_k \|^2-\frac{8}{\bar{\zeta}}\left(1-\frac{q_{k-1}}{q_k}\right)\|x_{k+1}\|^2\nonumber\\
		\geq 	&  -\frac{8}{\bar{\zeta}^2q_k}\|x_k-x_{k-1}\|^2 -\frac{8}{\bar{\zeta}}(1-\frac{q_{k-1}}{q_k}
			)\|x_k\|^2-\frac{8L_{21}^2}{\bar{\zeta} q_k\bar{a}_k}\|y_k-y_{k-1} \|^2\nonumber\\
			& -\frac{8\bar{a}_k}{\bar{\zeta} q_k}\| x_{k+1}-x_k \|^2 +\frac{8}{\bar{\zeta}}\|x_k-x_{k-1}\|^2.
		\end{align}
		Setting $\bar{a}_k=\frac{q_k}{2}$ in the above inequality, and using the definition of $\hat{\mathcal{S}}_{k+1}$ and \eqref{lem:4.6:1:para}, we have
		\begin{align}\label{lem:4.6:4}
			\hat{\mathcal{S}}_{k+1}-\hat{\mathcal{S}}_{k}
		\geq 	& \frac{8}{\bar{\zeta}}\left(\frac{q_k}{q_{k+1}}- \frac{q_{k-1}}{q_k} \right)\|x_{k+1}\|^2-\frac{16L_{21}^2}{\bar{\zeta} (q_k)^2}\|y_k-y_{k-1} \|^2\nonumber\\
			& -\frac{28}{5\bar{\zeta}}\| x_{k+1}-x_k \|^2 +\frac{8}{\bar{\zeta}}\|x_k-x_{k-1}\|^2.
		\end{align}
	The proof is completed by combining \eqref{lem:4.6:4} and \eqref{lem:4.5:1:conclu} in Lemma \ref{lem:4.5}, and using the definition of $\hat{\mathcal{F}}_{k+1}$.
	\end{proof}
	
	Note that the proof of Lemma \ref{lem:4.6} is different from that of Lemma \ref{lem:3.6}, since different potential functions, i.e.,  $\hat{\mathcal{F}}_{k+1}$ and $\mathcal{F}_{k+1}$ respectively, are used to establish the convergence of the proposed algorithms. The construction of these potential functions is the key step for our convergence analysis of AGP algorithms.
We are now ready to establish the iteration complexity for the AGP algorithm to achieve an $\varepsilon$-stationary point for solving \eqref{problem:1} under general convex-nonconcave setting.
	
	\begin{thm}\label{thm:4.2}
		Suppose that Assumptions \ref{sec:2:ass1:lip} holds. Let $\{\left(x_k,y_k\right)\}$ be a sequence generated by Algorithm \ref{alg:2.1}
		with parameter settings in \eqref{subsec:4.2:1:para}.
		If  $\bar{\zeta} \le \frac{1}{10L_{11}}$, $\tau>\max\{\frac{19^2(L_{22}+2\bar{\nu}-\bar{\zeta}L_{21}^2)}{20^2\cdot16\bar{\zeta}L_{21}^2}, 2\}$, $q_k=\frac{19}{20\bar{\zeta} k^{\text{1/}4}}, k\geq 1$, then for any given $\varepsilon>0$,
		$$ T(\varepsilon)
		\leq \max \left(	\left( \frac{2\cdot80^2 \bar{\zeta} (\tau -2) L_{21}^2\hat{d}_3\hat{d}_4}{19^2\varepsilon ^2}+2 \right) ^2, \frac{19^4\hat{\sigma} _x^4}{10^4\bar{\zeta} ^4\varepsilon ^4}\right),$$
		where $\hat{\sigma}_x:= \max\{\|x\| \mid x \in  \mathcal{X}\}$,
		$$\hat{d}_3:=\hat{\mathcal{F}}_0-\hat{\mathcal{F}}_2+ \frac{8\cdot 2^{1/4} }{\bar{\zeta}}\hat{\sigma}_x^2 + \frac{19}{40\bar{\zeta}} \hat{\sigma}_x^2,$$
		$\hat{d}_4:=\max\{\hat{D}_1,\frac{10+20\bar{\zeta}^2 L_{21}^2}{9\bar{\zeta}p_2}\}$  with $\hat{\mathcal{D}}_1 := \tfrac{16\tau^2}{(\tau-2)^2} + \tfrac{19^4(\bar{\zeta} L_{21}^2 - \bar{\nu})^2}{16\cdot 20^4(\tau-2)^2 L_{21}^4 \bar{\zeta}^2}$, $\hat{\mathcal{F}}_0 := \bar{f}+\left(2^{13/4}+15+\frac{19}{40}\right)\frac{\hat{\sigma}_x^2}{\bar{\zeta}}$ with $\bar{f} := \max_{(x,y) \in \mathcal{X} \times \mathcal{Y}} f(x,y)$.
	\end{thm}
	
	\begin{proof}
		By $\bar{\zeta}\le \frac{1}{10L_{11}}$,$q_k=\frac{19}{20\bar{\zeta} k^{\text{1/}4}},\forall k\geq 1$, let us denote $\bar{\gamma}_k=\bar{\zeta} L_{21}^2+\frac{16\tau L_{21}^2}{\bar{\zeta} (q_{k+1})^2}-2\bar{\nu}$, $p_k=\frac{8(\tau - 2)L_{21}^2}{\bar{\zeta}( q_{k+1})^2}$, it can be easily checked that the relations in \eqref{lem:4.6:1:para} are satisfied. It follows from the selection of $\bar{\gamma}_k$ and $p_k$ that
		$$\bar{\nu}+\frac{\bar{\gamma}_k}{2}-\frac{\bar{\zeta} L_{21}^2}{2}-\frac{16L_{21}^2}{\bar{\zeta}(q_{k+1})^2}=p_k.$$
		This observation, in view of Lemma \ref{lem:4.6}, then immediately implies that
		\begin{align}\label{thm:4.2:1}
			p_k\|y_{k+1}-y_k \|^2+\frac{9}{10\bar{\zeta}}\| x_{k+1}-x_k \|^2
			\leq& \hat{\mathcal{F}}_{k+1}-\hat{\mathcal{F}}_{k}+\frac{8}{\bar{\zeta}}\left(\frac{q_{k-1}}{q_k} -\frac{q_k}{q_{k+1}} \right)\|x_{k+1}\|^2\nonumber\\
			&+\frac{q_{k-1}-q_k}{2}\|x_{k+1}\|^2.
		\end{align}
		We can easily check from the definition of $f_{k} (x_k, y_k)$ that  $$\|\nabla {G}_k\| -\|\nabla \tilde{G}_k \| \leq q_{k}\| x_k \|.$$
	Similar to \eqref{thm:4.1:3:g} in the proof of Theorem \ref{thm:4.1}, by replacing $f$ with $f_{k}$, $\nabla {G}_k$ with $\nabla \tilde{G}_k$, $\zeta$ with $\bar{\zeta}$, $\nu$ with $\bar{\gamma}_k+\bar{\nu}$ respectively, we conclude that
		\begin{align}\label{thm:4.2:2:g}
			\| \nabla \tilde{G}_k \|^2
			\leq 2\left(\bar{\gamma}_k+\bar{\nu}\right)^2\| y_{k+1}-y_k \|^2+\left( \frac{1}{\bar{\zeta} ^2}+2L_{12}^{2} \right)\|x_{k+1}-x_k \|^2.
		\end{align}
		Since both $p_k$ and $\bar{\gamma}_k$ are in the same order when $k$ becomes large enough, it then follows from the definition of $\hat{D}_1$  that $\forall k\ge 1$,
		\begin{align}\label{thm:4.2:3:D}
			\hat{D}_1\geq \frac{2\left(\bar{\gamma}_k+\bar{\nu}\right)^2}{(p_k)^2}.
		\end{align}
		Combining the previous two inequalities in \eqref{thm:4.2:3:D} and \eqref{thm:4.2:2:g}, we obtain
		\begin{equation}\label{thm:4.2:4:newg}
			\| \nabla \tilde{G}_k \|^2
			\leq \hat{D}_1(p_k)^2\| y_{k+1}-y_k \|^2+\left( \frac{1}{\bar{\zeta}^2}+2L_{12}^{2} \right)\|x_{k+1}-x_k \|^2.
		\end{equation}
		Denote $\hat{d}_k^{(2)}=\frac{1}{\max \left\{ \hat{D}_1p_k,\frac{10+20\bar{\zeta}^2 L_{12}^2}{9\bar{\zeta}} \right\}}$.
		By multiplying $\hat{d}_k^{(2)}$ on the both sides of \eqref{thm:4.2:4:newg}, and using \eqref{thm:4.2:1}, we get
		\begin{align}\label{thm:4.2:5}
			\hat{d}_k^{(2)}\| \nabla \tilde{G}_k \|^2
			\leq \hat{\mathcal{F}}_{k+1}-\hat{\mathcal{F}}_{k}+\frac{8}{\bar{\zeta}}\left(\frac{q_{k-1}}{q_k} -\frac{q_k}{q_{k+1}} \right)\|x_{k+1}\|^2+\frac{q_{k-1}-q_k}{2}\|x_{k+1}\|^2.
		\end{align}
		Denote $\hat{T}(\varepsilon):=\min\{k \mid \lVert \nabla \tilde{G}(x_k,y_k) \rVert \leq \frac{\varepsilon}{2}, k\geq 2\}$. By summing both sides of \eqref{thm:4.2:5} from $k=2$ to $k=\hat{T}(\varepsilon)$, we then obtain
		\begin{align}\label{thm:4.2:6}
			&\sum_{k=2}^{\hat{T}(\varepsilon)}{\hat{d}_k^{(2)}\lVert \nabla \tilde{G}_k \rVert ^2}\nonumber\\
			\leq & \hat{\mathcal{F}}_{\hat{T}(\varepsilon)+1}-\hat{\mathcal{F}}_2+\frac{8}{\bar{\zeta}}\left( \frac{q_1}{q_2}-\frac{q_{\hat{T}(\varepsilon)}}{q_{\hat{T}(\varepsilon)+1} } \right)\hat{\sigma}_x^2+\frac{q_1-q_{\hat{T}(\varepsilon)}}{2}\hat{\sigma}_x^2\nonumber\\
			\leq & \hat{\mathcal{F}}_{\hat{T}(\varepsilon)+1}-\hat{\mathcal{F}}_2+\frac{8q_1}{\bar{\zeta} q_2}\hat{\sigma}_x^2+\frac{q_1}{2}\hat{\sigma}_x^2\nonumber\\
			= & \hat{\mathcal{F}}_{\hat{T}(\varepsilon)+1}-\hat{\mathcal{F}}_2+ \frac{8\cdot 2^{1/4} }{\bar{\zeta}}\hat{\sigma}_x^2 + \frac{19}{40\bar{\zeta}} \hat{\sigma}_x^2.
		\end{align}

Note that by the definition of $\hat{\mathcal{F}}_{k+1}$ in Lemma \ref{lem:4.6}, we have
	\begin{align*}
		\hat{\mathcal{F}}_{\hat{T}\left( \varepsilon \right)+1}&\le
		\bar{f}+\frac{8\hat{\sigma}_x^2}{\bar{\zeta}}\left(1+\frac{1}{\hat{T}\left( \varepsilon \right)}\right)^{1/4}
+\frac{15\hat{\sigma}_x^2}{\bar{\zeta}}+\frac{q_1\hat{\sigma}_x^2}{2}\\
&=\bar{f}+\left(2^{13/4}+15+\frac{19}{40}\right)\frac{\hat{\sigma}_x^2}{\bar{\zeta}}= \hat{\mathcal{F}}_0,
	\end{align*}
where $\bar{f} := \max_{(x,y) \in \mathcal{X} \times \mathcal{Y}} f(x,y)$. We then conclude from \eqref{thm:4.2:6} that
\begin{align}\label{thm:4.2:6n}
\sum_{k=2}^{\hat{T}(\varepsilon)}{\hat{d}_k^{(2)}\lVert \nabla \tilde{G}_k \rVert ^2}\le\hat{\mathcal{F}}_0-\hat{\mathcal{F}}_2+ \frac{8\cdot 2^{1/4} }{\bar{\zeta}}\hat{\sigma}_x^2 + \frac{19}{40\bar{\zeta}} \hat{\sigma}_x^2=\hat{d}_3.
\end{align}

		Note that $\hat{d}_4=\max\{\hat{D}_1,\frac{10+20\bar{\zeta}^2 L_{12}^2}{9\bar{\zeta}p_2}\}$. Since $p_k$ is an increasing sequence when $k\geq 2$, we have that $\hat{d}_4\ge\max\{\hat{D}_1,\frac{10+20\bar{\zeta}^2 L_{12}^2}{9\bar{\zeta}p_k}\}$, which implies that $\hat{d}_k^{(2)}\geq\frac{1}{\hat{d}_4p_k}$. By multiplying $\hat{d}_4$ on the both sides of \eqref{thm:4.2:6n}, and using the definition of $\hat{d}_3$, we have $	\tsum_{k=2}^{\hat{T}(\varepsilon)} \frac{1}{p_k} \lVert \nabla \tilde{G}_k \rVert ^2\leq \hat{d}_3\hat{d}_4$,
		which, by the definition of $\hat{T}(\varepsilon)$, implies that
		\begin{equation}\label{thm:4.2:7}
			\frac{\varepsilon^2}{4} \leq \frac{\hat{d}_3\hat{d}_4}{\tsum_{k=2}^{\hat{T}(\varepsilon)}{\frac{1}{p_k}}}.
		\end{equation}
		Using the assumptions $q_k=\frac{19}{20\bar{\zeta} k^{1/4}}$ and $p_k = \frac{2\cdot40^2\bar{\zeta} (\tau -2) L_{21}^2 \sqrt{k+1}}{19^2}$, \eqref{thm:4.2:7} and the fact $\sum_{k=2}^{\hat{T}(\varepsilon)}1/\sqrt{k+1}\geq \sqrt{\hat{T}(\varepsilon)}-2$ , we conclude that $\frac{\varepsilon^2}{4} \leq \frac{2\cdot40^2\bar{\zeta}(\tau -2) L_{21}^2\hat{d}_3\hat{d}_4}{19^2\left(\sqrt{\hat{T}(\varepsilon)}-2\right)}$ or equivalently,
		\begin{equation*}
			\hat{T}(\varepsilon) \le \left( \frac{2\cdot80^2\bar{\zeta} (\tau -2) L_{21}^2\hat{d}_3\hat{d}_4 }{19^2\varepsilon ^2}+2 \right) ^2.
		\end{equation*}
		On the other hand, if $k> \frac{19^4\hat{\sigma}_x^4}{10^4\bar{\zeta} ^4\varepsilon ^4}$, then $q_k=\frac{19}{20\bar{\zeta} k^{1/4}} \leq \frac{\varepsilon}{2\hat{\sigma}_x}$, this inequality together with the definition of $\hat{\sigma}_x$ then imply that $q_k\lVert x_k \rVert \leq \frac{\varepsilon}{2}$. Therefore, there exists a
		\begin{align*}
			T(\varepsilon) \leq \max (	\hat{T}(\varepsilon), \frac{19^4\hat{\sigma} _x^4}{10^4\bar{\zeta}^4\varepsilon ^4})
			\leq \max \left(	\left( \frac{2\cdot80^2 \bar{\zeta} (\tau -2) L_{21}^2\hat{d}_3\hat{d}_4}{19^2\varepsilon ^2}+2 \right) ^2, \frac{19^4\hat{\sigma} _x^4}{10^4\bar{\zeta}^4\varepsilon ^4}\right),
		\end{align*}
		such that
		$\lVert \nabla G_k\rVert \leq \lVert \nabla \tilde{G}_k \rVert + q_k\lVert x_k \rVert \leq \frac{\varepsilon}{2}+\frac{\varepsilon}{2} =\varepsilon$.
	\end{proof}


	Theorem \ref{thm:4.2} shows that the number of gradient evaluations performed by Algorithm \ref{alg:2.1} to obtain an $\varepsilon$-stationary point of $f$ which satisfies \eqref{lem:4.6:1:para} is bounded by $\mathcal{O}\left(L^4 \varepsilon ^{-4} \right)$ under the general convex-nonconcave setting. To the best of our knowledge, this is the first theoretical guarantee that has been obtained in the literature for solving this class of minimax problems. As mentioned in Section 1, it is difficult to extend existing algorithms, especially those nested-loop methods, for minimax optimization to the general convex-nonconcave setting. One possible approach would be to switch the order of the ``$\min$" and ``$\max$" operators. However, in general, $\min \limits_{x\in \mathcal{X}}\max \limits_{y\in \mathcal{Y}}\ f(x,y) \not= \max \limits_{y\in \mathcal{Y}}	\min \limits_{x\in \mathcal{X}}\ f(x,y)$  when $f(x,y)$ is nonconvex w.r.t. $x$ or nonconcave w.r.t. $y$. Even if the set of stationary points for the above two problems are the same, e.g., by using the same stationarity criterion as in this paper,
	the algorithm applied to these problems will converge to different solutions. This can be seen by applying  the unified AGP algorithm to solve $\min \limits_{x\in \mathcal{X}}\max \limits_{y\in \mathcal{Y}}\ f(x,y) $ and $\max \limits_{y\in \mathcal{Y}}\min \limits_{x\in \mathcal{X}}\ f(x,y)$. The AGP algorithm will likely converge to different stationary points because different potential functions (i.e.,  $\mathcal{F}_{k+1}$ and $\hat{\mathcal{F}}_{k+1}$ in Lemma 3.6 and Lemma 4.6) have to be applied for analyzing the convergence of AGP for solving these two problems.

For  nonconvex-strongly concave or strongly convex-nonconcave setting, the AGP algorithm
with properly chosen parameters is actually equivalent to the alternating GDA algorithm. By constructing suitable potential function, i.e., $F_{k+1}$, we were able to show the  $O(\epsilon^{-2})$ iteration complexity of  this method. Whereas for nonconvex-concave or convex-nonconcave setting,  where $b_k$ or $c_k$ is not equal to $0$, AGP algorithm is completely new and the selection of these parameters $b_k$ or $c_k$ plays a very crucial role to guarantee the convergence of the APG method. In these cases, the key step in our proof is also to construct a suitable potential function, i.e.,  $\mathcal{F}_{k+1}$ and $\hat{\mathcal{F}}_{k+1}$ with $b_k$ or $c_k$ playing a very crucial role  (see Lemma 3.6 and Lemma 4.6 respectively). To our best knowledge, this is the first time that such potential functions have been constructed for solving minimax problems.

\section{Block-wise Nonsmooth Nonconvex Minimax Problems}\label{sec:5}
In this section, we consider a more general block-wise nonsmooth minimax problem as follows.
	\begin{align}
		\min \limits_{x\in \mathcal{X}}\max \limits_{y\in \mathcal{Y}}\   l(x,y):=&f(x^{(1)},x^{(2)},\cdots,x^{(K_1)},y^{(1)},y^{(2)},\cdots,y^{(K_2)})+\sum_{i=1}^{K_1}h_i(x^{(i)})-\sum_{j=1}^{K_2}g_j(y^{(j)}),\tag{BP}\label{problem:block}\\
		\mbox{s.t.}\quad x^{(i)} \in& \mathcal{X}_{i},\quad i=1,\cdots,K_1,\nonumber\\
        \quad y^{(j)} \in& \mathcal{Y}_{j},\quad j=1,\cdots,K_2,\nonumber
	\end{align}
	where $f: \mathbb{R}^{d_{x}K_{1}+d_{y}K_{2}}\rightarrow \mathbb{R}$ is a continuously differentiable function, $h_{i}:\mathbb{R}^{d_{x}}\rightarrow \mathbb{R}$ and $g_{j}: \mathbb{R}^{d_{y}}\rightarrow \mathbb{R}$ are some convex continuous possibly nonsmooth functions, $x=\left[ x^{(1)},x^{(2)},\cdots,x^{(K_1)} \right]\in \mathcal{X}=\mathcal{X}_{1}\times \mathcal{X}_{2}\times \cdots \times\mathcal{X}_{K_{1}}\subset \mathbb{R}^{d_{x}K_{1}}$ and  $y=\left[ y^{(1)},y^{(2)},\cdots,y^{(K_2)} \right]\in \mathcal{Y}=\mathcal{Y}_{1}\times \mathcal{Y}_{2}\times \cdots \times\mathcal{Y}_{K_{2}}\subset \mathbb{R}^{d_{y}K_{2}}$ are the block variables, $\mathcal{X}_{i}(i=1,\cdots,K_{1})$ and $\mathcal{Y}_{j}~(j=1,\cdots,K_{2})$ are nonempty compact convex sets. For notational simplicity, we denote $f\left( x^{(1)},x^{(2)},\cdots,x^{(K_1)},y^{(1)},y^{(2)},\cdots,y^{(K_2)}\right)$ by $f\left( x,y\right)$.

Note that \eqref{problem:block} reduces to \eqref{problem:1} if we set $K_1=K_2=1$, $h_{i}(\cdot)=0$ and $g_{j}(\cdot)=0$,
	which means that  \eqref{problem:1}  is a special case of \eqref{problem:block}. There exist very few known existing nested-loop algorithms that are proposed to solve \eqref{problem:block} with multi-block structure. However, the minimax problems with block structure are important in machine learning and signal processing, e.g., distributed training \cite{Lu}.
	
	The difficulty to solve \eqref{problem:block} comes from two aspects. One is that for any given $j$, $f(x,y^{(1)},y^{(2)},\cdots,y^{(K_2)})$ is nonconcave w.r.t. $y^{(j)}$ for any given $x$ and other blocks of $y$, i.e., $y^{(1)},y^{(2)},\cdots, y^{(j-1)}, y^{(j+1)}$, $\cdots$, $y^{(K_2-1)}$ and $y^{(K_2)}$. For any given $x$, to solve the inner maximization subproblems with respect to $y^{(j)}$ is already NP-hard.
	Due to this reason, almost all the existing nested-loop algorithms will lose their theoretic guarantees since they need to solve the inner subproblem exactly, or approximately with an error proportional to the accuracy $\varepsilon$, unless we can exchange the order of ``$\min$" and ``$\max$" operators. However, in general, $\min \limits_{x\in \mathcal{X}}\max \limits_{y\in \mathcal{Y}}\ l(x,y) \not= \max \limits_{y\in \mathcal{Y}}	\min \limits_{x\in \mathcal{X}}\ l(x,y)$  when $f(x,y)$ is nonconvex w.r.t. $x$ or nonconcave w.r.t. $y$. Although $\min \limits_{x\in \mathcal{X}}\max \limits_{y\in \mathcal{Y}}\ l(x,y) $ and $\max \limits_{y\in \mathcal{Y}}	\min \limits_{x\in \mathcal{X}}\ l(x,y)$ share the same stationary point set, when we use the same algorithm to solve them, it may converge to different stationary points.
	Another difficulty comes from  the multi-block structure. To the best of our knowledge, there are very few  nested-loop algorithms (and complexity results) for solving the aforementioned multi-block structure nonconvex minimax problems before our work. On the other hand, single loop algorithms do not need to solve the inner subproblem, and can be easily generalized to handle the multi-block structure.
	
	Similar to the idea of AGP algorithm, we propose a block alternating proximal gradient algorithm (BAPG) to solve \eqref{problem:block}. Instead of the original function $l\left(x,y \right) $, the BAPG algorithm uses the gradient of a regularized version of the original function, i.e.,
	\begin{equation}\label{sec5:1}
	 	\bar{l}_k\left(x,y \right):=  \bar{f}_k\left(x,y \right)+\sum_{i=1}^{K_1}h_i(x^{(i)})-\sum_{j=1}^{K_2}g_j(y^{(j)}),
	\end{equation}
	where $ \bar{f}_k\left(x,y \right) =f\left(x,y \right)+\frac{\hat{b}_k}{2}\sum_{i=1}^{K_1}\|x^{(i)}\|^2 - \frac{\hat{c}_k}{2}\sum_{j=1}^{K_2}\|y^{(j)}\|^2 $ with regularization parameters $\hat{b}_k\geq 0$ and $\hat{c}_k\geq 0$ at the $k$th iteration. Before presenting the detailed algorithm, we give some notations as follows.	We denote two proximity operators for $x^{(i)}$ and $y^{(j)}$ as follows,
	\begin{align}
		\operatorname{Prox}_{h_i,\mathcal{X}_i}^{\beta^{(i)}}(\upsilon^{(i)})&:=\arg\min\limits_{x^{(i)}\in \mathcal{X}_i} h_i(x^{(i)})+\frac{\beta^{(i)}}{2}\|x^{(i)}-\upsilon^{(i)} \|^2,\label{proximal-x}\\
		\operatorname{Prox}_{g_j,\mathcal{Y}_j}^{\xi^{(j)}}(\omega^{(j)})&:=\arg\max\limits_{y^{(j)}\in \mathcal{Y}_j}-g_j(y^{(j)}) -\frac{\xi^{(j)}}{2}\| y^{(j)}-\omega^{(j)} \|^2.\label{proximal-y}
	\end{align}
	Let $k$ be the number of iteration. Denote $x_k=\left[ x_k^{(1)},x_k^{(2)},\cdots,x_k^{(K_1)} \right]$, $y_k=\left[ y_k^{(1)},y_k^{(2)},\cdots,y_k^{(K_2)} \right]$ and define
	\begin{align*}
		&v^{(i)}_{k+1}:=\left[x^{(1)}_{k+1}, x^{(2)}_{k+1}, \cdots, x^{(i-1)}_{k+1}, x^{(i)}_{k}, \cdots x^{(K_1)}_{k}\right],\\
		&v^{(-i)}_{k+1}:=\left[x^{(1)}_{k+1}, x^{(2)}_{k+1}, \cdots, x^{(i-1)}_{k+1}, x^{(i+1)}_{k}, \cdots x^{(K_1)}_{k}\right],\\
	&w^{(j)}_{k+1}:=\left[y^{(1)}_{k+1}, y^{(2)}_{k+1}, \cdots, y^{(j-1)}_{k+1}, y^{(j)}_{k}, \cdots y^{(K_2)}_{k}\right],\\
	&w^{(-j)}_{k+1}:=\left[y^{(1)}_{k+1}, y^{(2)}_{k+1}, \cdots, y^{(j-1)}_{k+1}, y^{(j+1)}_{k}, \cdots y^{(K_2)}_{k}\right].
	\end{align*}

Each iteration of the proposed BAPG algorithm conducts two proximal gradient steps for updating both $x$ and $y$. More specifically, at the $k$th iteration, it updates $x_k$ by minimizing a linearized approximation of $\bar{l}_k\left( x,  y_k\right)$ with the gradient at point $\left( v_{k+1}^{(i)},y_k \right)$, i.e., for each $i=1,\cdots,K_1$,
	\begin{align*}
		x^{(i)}_{k+1}=&\arg\min\limits_{x^{(i)}\in \mathcal{X}_i}\left\langle \nabla _{x^{(i)}}\bar{f}_k(v^{(i)}_{k+1},y_k) ,x^{(i)}-x^{(i)}_{k} \right\rangle +h_i(x^{(i)})+\frac{\beta^{(i)}_k}{2}\|x^{(i)}-x^{(i)}_{k} \|^2\nonumber\\
		=& \operatorname{Prox}_{h_i,\mathcal{X}_i}^{\beta^{(i)}_k}\left( x^{(i)}_{k}-\frac{1}{\beta^{(i)}_k}\nabla _{x^{(i)}} f( v^{(i)}_{k+1},y_k )-\frac{1}{\beta^{(i)}_k}\hat{b}_kx^{(i)}_{k}\right),
	\end{align*}
	where  $\operatorname{Prox}_{h_i,\mathcal{X}_i}^{\beta^{(i)}_k}$ is the proximal operator which is defined in \eqref{proximal-x} and $\beta^{(i)}_k > 0$. Similarly, it updates $y_k$ by maximizing a linearized approximation of $\bar{l}_k\left( x_{k+1},y \right)$ minus some regularized terms, i.e.,
	\begin{align*}
		y^{(j)}_{k+1}=&\arg\max\limits_{y^{(j)}\in \mathcal{Y}_j}\langle \nabla _{y^{(j)}} \bar{f}_k( x_{k+1},w^{(j)}_{k+1}) ,y^{(j)}-y^{(j)}_k \rangle-g_j(y^{(j)}) -\frac{\xi^{(j)}_k}{2}\| y^{(j)}-y^{(j)}_k \|^2\nonumber\\
		=&\operatorname{Prox}_{g_j,\mathcal{Y}_j}^{\xi^{(j)}_k}\left( y^{(j)}_k+ \frac{1}{\xi^{(j)}_k} \nabla _{y^{(j)}}f(x_{k+1},w^{(j)}_{k+1})-\frac{1}{\xi^{(j)}_k}\hat{c}_ky^{(j)}_k \right),
	\end{align*}
	where  $\operatorname{Prox}_{g_j,\mathcal{Y}_j}^{\xi^{(j)}_k}$ is the proximal operator which is defined in \eqref{proximal-y} and $\xi^{(j)}_k >0$.  The proposed BAPG algorithm is formally stated in Algorithm \ref{alg:BAPG}, where sequences $\beta^{(i)}_k$, $\hat{b}_k$, $\xi^{(j)}_k$, $\hat{c}_k$ and the stopping rule in Step 4 will be specified later in each of the different problem settings to be studied. Note that it reduces to the AGP algorithm when $K_1=K_2=1$, $h_{i}(\cdot)=0$ and $g_{j}(\cdot)=0$.
	\begin{algorithm}
		\caption{(BAPG Algorithm)}
		\label{alg:BAPG}
		\begin{algorithmic}
			\STATE{\textbf{Step 1}:  Input: $x_{1}, y_{1}$, $\hat{b}_1$, $\hat{c}_1$, $\beta^{(i)}_1$ $(i=1,\cdots,K_1)$,  $\xi^{(j)}_1$ $(j=1,\cdots,K_2)$; Set $k=1$.}
			\STATE{\textbf{Step 2}:~Calculate $\hat{b}_k$, $\beta^{(i)}_k$, and for $i=1,\ldots, K_1$, perform the following update for $x^{(i)}_k$:  	\begin{align}\label{BAPG:upx}
					x^{(i)}_{k+1}
					&=\operatorname{Prox}_{h_i,\mathcal{X}_i}^{\beta^{(i)}_k}\left( x^{(i)}_{k}-\frac{1}{\beta^{(i)}_k}\nabla _{x^{(i)}} f( v^{(i)}_{k+1},y_k )-\frac{1}{\beta^{(i)}_k}\hat{b}_kx^{(i)}_{k}\right).
			\end{align}}		
			\STATE{\textbf{Step 3}:~Calculate $\hat{c}_k$, $\xi^{(j)}_k$, and for $j=1,\ldots, K_2$ and perform the following update for $y^{(j)}_k$:
				\begin{align}\label{BAPG:upy}
					y_{k+1}
					&=\operatorname{Prox}_{g_j,\mathcal{Y}_j}^{\xi^{(j)}_k}\left( y^{(j)}_k+ \frac{1}{\xi^{(j)}_k} \nabla _{y^{(j)}}f(x_{k+1},w^{(j)}_{k+1})-\frac{1}{\xi^{(j)}_k}\hat{c}_ky^{(j)}_k \right).
			\end{align}}
			\STATE{\textbf{Step 4}:~If some stationary condition is satisfied, stop; otherwise, set $k=k+1, $ go to Step 2.}
		\end{algorithmic}
	\end{algorithm}

Before analyzing the convergence of Algorithm \ref{alg:BAPG}, we define the stationarity gap as the termination criterion as follows.
	
	\begin{defi}\label{nblock:G}
		At each iteration of Algorithm \ref{alg:BAPG}, the stationarity gap for problem \eqref{problem:block} w.r.t. $l(x,y)$ is defined as:
		\begin{equation*}
			\nabla \mathcal{G}\left( x_k,y_k \right) :=\left[\begin{array}{c}
				{\beta^{(1)}_k\left(x^{(1)}_k-\operatorname{Prox}_{h_1,\mathcal{X}_1}^{\beta^{(1)}_k}(x^{(1)}_k-\frac{1}{\beta^{(1)}_k} \nabla_{x^{(1)}} f(x_k ,y_k))\right)} \\
				{\vdots} \\
				{\beta^{(K_1)}_k\left(x^{(K_1)}_k-\operatorname{Prox}_{h_{K_1},\mathcal{X}_{K_1}}^{\beta^{(K_1)}_k}(x^{(K_1)}_k-\frac{1}{\beta^{(K_1)}_k} \nabla_{x^{(K_1)}} f(x_k ,y_k))\right)}\\
				{\xi^{(1)}_k\left(y^{(1)}_k-\operatorname{Prox}_{g_1, \mathcal{Y}_1}^{\xi^{(1)}_k}\left(y^{(1)}_k+\frac{1}{\xi^{(1)}_k} \nabla_{y^{(1)}} f(x_k ,y_k)\right)\right)} \\
{\vdots} \\
{\xi^{(K_2)}_k\left(y^{(K_2)}_k-\operatorname{Prox}_{g_{K_2},\mathcal{Y}_{K_2}}^{\xi^{(K_2)}_k}\left(y^{(K_2)}_k
+\frac{1}{\xi^{(K_2)}_k} \nabla_{y^{(K_2)}} f(x_k ,y_k)\right)\right)}
			\end{array}\right].
		\end{equation*}
		For simplicity, we denote $\nabla {\mathcal{G}}_k:= \nabla {\mathcal{G}}(x_k, y_k)$, $(\nabla \mathcal{G}_k)_{x^{(i)}}:= \beta^{(i)}_k( x^{(i)}_k-\operatorname{Prox}_{h_i,\mathcal{X}_i}^{\beta^{(i)}_k}( x^{(i)}_k- \frac{1}{\beta^{(i)}_k} \nabla _{x^{(i)}}f( x_k,y_k)) )$ and $	(\nabla \mathcal{G}_k)_{y^{(j)}}:=\xi^{(j)}_k\left( y^{(j)}_k-\operatorname{Prox}_{g_j,\mathcal{Y}_j}^{\xi^{(j)}_k}( y^{(j)}_k+ \frac{1}{\xi^{(j)}_k} \nabla _{y^{(j)}}f\left( x_k,y_k \right)) \right)$.
	\end{defi}

	\begin{defi}\label{nblock:Gn}
		At each iteration of Algorithm \ref{alg:BAPG}, the stationarity gap for problem \eqref{problem:block} w.r.t. $\bar{l}_k(x,y)$, denoted by $\nabla \bar{\mathcal{G}}\left( x_k,y_k \right) $, is defined almost the same as in Definition \ref{nblock:G} except by replacing $f$ with $\bar{f}_k$.
		For simplicity, we denote $\nabla {\bar{\mathcal{G}}}_k:=\nabla {\bar{\mathcal{G}}}(x_k, y_k)$, $(\nabla \bar{\mathcal{G}}_k)_{x^{(i)}}:= \beta^{(i)}_k( x^{(i)}_k-\operatorname{Prox}_{h_i,\mathcal{X}_i}^{\beta^{(i)}_k}( x^{(i)}_k- \frac{1}{\beta^{(i)}_k} \nabla _{x^{(i)}}\bar{f}_k\left( x_k,y_k \right)))$ and $(\nabla \bar{\mathcal{G}}_k)_{y^{(j)}}:=\xi^{(j)}_k\left( y^{(j)}_k-\operatorname{Prox}_{g_j,\mathcal{Y}_j}^{\xi^{(j)}_k}( y^{(j)}_k+ \frac{1}{\xi^{(j)}_k} \nabla _{y^{(j)}}\bar{f}_k\left( x_k,y_k \right)) \right)$.
	\end{defi}

%
%
%
%
%

\subsection{Complexity Analysis for Nonconvex-(Strongly) Concave Setting}
\subsubsection{Nonconvex-Strongly Concave Setting}
In this subsection, we analyze the iteration complexity of Algorithm \ref{alg:BAPG} for solving the
nonsmooth one-sided block-wise nonconvex-strongly concave minimax optimization problem \eqref{problem:block} with $K_2=1$, i.e., $f(x,y)$ is nonconvex w.r.t. $x$ for any fixed $y\in \mathcal{Y}$, and $\hat{\mu}$-strongly concave  w.r.t. $y$ for any given $x\in \mathcal{X}$. Under this setting,  $\forall k\ge 1$, we set
	\begin{equation}\label{sec5.1:1}
		\beta^{(i)}_k=\hat{\eta}, \xi^{(j)}_k=\tfrac{1}{\hat{\rho}}, \ \mbox{and} \ \hat{b}_k=\hat{c}_k=0.
	\end{equation}
Lemma \ref{nlem:5.1} below shows a descent result for the $x_k$ update.
	\begin{lem}\label{nlem:5.1}
		Suppose that Assumption \ref{sec:2:ass1:lip} holds. Let $\{(x_k,y_k)\}$ be a sequence generated by Algorithm \ref{alg:BAPG} with parameter settings in \eqref{sec5.1:1}. If $\forall k,\hat{\eta} > L_{11}$,  we have
		\begin{equation}\label{nlem:5.1:1}
			l(x_{k+1},y_{k})-l(x_k,y_k)\le -\frac{\hat{\eta}}{2}\| x_{k+1}-x_k \|^2.
		\end{equation}
	\end{lem}
	
	\begin{proof}
		By the optimality condition for \eqref{BAPG:upx}, $\forall i$ we have
		\begin{equation}\label{nlem:5.1:2}
			\langle  \nabla _{x^{(i)}}f( v^{(i)}_{k+1},y_k)+\partial h_i(x^{(i)}_{k+1})+\hat{\eta}( x^{(i)}_{k+1}-x^{(i)}_k) ,x^{(i)}_k-x^{(i)}_{k+1} \rangle  \ge 0.
		\end{equation}
		By Assumption \ref{sec:2:ass1:lip} and the convexity of $h_i(x)$, we obtain
		\begin{align}\label{nlem:5.1:3}
			&l(x^{(i)}_{k+1},v^{(-i)}_{k+1},y_k)-l( x^{(i)}_k,v^{(-i)}_{k+1},y_k) \nonumber\\
			\le & \langle \nabla _{x^{(i)}}f(v^{(i)}_{k+1},y_k) +\partial h_i(x^{(i)}_{k+1}),x^{(i)}_{k+1}
			-x^{(i)}_k \rangle +\frac{L_{11}}{2}\| x^{(i)}_{k+1}-x^{(i)}_k \|^2.
		\end{align}
		Adding \eqref{nlem:5.1:2} and \eqref{nlem:5.1:3}, and summing up the inequality from $i=1$ to $i=K_1$, we have
		\begin{equation*}
			l( x_{k+1},y_{k}) -l\left( x_k,y_k \right)\le -\sum_{i=1}^{K_1}\left(\hat{\eta} -\frac{L_{11}}{2} \right) \| x^{(i)}_{k+1}-x^{(i)}_k \|^2.
		\end{equation*}
		By using the assumption that $\hat{\eta} > L_{11}$, we complete the proof.
	\end{proof}

The rest of the proof is almost the same as that of the AGP algorithm shown in Section \ref{sec:3} since the ascent step of $y$'s update can be similarly estimated.
By replacing $\nabla _yf(x_{k+1},y_k)$ with $\nabla _yf(x_{k+1},y_k)-\partial g_1(y_{k+1})$, and using $\langle \partial g_1(y_{k+1})-\partial g_1(y_k),y_{k+1}-y_k \rangle \geq0$, and the fact that  $\| x^{(i)}_{k+1}-x^{(i)}_k \|\le\| x_{k+1}-x_k \|$ and $\|v^{(i)}_{k+1}-x_k \|\le\| x_{k+1}-x_k \|$, we can show similar results to those in Lemmas \ref{lem:3.2:blostr:y} and \ref{lem:3.3:blostr:F}. Then, similar to the proof of Theorem \ref{thm:3.1:blockstr}, we obtain the following convergence result by some parameters replacement, e.g., $\eta$ with $\hat{\eta}+L_{11}$ in \eqref{thm:3.1:1:gx} and $\eta^2$ with $K_1\left(\hat{\eta}+L_{11}\right)^2$ in \eqref{thm:3.1:3:g}. We omit the proof details here for simplicity.

	In particular,
letting $\nabla \mathcal{G}\left( x_k,y_k \right)$ be defined as  in Definition \ref{nblock:G} and $\varepsilon>0$ be
a given target accuracy, we provide
a bound on $\mathcal{T}(\varepsilon)$, the first iteration index to achieve an $\varepsilon$-stationary point, i.e., $ \| \nabla \mathcal{G}(x_k,y_k)\| \leq \varepsilon$, which is equivalent to
\begin{equation}\label{defi:T}
	\mathcal{T}(\varepsilon):=\min\{k\mid \|\nabla \mathcal{G}(x_k,y_k)\|\leq \varepsilon \}.
\end{equation}
\begin{thm}\label{nthm:5.1}
		Suppose that Assumption \ref{sec:2:ass1:lip} holds. Let $\{\left(x_k,y_k\right)\}$ be a sequence generated by Algorithm \ref{alg:BAPG}
		with parameter settings in \eqref{sec5.1:1}. If $$
			\hat{\eta}>L_{11}, \hat{\eta}>L_{12}^2\hat{\rho}+\frac{4L_{12}^2}{\hat{\rho} \hat{\mu}^2}, \ \mbox{and} \ \hat{\rho}\leq \frac{\hat{\mu}}{4L_{22}^2},
$$
		then it holds that $$ \mathcal{T}\left( \varepsilon \right) \le \frac{L_1-\underline{L}}{r_1\varepsilon ^2},$$
		where $	r_1 :=\frac{\min \left\{\frac{\hat{\eta}}{2}-\frac{\hat{\rho} L_{12}^2}{2}-\frac{2L_{12}^2}{\hat{\rho} \hat{\mu}^2},\frac{3\hat{\mu}-\hat{\rho} L_{22}^2}{2}+\frac{\hat{\mu}-4\hat{\rho} L_{22}^2}{2\hat{\rho}\hat{\mu}} \right\}}{\max \left\{ K_1\left(\hat{\eta}+L_{11}\right)^2+2L_{12}^2,\frac{2}{\hat{\rho}^2} \right\}}$, $L_1:=l(x_1,y_1)+\frac{2\delta_y^2}{\hat{\rho}^2\hat{\mu}}$ and  $ 			\underline{L}:=\underline{l}-(\hat{\mu}+\frac{7}{2\hat{\rho}}-\frac{\hat{\rho} L_{22}^2}{2}-\frac{2L_{22}^2}{\hat{\mu}})\delta_y^2$ with $\underline{l} := \min_{(x,y) \in \mathcal{X} \times \mathcal{Y}} l(x,y)$ and $\delta _y:= \max\{\|y_1-y_2\| \mid y_1,y_2 \in  \mathcal{Y}\}$.
	\end{thm}

	Theorem \ref{nthm:5.1} implies that the iteration complexity of Algorithm \ref{alg:BAPG} to obtain an $\varepsilon$-stationary point for solving general block-wise nonsmooth nonconvex-strongly concave minimax problems \eqref{problem:block} is bounded by $\mathcal{O}\left( L^2\varepsilon^{-2} \right)$.
%
\subsubsection{Nonconvex-Concave Setting}
We analyze the iteration complexity of Algorithm \ref{alg:BAPG} for solving \eqref{problem:block} in the nonconvex-concave setting. Under this setting, let $K_2=1$, and $\forall k\ge 1$,  we set
	\begin{equation}\label{sec5.1.2}
		  \beta^{(i)}_k=\iota+\vartheta_k,\quad \xi^{(j)}_k=\frac{1}{\xi},\ \mbox{and} \ \hat{b}_k=0,
	\end{equation}
	where $\vartheta_k$ is stepsize parameter to be defined later.

	\begin{ass}\label{ass5.2}
		$\{\hat{c}_k\}$ is a nonnegative monotonically decreasing sequence.
	\end{ass}

	\begin{lem}\label{nlem:5.2}
		Suppose that Assumption \ref{sec:2:ass1:lip} holds. Let $\{(x_k,y_k)\}$ be a sequence generated by Algorithm \ref{alg:BAPG} with parameter settings in \eqref{sec5.1.2}. If $\forall k,\vartheta_k>L_{11}$,  we have
		\begin{equation}\label{nlem:5.2:1}
			l(x_{k+1},y_k)-l(x_k,y_k)\le -\left( \iota +\frac{\vartheta_k}{2} \right) \| x_{k+1}-x_k \|^2.
		\end{equation}
	\end{lem}
	\begin{proof}
 The proof is the same with that of Lemma~\ref{nlem:5.1} except replacing $\beta^{(i)}_k$ by $\iota +\vartheta_k$. We omit the details here.
\end{proof}

The rest of the proof is almost the same as that of the AGP algorithm shown in Section \ref{sec:3} since the ascent step of $y$'s update can be similarly estimated.
By replacing  $\nabla _yf_k(x_{k+1},y_k)$ with $\nabla _y\bar{f}_k(x_{k+1},y_k)-\partial g_1(y_{k+1})$, and using $\langle \partial g_1(y_{k+1})-\partial g_1(y_k),y_{k+1}-y_k \rangle \geq0$, and the fact that  $\| x^{(i)}_{k+1}-x^{(i)}_k \|\le\| x_{k+1}-x_k \|$ and $\|v^{(i)}_{k+1}-x_k \|\le\| x_{k+1}-x_k \|$, we can prove similar results to those in Lemmas \ref{lem:3.5} and \ref{lem:3.6}. Then, similar to the proof of Theorem \ref{thm3.2}, we obtain the following result only through some parameters replacement, e.g., $\bar{\eta}+\bar{\beta}_k$ with $\iota+\vartheta_k+L_{11}$ in \eqref{thm3.2:2:gx} and $(\bar{\eta}+\bar{\beta}_k)^2$ with $K_1\left(\iota+\vartheta_k+L_{11}\right)^2$ in \eqref{thm3.2:4:g}. 

\begin{thm}\label{nthm:5.2}
		Suppose that Assumption \ref{sec:2:ass1:lip} holds. Let $\{(x_k,y_k)\}$ be a sequence generated by Algorithm \ref{alg:BAPG} with parameter settings in \eqref{sec5.1.2}.
		If $\xi\le \frac{1}{10L_{22}}$, $\hat{c}_k=\frac{19}{20\xi k^{\text{1/}4}}, k\geq 1$, $\tau>\max\{\frac{\xi\hat{c}_1^2(L_{11}+2\iota-\xi L_{12}^2)}{16L_{12}^2}, 2\}$, then for any given $\varepsilon>0$,
		$$\mathcal{T}(\varepsilon)
		\leq \max \left(\left( \frac{2\cdot80^2\xi (\tau -2) L_{12}^2r_2r_3 }{19^2\varepsilon ^2}+1 \right) ^2, \frac{19^4\hat{\delta}_y^4}{10^4\xi ^4\varepsilon ^4}\right),$$
		where $r_2=\mathcal{L}_2-\underline{\mathcal{L}}+ (8\cdot 2^{1/4} + \frac{19}{40} )\frac{\hat{\delta}_y^2}{\xi}$, $r_1=\frac{8K_1\tau^2}{(\tau-2)^2} + \frac{19^4\left(2K_1\left(\xi L_{12}^2-\iota+L_{11}\right)^2 + 2L_{12}^2\right)}{64\cdot20^4 \xi^2 (\tau-2)^2 L_{12}^4}$, 	$r_3=\max\{r_1,\frac{19^2}{1440 \sqrt{2}(\tau -2 )\xi^2 L_{12}^2}\}$, $\hat{\delta}_y:= \max\{\|y\| \mid y \in  \mathcal{Y}\}$, $\mathcal{L}_2:=l(x_2,y_2)+\frac{16\hat{\delta}_y^2}{\xi^2\hat{c}_2}$,
		$\underline{\mathcal{L}} := \underline{l}-\left(2^{13/4}+15+\frac{19}{40}\right)\frac{\hat{\delta}_y^2}{\xi}$ with $\underline{l} := \min_{(x,y) \in \mathcal{X} \times \mathcal{Y}} l(x,y)$.
	\end{thm}

	By setting $\xi = \frac{1}{10L_{22}}$, from Theorem \ref{nthm:5.2}, we conclude that for any given $\varepsilon \in (0,1)$,
	\[\mathcal{T}(\varepsilon)\leq \max \left(\left( \frac{1280(\tau -2) L_{12}^2r_2r_3 }{19^2L_{22}\varepsilon ^2}+1 \right) ^2, \frac{19^4\hat{\delta}_y^4L_{22}^4}{\varepsilon ^4}\right)= \mathcal{O}\left( L^4\varepsilon^{-4} \right).\] This implies that the iteration complexity of the proposed BAPG algorithm to obtain a point that satisfies $\|\nabla \mathcal{G}(x_k, y_k)\|\leq \varepsilon$ with $\nabla \mathcal{G}(x_k, y_k)$ being defined in Definition \ref{nblock:G}  for nonsmooth block-wise nonconvex-concave minimax problems \eqref{problem:block} is bounded by $\mathcal{O}\left( L^4\varepsilon^{-4} \right)$.

%

\subsection{Complexity Analysis for (Strongly)  Convex-Nonconcave Setting}
In this subsection we analyze the iteration complexity of Algorithm \ref{alg:BAPG} for solving \eqref{problem:block} with $K_1=1$ under the convex-nonconcave setting.
We first consider the strongly convex-nonconcave setting. Under this setting,   $\forall k\ge 1$, we set
\begin{equation}
\beta^{(i)}_k=\tfrac{1}{\varrho}, \xi^{(j)}_k=\bar{\xi},  \hat{b}_k=\hat{c}_k=0,\label{sec5.2:par}
\end{equation}

\begin{lem}\label{nlem:5.3}
		Suppose that Assumption \ref{sec:2:ass1:lip} holds. Let $\{(x_k,y_k)\}$ be a sequence generated by Algorithm \ref{alg:BAPG} with parameter settings in \eqref{sec5.2:par}. If $\forall k,\bar{\xi}> L_{22}$, then  we have
		\begin{equation}\label{nlem:5.3:1}
			l( x_{k+1},y_{k+1}) -l( x_{k+1},y_k )\ge \frac{\bar{\xi}}{2}\| y_{k+1}-y_k \|^2.
		\end{equation}
	\end{lem}
	
	\begin{proof}
		By the optimality condition for \eqref{BAPG:upy}, $\forall j$, we have
		\begin{equation}\label{nlem:blocvx:y:optx}
			\langle  \nabla _{y^{(j)}}f( x_{k+1},w^{(j)}_{k+1})-\partial g_j(y^{(j)}_{k+1})-\bar{\xi}( y^{(j)}_{k+1}-y^{(j)}_k) ,y^{(j)}_k-y^{(j)}_{k+1} \rangle  \leq 0.
		\end{equation}
		By  the convexity of $g_j(\cdot)$ and Assumption \ref{sec:2:ass1:lip}, the gradient of $f$ is Lipschitz continuous,  implying that
		\begin{align}\label{nlem:blocvx:y:des}
			&l(x_{k+1},w^{(-j)}_{k+1},y^{(j)}_{k+1})-l( x_{k+1},w^{(-j)}_{k+1},y^{(j)}_{k})\nonumber\\
			\ge & \langle \nabla _{y^{(j)}}f(x_{k+1},w^{(j)}_{k+1}) -\partial g_j(y^{(j)}_{k+1}),y^{(j)}_{k+1}
    -y^{(j)}_k \rangle -\frac{L_{22}}{2}\| y^{(j)}_{k+1}-y^{(j)}_k \|^2.
		\end{align}
		Adding \eqref{nlem:blocvx:y:optx} and \eqref{nlem:blocvx:y:des}, and then summing it up from $j=1$ to $K_2$, we have
		\begin{equation*}
			l( x_{k+1},y_{k+1}) -l( x_{k+1},y_k )\ge \sum_{j=1}^{K_2}( \bar{\xi} -\frac{L_{22}}{2} ) \| y^{(j)}_{k+1}-y^{(j)}_k \|^2.
		\end{equation*}
	The result then follows by the assumption that $\bar{\xi}> L_{22}$.
	\end{proof}

By Lemma \ref{nlem:5.3}, similar to the proof of Lemma \ref{lem:4.2}-\ref{lem:4.3} and Theorem \ref{thm:4.1} in Subsection \ref{subsec4}, we can prove the following theorem and we omit the details here.
\begin{thm}\label{nthm:5.3}
	Suppose that Assumption \ref{sec:2:ass1:lip} holds. Let $\{\left(x_k,y_k\right)\}$ be a sequence generated by Algorithm \ref{alg:BAPG}
	with parameter settings in \eqref{sec5.2:par}. If $$
			\bar{\xi}>L_{22}, \bar{\xi}>L_{21}^2\varrho+\frac{4L_{21}^2}{\varrho \hat{\theta}^2}, \ \mbox{and} \ \varrho\leq \frac{\hat{\theta}}{4L_{11}^2},
$$
	then $\forall \varepsilon>0$,  it holds that $$\mathcal{ T}\left( \varepsilon \right) \le
	\tfrac{\overline{L}-\hat{L}_1}{\hat{r}_1\varepsilon ^2},$$
	where
	$\hat{r}_1:=\tfrac{\min \left\{\tfrac{\bar{\xi}}{2}-\tfrac{\varrho L_{21}^2}{2}-\tfrac{2 L_{21}^2}{\varrho \hat{\theta}^2},\tfrac{3\hat{\theta}-\varrho L_{11}^2}{2}+\tfrac{\hat{\theta}-4\varrho L_{11}^2}{2\varrho\hat{\theta}} \right\}}{\max \left\{ \tfrac{1}{\varrho^2}+2K_2L_{12}^{2},2K_2(\bar{\xi}+L_{22})^2 \right\}}
	$
	and $\overline{L}=\hat{l}+(\hat{\theta}+\frac{7}{2\varrho}-\frac{\varrho L_{11}^2}{2}-\frac{2L_{11}^2}{\hat{\theta}})\delta_x^2$, $\hat{L}_1:=\hat{l}-\frac{4\delta _x^2}{\varrho^2\hat{\theta}}-(\frac{L_{21}\varrho}{2}+\frac{2L_{21}^2}{\hat{\theta}^2\varrho})\|y_{2}-y_1\|^2$ with $\hat{l} := \max\limits_{(x,y) \in \mathcal{X} \times \mathcal{Y}} l(x,y)$ and $\delta _x= \max\{\|x_1-x_2\| \mid \forall x_1, x_2\in \mathcal{X}\}$.
\end{thm}

Theorem \ref{nthm:5.3} implies that the iteration complexity of the proposed BAPG algorithm to obtain an $\varepsilon$-stationary point for smooth strongly convex-nonconcave minimax problems \eqref{problem:block} is bounded by $\mathcal{O}\left(L^2 \varepsilon^{-2} \right)$.


We can also analyze the iteration complexity of Algorithm \ref{alg:BAPG} applied to
	the general convex-nonconcave setting. 	Under this setting,  $\forall k\ge 1$, we set
	\begin{equation}
		 \beta^{(i)}_k=\frac{1}{\bar{\varrho}},\quad \xi^{(j)}_k=\bar{\iota}+\bar{\vartheta}_k,\quad \hat{c}_k=0,\label{subsec:5.2:6:para}
	\end{equation}
	where $\xi^{(j)}_k$ is stepsize parameter to be defined later.

\begin{ass}\label{ass5.4}
$\{\hat{b}_k\}$ is a nonnegative monotonically decreasing sequence.
\end{ass}

Similar to the proof of Lemma \ref{lem:4.5}-\ref{lem:4.6} and Theorem \ref{thm:4.2} in Subsection \ref{subsec4-s}, we can prove the following theorem.
\begin{thm}\label{nthm:5.4}
		Suppose that Assumption \ref{sec:2:ass1:lip} holds. Let $\{\left(x_k,y_k\right)\}$ be a sequence generated by Algorithm \ref{alg:BAPG}
		with parameter settings in \eqref{subsec:5.2:6:para}.
		If  $\bar{\varrho} \le \frac{1}{10L_{11}}$,$\hat{b}_k=\frac{19}{20\bar{\varrho} k^{\text{1/}4}}, k\geq 1$, $\tau>\max\{\frac{\bar{\varrho}\hat{b}_2^2(L_{22}+2\bar{\iota}-\bar{\varrho}L_{21}^2)}{16L_{21}^2}, 2\}$, then for any given $\varepsilon>0$,
		$$ \mathcal{T}(\varepsilon)
		\leq \max \left(	\left( \frac{2\cdot80^2 \bar{\varrho} (\tau -2) L_{21}^2\hat{r}_2\hat{r}_3}{19^2\varepsilon ^2}+2 \right) ^2, \frac{19^4\hat{\delta} _x^4}{10^4\bar{\varrho} ^4\varepsilon ^4}\right),$$
		where $\hat{\delta}_x:= \max\{\|x\| \mid x \in  \mathcal{X}\}$,
		$$\hat{r}_2:=\hat{\mathcal{L}}_0-\hat{\mathcal{L}}_2+ \frac{8\cdot 2^{1/4} }{\bar{\varrho}}\hat{\delta}_x^2 + \frac{19}{40\bar{\varrho}} \hat{\delta}_x^2,$$
		$\hat{r}_3:=\max\{\hat{d}_1,\frac{10+20K_2\bar{\varrho}^2 L_{12}^2}{9\bar{\varrho}\alpha_2}\}$ with $\hat{d}_1=\frac{16K_2\tau^2}{(\tau-2)^2} + \frac{19^4\left(K_2\left(\bar{\varrho} L_{21}^2-\bar{\iota}+L_{22}\right)^2 \right)}{16\cdot20^4 \bar{\varrho}^2 (\tau-2)^2 L_{21}^4}$, $\hat{\mathcal{L}}_0 :=\hat{l}+\left(2^{13/4}+15+\frac{19}{40}\right)\frac{\hat{\delta}_x^2}{\bar{\varrho}}$ with $\hat{l} := \max_{(x,y) \in \mathcal{X} \times \mathcal{Y}} l(x,y)$, $\hat{\mathcal{L}}_2:=L(x_2,y_2)-\frac{16\hat{\delta}_x^2}{\bar{\varrho}^2\hat{b}_{2}}-\frac{8\hat{\delta}_x^2}{\bar{\varrho}}-\left( \frac{\bar{\varrho} L_{21}^2}{2}+\frac{16L_{21}^2}{\bar{\varrho} (\hat{b}_{2})^2} \right)\| y_{2}-y_1 \| ^2$.
	\end{thm}

Theorem \ref{nthm:5.4} implies that the iteration complexity of Algorithm \ref{alg:BAPG} to obtain an $\varepsilon$-stationary point for general block-wise nonsmooth convex-nonconcave minimax problems \eqref{problem:block} is bounded by $\mathcal{O}\left( L^4\varepsilon^{-4} \right)$.


\section{Numerical results}
In this section, we compare the numerical performance of the proposed AGP algorithm with the gradient descent ascent algorithm (GDA),  the alternating gradient descent ascent algorithm (AGDA) and the state-of-the-art nested-looped algorithm, i.e., MINIMAX-PPA in \cite{Lin2020} through two representative test problems.  The first numerical test is implemented in Python 3.9 and run with an Apple M1 processor, while the second one is carried out on an NVIDIA Tesla P100 GPU.

\subsection{\textbf{Dirac-GAN Problem}}
The Dirac-GAN problem \cite{Mescheder2018} can be formulated as the following nonconvex-concave minimax problem :
\begin{align}
	\min\limits_{x}\max\limits_{y} L(x,y)=-\log(1+exp(-x y))+\log 2,
\end{align}
where $(0,0)$ is the unique stationary point.

\textbf{Experimental setup.} Let $\alpha_{x}$ and $\beta_{y}$ be the stepsize of x and y, respectively, and $T$ denote the outer-loop iteration number for AGP, GDA, AGDA and MINIMAX-PPA algorithms. The initial point of all algorithms is chosen as $(1,1).$

For the AGP algorithm, we set $T=72,\alpha_{x}= \frac{0.8}{\sqrt{k}}$, $\beta_{y}=0.3$ and $c_{k}=\frac{0.5}{k^{1/4}}$. For the GDA and AGDA algorithms, we set $T=100$, $\alpha_{x}= 0.3$ and $\beta_{y}=0.3$. For MINIMAX-PPA algorithm, we set $T=30,\epsilon=0.01,l=2,\mu_{y}=0.000625$ and  $g(x,y)=L(x,y)-0.00125\|y-1\|^{2}$. Moreover, for MINIMAX-PPA algorithm, we choose $\delta=1e-6$ and $\tilde{\epsilon}=5.8e-9$ through grid search in the sub-rountines AG2 and AGD respectively, since the theoretical settings for these two parameters are $\delta=2.38e-27$ and $\tilde{\epsilon}=4.06e-60$ respectively, which will cause the algorithm to converge very slowly in practice.

\begin{figure}
	\centering
	\includegraphics[width=160pt]{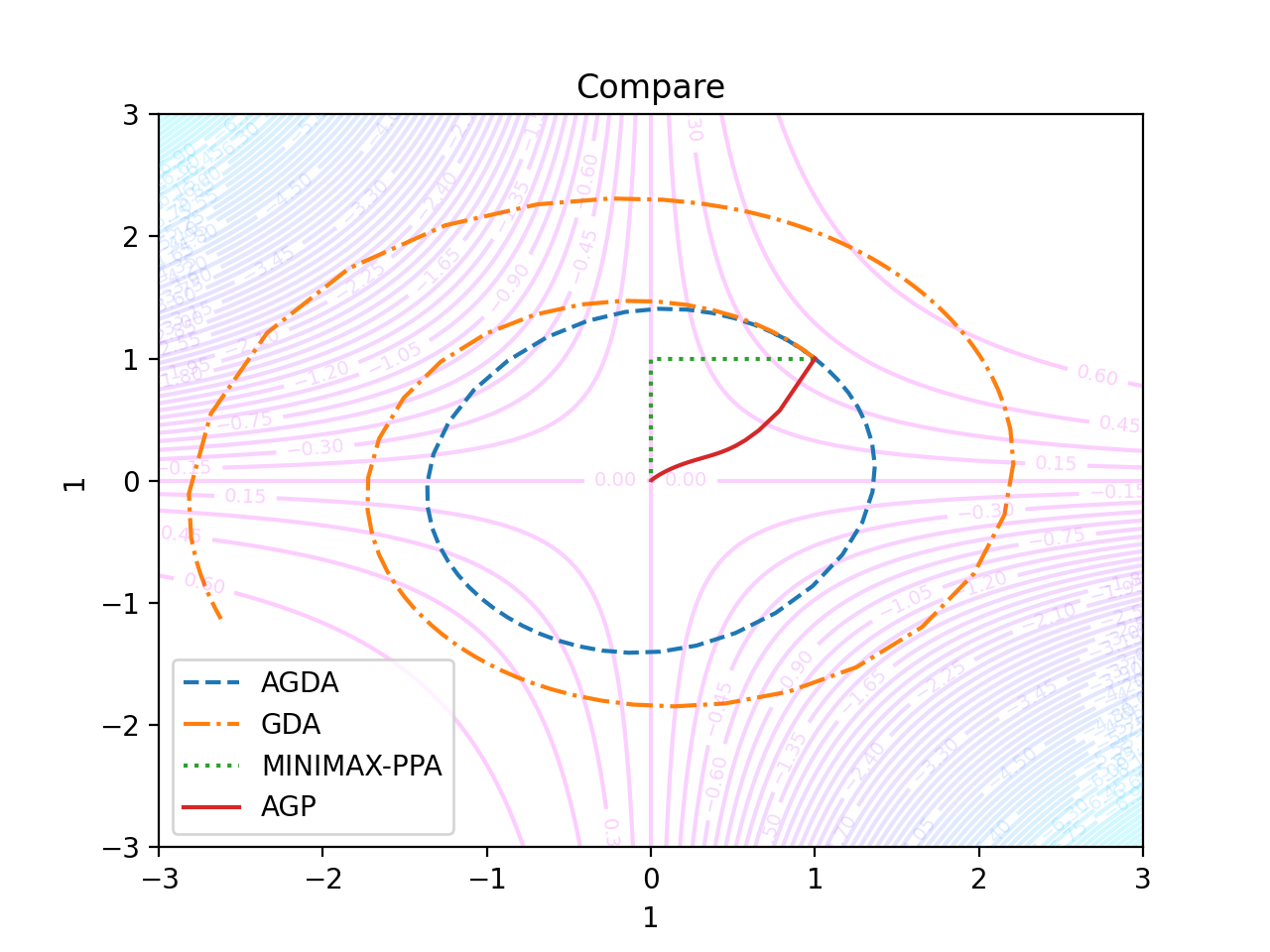}
	\caption{Comparison of four algorithms for solving the dirac-GAN problem}
	\label{fig_a}
\end{figure}

\textbf{Results.} Fig \ref{fig_a} shows the sequences generated by these four algorithms. We find that GDA fails to converge, and AGDA falls into a limit cycle no matter what stepsize is chosen, which  is shown in Lemma 2.4 in \cite{Mescheder2018}. Numerical results show that both AGP and MINIMAX-PPA algorithms approach the unique stationary point. Note that there are much more parameters needed to be adjusted in advance in MINIMAX-PPA algorithm than the AGP algorithm. Table \ref{tab1} shows the total number of iterations and CPU times for both algorithms in solving the dirac-GAN problem. We count the total number of iterations in all the sub-rountines in MINIMAX-PPA in \cite{Lin2020}. The ``Distance" column reports the Euclidean distance between the point obtained by both algorithms and the unique stationary point. Note that the theoretical settings for $\tilde{\epsilon}$ in the sub-rountine AGD  in MINIMAX-PPA is actually $4.06e-60$ which is almost unacceptable in training a neural network with millions of parameters. In this test, we reduce the precision for the sub-rountine AGD  in MINIMAX-PPA algorithm in order to obtain a meaningful solution, otherwise it will get stuck. To obtain a solution near $(0,0)$ with the accuracy within $0.01$, MINIMAX-PPA algorithm needs millions number of iterations totally, and hence takes much more CPU time than the AGP algorithm.
\begin{table}
	\centering
	\caption{Comparison of iteration number and time for AGP algorithm and MINIMAX-PPA algorithm}
	\begin{tabular}{lllll}
		\hline
		& Distance  	& Total iteration    & Time(s)\\ \hline
		MINIMAX-PPA & 0.01 & 18812            & 0.0313   \\
		AGP & 0.01 & 144                   & 0.0001 \\ \hline
		
	\end{tabular}
	\label{tab1}
\end{table}

%
%
%

\subsection{\textbf{Robust learning over multiple domains}}
In this subsection,  we perform some numerical tests for solving a robust learning problem over multiple domains \cite{Qian}, formulated as a nonconvex-linear minimax problem:
\begin{equation}
	\min _{x} \max _{y \in \Delta} y^{T} F(x),
\end{equation}
where $F(x):=\left[f_{1}(x) ; \ldots ; f_{M}(x)\right] \in \mathbb{R}^{M \times 1}$ with $f_{m}(x)=\frac{1}{\left|\mathcal{S}_{m}\right|} \sum_{i=1}^{\left|\mathcal{S}_{m}\right|} \ell\left(d_{i}^{m}, l_{i}^{m}, x\right)$ and $\ell$ can be any non-negative loss function,  $ M $ is the number of tasks and $ x $ denotes the network parameters. $ y \in \Delta$ describes the weights over different tasks and $ \Delta $ is the simplex, i.e., $\Delta=\left\{(y_1,\cdots,y_M)^{T} \mid \sum_{m=1}^{M} y_{m}=1,  y_{m} \geq 0\right\}$. We compare the proposed AGP algorithm with the  GDA algorithm and a heuristic training algorithm with even weights, i.e., fixing $y_1=y_2=0.5$ when $M=2$, which is also used as a baseline algorithm in many other existing works \cite{Lu}.
We consider two image classification problems with MNIST \cite{yangMNIST} and CIFAR10 datasets\cite{Alex}. MNIST consists of images about handwritten digits while CIFAR10 contains 10 different classes such as airplanes, cars, birds etc. Our goal is to train a neural network that works on these two completely unrelated problems simultaneously. Since cross-entropy loss is popular in multi-class classification problem, we use it as our loss function. We should point out that the quality of the algorithm for solving a robust learning problem over multiple domains is measured by the worst case accuracy over all domains~\cite{Qian}.

\textbf{Data augmentation}. Since the data in MNIST are $ 28\times28 $ gray images while those in CIFAR10 are $32\times32$ RGB images, we first repeat gray channel to RGB channels and resize the examples in MNIST from $ 28\times28 $ to $ 32\times32 $ using bilinear interpolation. We adopt AlexNet \cite{Alexnet} as our base model which is the same as the one used in \cite{Lu}. Then, we convert $ 32\times32 $ to $ 224\times224 $ with the same method to fit the input of AlexNet.

\textbf{Experiment setup.} We set $ \frac{1}{\beta_{k}} $, $ \gamma_{k} $, and $ c_{k} $ to be $ \frac{2}{2+\sqrt{k}} $, $ 100 $, $ \frac{1}{10+k^{\frac{1}{4}}} $ respectively for the AGP algorithm. In the GDA algorithm, we set $ \alpha=0.02$ and $\beta=0.01$. In the Heuristic Algorithm, we set the stepsize of $ x $ to be 0.02, but fix $ y $ as $ (0.5,0.5)^{T} $. We set the batch size in all tasks to be $ 128 $ and run $ 50 $ epochs for all algorithms. Moreover, we sample our results every epoch and calculate the average accuracy on these two testing sets to evaluate the performance of different algorithms.

\begin{figure}[h]
	\centering
	\subfigure[]{
		\label{fig_c}
		\includegraphics[width=160pt]{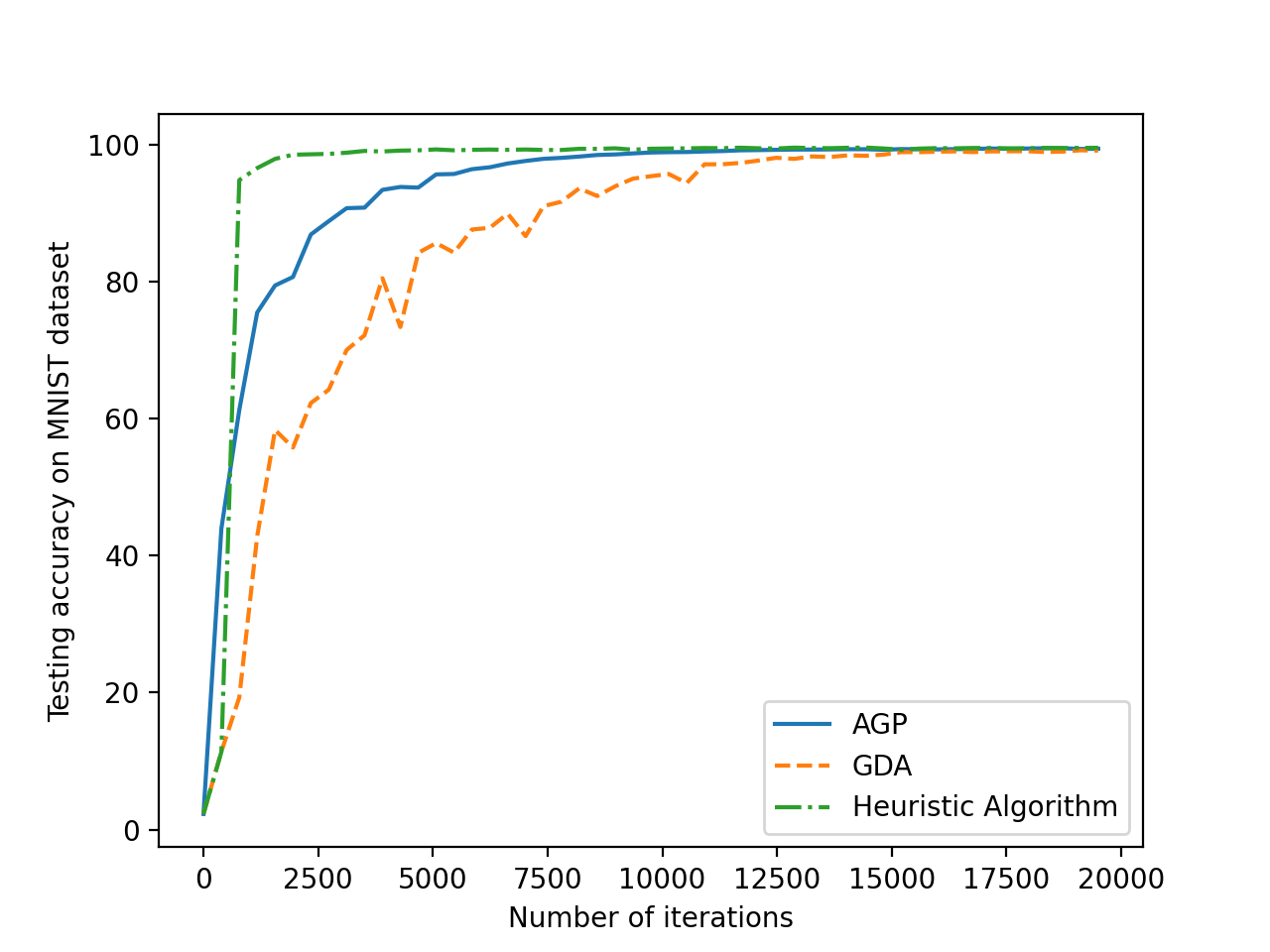}}
	\subfigure[]{
		\label{fig_d}
		\includegraphics[width=160pt]{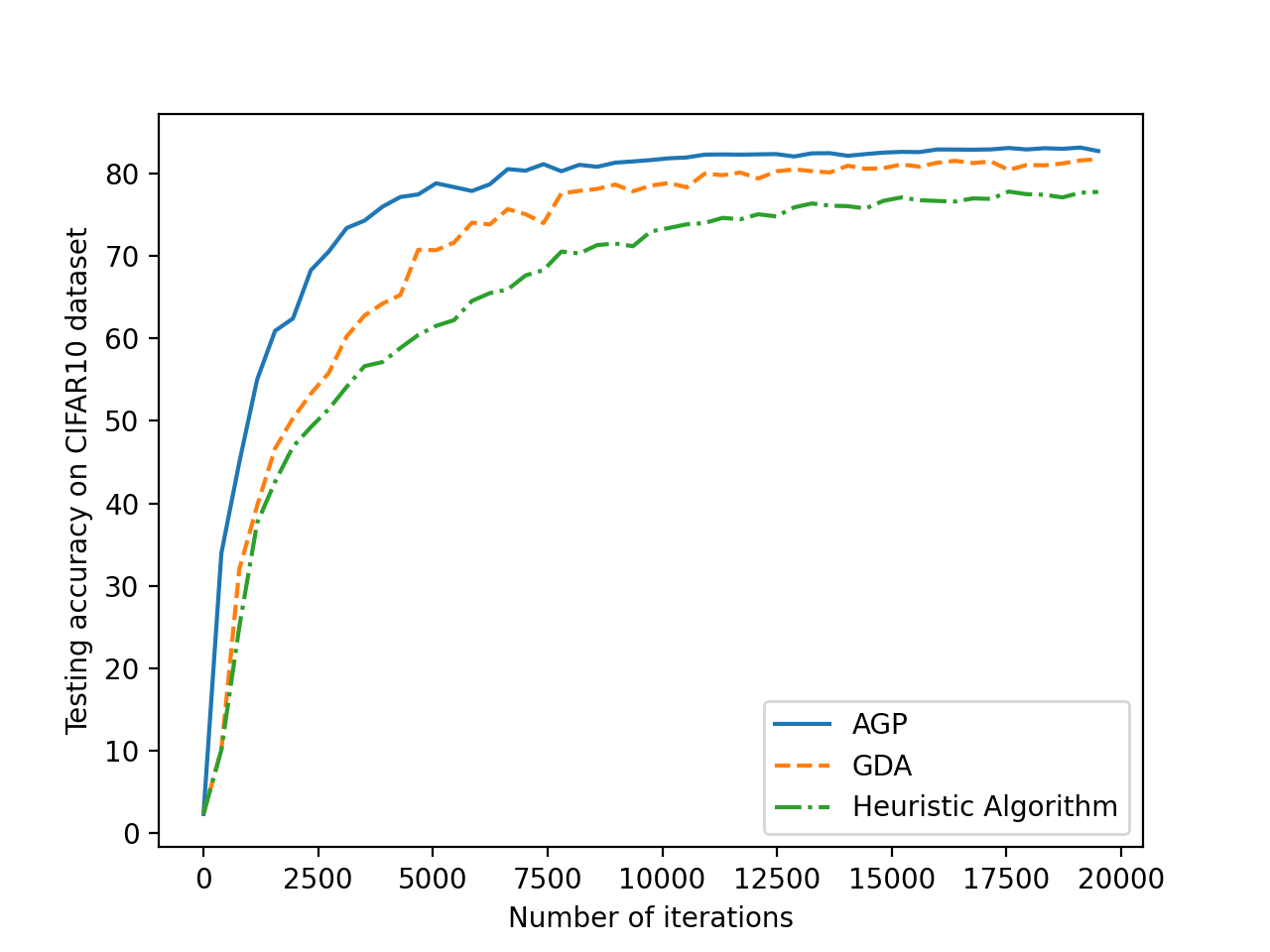}}
	\label{img1}
	\caption{Performance of three algorithms for solving robust multi-task learning problem.}
\end{figure}

\begin{table}[h]
	\centering
	\begin{tabular}{|l|l|l|l|l|}
		\hline
		\multirow{2}{*}{Method} & \multicolumn{2}{l|}{MNIST} & \multicolumn{2}{l|}{CIFAR10} \\ \cline{2-5}
		& Training set & Testing set & Training set  & Testing set  \\ \hline
		Heuristic Algorithm                    & 99.976       & 99.479      & 99.868        & 77.694       \\ \hline
		GDA                     & 99.355       & 99.128      & 99.994        & 81.611       \\ \hline
		AGP                     & 99.954       & 99.329      & 99.998        & 83.173       \\ \hline
	\end{tabular}
	\caption{Comparison the test accuracies for three algorithms on robust multi-task learning}
	\label{tab2}
\end{table}

\textbf{Results.}  Figure \ref{fig_c} shows the testing accuracy on MNIST dataset. All three algorithms achieve high precision, while AGP algorithm takes less time than that of GDA. Figure \ref{fig_d} shows the testing accuracy on CIFAR10 dataset. AGP algorithm still performs slightly better than other two algorithms.

Table \ref{tab2} shows the accuracies of all algorithms. MNIST and CIFAR10 indicate training only on MNIST and CIFAR10 dataset respectively. GDA and Heuristic Algorithm can achieve good performance while AGP algorithm can further improve the performance and provide a more reliable model for multi-task learning.

\textbf{Remark.} Due to hardware (especially memory) limitations, we use mini-batch randomized gradient instead of the exact gradient at each iteration in our numerical experiment, which is the same as in \cite{Lu}. Since at each iteration solving the inner subproblem will take huge amounts of time in the nested loop algorithms, e.g., MINIMAX-PPA algorithm, as the gradient calculation needs to go through the whole dataset and they need to calculate gradient several times per iteration. Moreover,   high requirements for precision in solving subproblems will occupy a lot of memory. Hence, we were not able to compare AGP algorithm with nested-loop algorithm for this test problem.

	\section{Conclusions and Discussion}
	In this paper, we propose a unified single-loop algorithm for general smooth one block nonconvex-concave or convex-nonconcave minimax problems. At each iteration, only one gradient projection step is employed for updating $x$ and $y$ respectively. The gradient complexity of the proposed unified AGP algorithm under four different settings are established. We prove that an $\varepsilon$-first order stationary point of $f$ can be obtained in $\mathcal{O}\left( \varepsilon ^{-4} \right)$ (resp. $\mathcal{O}\left( \varepsilon ^{-2} \right)$) iterations for nonconvex-concave (resp. nonconvex-strongly concave) setting. To the best of our knowledge, they are the best complexity bounds among single-loop algorithms in general nonconvex-(strongly) concave settings, and very simple to be implemented. Numerical results show the efficiency of the proposed AGP algorithm. Nonetheless, there is still a certain gap between the iteration complexity of the proposed AGP algorithm and the best complexity of $\mathcal{O}(\varepsilon^{-2.5})$ among nested-looped algorithms in nonconvex-concave setting, which is also an open problem worth studying in the future.

	Moreover, we consider the (strongly) convex-nonconcave setting of \eqref{problem:1}. 
	Under this setting, the whole problem is convex. However, for any given $x$, to solve the inner max subproblem, i.e., $\max_{y\in \mathcal{Y}} f(x, y)$, is already NP-hard. Due to this reason, almost all the existing nested-loop algorithms will lose the theoretic guarantee since they all need to solve the inner subproblem exactly, or inexactly but only an error proportional to the accuracy $\varepsilon$ is allowed. 
	We show that the proposed unified single-loop AGP algorithm can deal with general convex-nonconcave setting.
	More specifically, the gradient complexity to obtain an $\varepsilon$-first-order stationary point of $f$ is $\mathcal{O}\left( \varepsilon ^{-2} \right)$ (resp., $\mathcal{O}\left( \varepsilon ^{-4} \right)$) under the strongly convex-nonconcave (resp.,  convex-nonconcave ) setting. To the best of our knowledge, 	these theoretical performance guarantees under this setting have not been obtained before in the literature.
	
 Furthermore, we consider more general nonsmooth multi-block nonconvex-(strongly) concave and (strongly) convex-nonconcave minimax problems, which include smooth one block setting as a special case. We propose a block alternating proximal gradient algorithm (BAPG) to solve it.  At each iteration, only simple proximal gradient steps are employed for updating one block variable, and for each block of multi-block variables alternatively. We prove the gradient complexity of the proposed BAPG algorithm under four different settings. To the best of our knowledge, these are the state-of-the-art single loop algorithms under nonconvex-concave setting and the first two theoretically guaranteed convergence results reported in the literature under (non)smooth multi-blocks (strongly) convex-nonconcave minimax problems.  
Our development shows that single loop algorithms do not need to solve the inner subproblem, and thus can be easily generalized.
This is also the reason that the iteration complexity analysis will be more difficult than that of nested-looped algorithms. It will also be interesting to study
 whether the iteration complexity of the proposed AGP algorithm is tightest or not among single loop algorithms.

\begin{appendix}
		\section{Proof of Theorem 5.1}	
		We prove the following two lemmas before giving the proof of Theorem 5.1.
		
		\begin{lem}\label{lem:5.2}
			Suppose that Assumption 2.1 holds. Let $\{\left(x_k,y_k\right)\}$ be a sequence generated by Algorithm 2
			with parameter settings in (5.6). If $\hat{\eta} > L_{11}$, then we have
			\begin{align}\label{lem:5.2:1}
				&l( x_{k+1},y_{k+1}) -l\left( x_k,y_k \right)  \nonumber \\
				\leq& -\left(\tfrac{\hat{\eta}}{2}-\tfrac{L_{12}^{2}\hat{\rho}}{2}\right)\| x_{k+1}-x_{k} \| ^2-\left(\tfrac{\hat{\mu}}{2}-\tfrac{1}{\hat{\rho}}\right)\| y_{k+1}-y_k \|^2-\left(\hat{\mu}-\tfrac{1}{2\hat{\rho}}-\tfrac{\hat{\rho }L_{22}^2}{2}\right)\|y_k-y_{k-1} \|^2.
			\end{align}
		\end{lem}
		
		\begin{proof}
			The optimality condition for $y_k$ in (5.5) implies that $\forall y\in \mathcal{Y}$ and $\forall k\geq 1$,
			\begin{equation}\label{lem:5.2:2}
				\langle \nabla _yf(x_{k+1},y_k)-\frac{1}{\hat{\rho}}(y_{k+1}-y_k)-\partial g_1(y_{k+1}),y-y_{k+1} \rangle \le 0.
			\end{equation}
			By choosing $y=y_k$ in \eqref{lem:5.2:2}, we have
			\begin{equation}\label{lem:5.2:3}
				\langle \nabla _yf(x_{k+1},y_k)-\frac{1}{\hat{\rho}}(y_{k+1}-y_k)-\partial g_1(y_{k+1}),y_k-y_{k+1} \rangle \le 0.
			\end{equation}
			On the other hand, by replacing $k$ with $k-1$ and choosing $y=y_{k+1}$ in \eqref{lem:5.2:2}, we obtain
			\begin{equation}\label{lem:5.2:4}
				\langle \nabla _yf(x_{k},y_{k-1})-\frac{1}{\hat{\rho}}(y_k-y_{k-1})-\partial g_1(y_{k}),y_{k+1}-y_k \rangle \le 0,
			\end{equation}
			which, in view of the fact that $f\left( x,y  \right)$ is $\hat{\mu} $-strongly concave w.r.t. $y$ for any given $x\in \mathcal{X}$ and $g_1(y)$ is convex,
			then implies that
			\begin{align}\label{lem:5.2:5}
				&l( x_{k+1},y_{k+1}) -l( x_{k+1},y_k ) \nonumber\\
				\leq &  \langle  \nabla _yf( x_{k+1},y_k )-\partial g_1(y_k),y_{k+1}-y_k\rangle  -\frac{\hat{\mu}}{2}\| y_{k+1}-y_k \|^2\  \nonumber\\
				\leq &\langle  \nabla _yf( x_{k+1},y_k ) -\nabla _yf( x_k,y_{k-1}) ,y_{k+1}-y_k \rangle +\frac{1}{\hat{\rho}}\langle  y_k-y_{k-1},y_{k+1}-y_k\rangle  -\frac{\hat{\mu}}{2}\| y_{k+1}-y_k \|^2.
			\end{align}
			Denoting $v_{k+1}:=\left(y_{k+1}-y_k\right)-\left(y_k-y_{k-1}\right)$, we can write the first inner product term in the r.h.s.
			of \eqref{lem:5.2:5} as
			\begin{align}\label{lem:5.2:6}
				&\langle  \nabla _y f ( x_{k+1},y_k ) -\nabla _y f ( x_k,y_{k-1}) ,y_{k+1}-y_k \rangle  \nonumber\\
				=&\langle  \nabla _y f ( x_{k+1},y_k ) -\nabla _y f \left( x_k,y_{k} \right) ,y_{k+1}-y_k \rangle  + \langle  \nabla _y f \left( x_k,y_{k} \right) -\nabla _y f ( x_k,y_{k-1}) ,v_{k+1} \rangle  \nonumber\\
				&+\langle  \nabla _y f \left( x_k,y_{k}\right) -\nabla _y f ( x_k,y_{k-1}) ,y_k-y_{k-1} \rangle .
			\end{align}
			Next, we estimate the three terms in the right hand side of \eqref{lem:5.2:6} respectively.
			By Assumption 2.1 and the Cauchy-Schwarz inequality, we can bound the first two terms according to
			\begin{align}\label{lem:5.2:7}
				&\langle  \nabla _yf( x_{k+1},y_k ) -\nabla _yf\left( x_k,y_{k} \right) ,y_{k+1}-y_k \rangle
				\leq \frac{L_{12}^{2}\hat{\rho}}{2}\| x_{k+1}-x_k \|^2+\frac{1}{2\hat{\rho}}\| y_{k+1}-y_k \|^2,
			\end{align}
			and
			\begin{align}\label{lem:5.2:8}
				&\langle  \nabla _yf\left(  x_k,y_{k} \right) -\nabla _yf( x_k,y_{k-1}) ,v_{k+1} \rangle  \leq \frac{\hat{\rho} L_{22}^{2}}{2}\|y_k-y_{k-1} \|^2+\frac{1}{2\hat{\rho}}\|v_{k+1} \|^2.
			\end{align}
			For the third term, by $\hat{\mu} $-strongly-concavity  of $f$ with respect to $y$,
			\begin{align}\label{lem:5.2:9}
				\quad \langle  \nabla _yf\left( x_k,y_{k} \right) -\nabla _yf( x_k,y_{k-1}) ,y_k-y_{k-1} \rangle  \le -\hat{\mu} \| y_k-y_{k-1} \|^2 .
			\end{align}
			Moreover, it can be easily checked that
			\begin{align}\label{lem:5.2:10:3points}
				&\langle  y_k-y_{k-1},y_{k+1}-y_k \rangle =\frac{1}{2}\| y_k-y_{k-1} \|^2+\frac{1}{2}\| y_{k+1}-y_k \|^2-\frac{1}{2}\| v_{k+1} \|^2.
			\end{align}
			Plugging \eqref{lem:5.2:6}-\eqref{lem:5.2:10:3points} into \eqref{lem:5.2:5} and rearranging the terms, we conclude that
			\begin{align}\label{lem:5.2:11}
				&l( x_{k+1},y_{k+1}) -l( x_{k+1},y_{k})  \nonumber \\
				\leq& \frac{L_{12}^{2}\hat{\rho}}{2}\| x_{k+1}-x_k \|^2-(\hat{\mu}-\frac{1}{2\hat{\rho}}-\frac{\hat{\rho} L_{22}^2}{2})\| y_k-y_{k-1}\|^2-(\frac{\hat{\mu}}{2}-\frac{1}{\hat{\rho}})\|y_{k+1}-y_k \|^2.
			\end{align}
			The proof is completed by combining \eqref{lem:5.2:11} with (5.7) in Lemma 5.1.
		\end{proof}

		\begin{lem}\label{lem:5.3}
			Suppose that Assumption 2.1 holds. Let $\{\left(x_k,y_k\right)\}$ be a sequence generated by Algorithm 2
			with parameter settings in (5.6).
			Also let us denote
			\begin{align*}
				l_{k+1} := l( x_{k+1},y_{k+1}),~
				S_{k+1} :=\frac{2}{\hat{\rho}^2\hat{\mu}}\|y_{k+1}-y_k \|^2, \\
				F_{k+1} :=l_{k+1}+S_{k+1}- (\hat{\mu}+\frac{7}{2\hat{\rho}}-\frac{\hat{\rho} L_{22}^2}{2}-\frac{2L_{22}^2}{\hat{\mu}})\|y_{k+1}-y_k\|^2.
			\end{align*}
			If $\hat{\eta}>L_{11}$,
			then
			\begin{align}\label{lem:5.3:2}
				F_{k+1}-F_k\leq& -\left(\frac{\hat{\eta}}{2}-\frac{\hat{\rho} L_{12}^2}{2}-\frac{2L_{12}^2}{\hat{\rho} \hat{\mu}^2}\right)\| x_{k+1}-x_{k}\|^2 -\left(\frac{3\hat{\mu}-\hat{\rho}L_{22}^2}{2}+\frac{\hat{\mu}-4\hat{\rho} L_{22}^2}{2\hat{\rho}\hat{\mu}}\right)\| y_{k+1}-y_k \|^2.
			\end{align}
		\end{lem}
		
		\begin{proof}
			First by \eqref{lem:5.2:3} and \eqref{lem:5.2:4}, we have
			\begin{align}\label{lem:3.3:3}
				\frac{1}{\hat{\rho}}\langle v_{k+1},y_{k+1} -y_k \rangle
				\le \langle \nabla _yf( x_{k+1},y_k ) -\nabla _yf\left( x_{k},y_{k-1} \right)-\partial g_1(y_{k+1})+\partial g_1(y_k),y_{k+1}-y_k \rangle,
			\end{align}
			By the $\langle \partial g_1(y_{k+1})-\partial g_1(y_k),y_{k+1}-y_k \rangle \geq0$, and  together with \eqref{lem:5.2:6} then imply that
			\begin{align}\label{lem:5.3:4}
				\frac{1}{\hat{\rho}}\langle v_{k+1},y_{k+1} -y_k  \rangle  \leq&\langle  \nabla _yf( x_{k+1},y_k ) -\nabla _yf\left( x_{k},y_k \right) ,y_{k+1}-y_k \rangle  \nonumber\\
				& + \langle  \nabla _yf\left( x_{k},y_k \right) -\nabla _yf\left( x_{k},y_{k-1} \right) ,v_{k+1} \rangle  \nonumber\\
				&+\langle  \nabla _yf\left( x_{k},y_k \right) -\nabla _yf\left( x_{k},y_{k-1} \right) ,y_{k}-y_{k-1} \rangle .
			\end{align}
			Similar to \eqref{lem:5.2:7}, we can easily see that
			\begin{align}\label{lem:5.3:5}
				\left\langle  \nabla _yf( x_{k+1},y_k ) -\nabla _yf\left( x_{k},y_k  \right) ,y_{k+1}-y_k \right\rangle
				\leq &\frac{L_{12}^{2}}{2 \hat{\mu}}\|x_{k+1}-x_{k} \|^2+\frac{\hat{\mu}}{2}\|y_{k+1}-y_k \|^2.
			\end{align}
			By plugging \eqref{lem:5.2:8},\eqref{lem:5.2:9},\eqref{lem:5.3:5} into \eqref{lem:5.3:4},
			and using the identity
			$\frac{1}{\hat{\rho}}\langle v_{k+1}, y_{k+1}-y_k \rangle  = \frac{1}{2\hat{\rho}}\|y_{k+1}-y_k\|^2+\frac{1}{2\hat{\rho}}\| v_{k+1}\| ^2-\frac{1}{2\hat{\rho}}\|y_{k}-y_{k-1}\|^2$,
			we conclude that
			\begin{align}\label{lem:5.3:6}
				&\frac{1}{2\hat{\rho}}\|y_{k+1}-y_k\|^2+\frac{1}{2\hat{\rho}}\| v_{k+1}\|^2-\frac{1}{2\hat{\rho}}\|y_{k}-y_{k-1}\|^2\nonumber\\
				\leq& \frac{L_{12}^{2}}{2\hat{\mu}}\| x_{k+1}-x_{k}\|^2+\frac{\hat{\mu}}{2}\| y_{k+1}-y_k \|^2+\frac{\hat{\rho} L_{22}^2}{2}\|y_{k}-y_{k-1} \|^2+\frac{1}{2\hat{\rho}}\| v_{k+1}\|^2 -\hat{\mu}\|y_{k}-y_{k-1} \|^2.
			\end{align}
			Rearranging the terms of \eqref{lem:5.3:6}, we have
			\begin{align}\label{lem:5.3:7}
				&\frac{1}{2\hat{\rho}}\|y_{k+1}-y_k\|^2-\frac{1}{2\hat{\rho}}\|y_{k}-y_{k-1}\|^2\nonumber\\
				\leq& \frac{L_{12}^{2}}{2\hat{\mu}}\| x_{k+1}-x_{k}\|^2+\frac{\hat{\mu}}{2}\| y_{k+1}-y_k \|^2-\left(\hat{\mu}-\frac{\hat{\rho }L_{22}^2}{2}\right)\|y_{k}-y_{k-1} \|^2.
			\end{align}
			Multiplying $\frac{4}{\hat{\rho}\hat{\mu}}$ on both sides of \eqref{lem:5.3:7} and using the definition of $S_{k+1}$, we obtain
			\begin{align*}
				S_{k+1}-S_k \leq & \frac{2L_{12}^{2}}{\hat{\mu}^2\hat{\rho}}\| x_{k+1}-x_{k}\|^2+\frac{2}{\hat{\rho}}\| y_{k+1}-y_k \|^2-\left(\frac{4}{\hat{\rho}}-\frac{2 L_{22}^2}{\hat{\mu}}\right)\|y_{k}-y_{k-1} \|^2.
			\end{align*}
			It then follows from \eqref{lem:5.2:1} in Lemma \ref{lem:5.2} and the definition of $F_k$ that
			\begin{align*}
				F_{k+1}-F_k\leq& -\left(\frac{\hat{\eta}}{2}-\frac{\hat{\rho} L_{12}^2}{2}-\frac{2L_{12}^2}{\hat{\rho} \hat{\mu}^2}\right)\| x_{k+1}-x_{k}\|^2 -\left(\frac{3\hat{\mu}-\hat{\rho}L_{22}^2}{2}+\frac{\hat{\mu}-4\hat{\rho} L_{22}^2}{2\hat{\rho}\hat{\mu}}\right)\| y_{k+1}-y_k \|^2.
			\end{align*}
		\end{proof}

		{\bf The proof of Theorem 5.1}
		\begin{proof}
			Noting that (5.4) is equivalent to $ x_{k+1}^{(i)}= \operatorname{ Prox}_{h_i,\mathcal{X}_i}^{\hat{\eta}}\left( x^{(i)}_k-\frac{1}{\hat{\eta}} \nabla _{x^{(i)}}f\left(v^{(i)}_{k+1}, y_{k} \right) \right)$, we immediately obtain
			\begin{align}\label{thm:blockstr:gx}
				&\|(\nabla \mathcal{G}_k)_{x^{(i)}} \|\nonumber\\
				\leq&\hat{\eta}\|\operatorname{ Prox}_{h_i,\mathcal{X}_i}^{\hat{\eta}}\left( x^{(i)}_k-\frac{1}{\hat{\eta}} \nabla _{x^{(i)}}f\left( v^{(i)}_{k+1},y_{k} \right) \right)-\operatorname{ Prox}_{h_i,\mathcal{X}_i}^{\hat{\eta}}\left( x^{(i)}_k- \frac{1}{\hat{\eta}} \nabla _{x^{(i)}}f\left( x_k,y_k \right) \right) \|\nonumber\\
				&+\hat{\eta}\| x^{(i)}_{k+1}-x^{(i)}_k \|\nonumber\\
				\leq& \hat{\eta}\| x^{(i)}_{k+1}-x^{(i)}_k \| +\|\nabla _{x^{(i)}}f\left( v^{(i)}_{k+1},y_{k} \right)-\nabla _{x^{(i)}}f\left( x_k,y_k \right) \|\nonumber\\
				\le&\hat{\eta}\| x^{(i)}_{k+1}-x^{(i)}_k \|+L_{11}\|v^{(i)}_{k+1}-x_k \|\nonumber\\
				\le& \left(\hat{\eta}+L_{11}\right)\| x_{k+1}-x_k \|,
			\end{align}
			where the second inequality holds by the nonexpansiveness of the proximal operator $\operatorname{ Prox}_{h_i,\mathcal{X}_i}^{\hat{\eta}}$, the last inequality hold since $\|x^{(i)}_{k+1}-x^{(i)}_k \|\le\| x_{k+1}-x_k \|$, $\|v^{(i)}_{k+1}-x_k \|\le\| x_{k+1}-x_k \|$ by the definition of $x^{(i)}$ and $v^{(i)}_{k+1}$.
			On the other hand, since (5.5) is equivalent to $ y_{k+1}= \operatorname{ Prox}_{g_1, \mathcal{Y}}^{1/\hat{\rho}}\left( y_k+\hat{\rho} \nabla _{y}f\left( x_{k+1} ,y_{k} \right) \right)$, we conclude from the triangle inequality that
			\begin{align}\label{thm:blockstr:gy}
				&\| (\nabla \mathcal{G}_k)_y \|\nonumber\\
				\le&\frac{1}{\hat{\rho}}\|\operatorname{ Prox}_{g_1, \mathcal{Y}}^{1/\hat{\rho}}\left( y_k+\hat{\rho} \nabla _{y}f\left( x_{k+1} ,y_{k} \right) \right)-\operatorname{ Prox}_{g_1, \mathcal{Y}}^{1/\hat{\rho}}\left( y_k+\hat{\rho} \nabla _{y}f\left( x_k ,y_k \right) \right)\|\nonumber\\
				&+\frac{1}{\hat{\rho}} \| y_{k+1}-y_k \|\nonumber\\
				\le&\frac{1}{\hat{\rho}} \| y_{k+1}-y_k \|+L_{12}\|x_{k+1}-x_k\|.
			\end{align}
			By combining \eqref{thm:blockstr:gx}-\eqref{thm:blockstr:gy}, and using Cauchy-Schwarz inequality, we obtain
			\begin{align}\label{thm:blockstr:g}
				\| \nabla \mathcal{G}_k \|^2
				\leq& \left(K_1\left(\hat{\eta}+L_{11}\right)^2+2L_{12}^2\right)\| x_{k+1}-x_k \|^2+\frac{2}{\hat{\rho} ^2}\|y_{k+1}-y_k \|^2.
			\end{align}
			Observing that $r_1 > 0$. Multiplying both sides of \eqref{thm:blockstr:g} by $r_1$, and using \eqref{lem:5.3:2} in Lemma \ref{lem:5.3},
			we have
			\begin{align}\label{thm:blockstr:d1}
				r_1\|\nabla \mathcal{G}_k \|^2\le F_k -F_{k+1}.
			\end{align}
			Summing up the above inequalities from $k=1$ to $k=\mathcal{T}(\varepsilon)$ and using the definition of $L_1$, we obtain
			\begin{align}\label{thm:blockstr:10}
				\sum_{k=1}^{\mathcal{T}\left( \varepsilon \right)}{r_1\|\nabla \mathcal{G}_k \| ^2}\le F_1 -F_{\mathcal{T}\left( \varepsilon \right)+1}\le L_1-F_{\mathcal{T}\left( \varepsilon \right)+1}.
			\end{align}
			Note that by the definition of $F_{k+1}$ in Lemma \ref{lem:5.3}, we have
			\begin{align*}
				F_{\mathcal{T}\left( \varepsilon \right)+1}&=
				l_{\mathcal{T}\left( \varepsilon \right)+1}+S_{\mathcal{T}\left( \varepsilon \right)+1}-(\hat{\mu}+\frac{7}{2\hat{\rho}}-\frac{\hat{\rho} L_{22}^2}{2}-\frac{2L_{22}^2}{\hat{\mu}})\|y_{\mathcal{T}\left( \varepsilon \right)+1}-y_{\mathcal{T}\left( \varepsilon \right)} \|^2\\
				&\ge \underline{l}-(\hat{\mu}+\frac{7}{2\hat{\rho}}-\frac{\hat{\rho} L_{22}^2}{2}-\frac{2L_{22}^2}{\hat{\mu}})\delta_y^2 = \underline{L},
			\end{align*}
			where the inequality follows from the definitions of $\underline{l}$ and  $\delta _y$, and the
			facts that $S_k\ge 0$ ($\forall k\ge 1$) and $\hat{\mu}+\frac{7}{2\hat{\rho}}-\frac{\hat{\rho} L_{22}^2}{2}-\frac{2L_{22}^2}{\hat{\mu}} \ge 0 $ due to
			the selection of $\hat{\rho}$.
			We then conclude from \eqref{thm:blockstr:10} that
			\begin{align*}
				\sum_{k=1}^{\mathcal{T}\left( \varepsilon \right)}{r_1\| \nabla \mathcal{G}_k \|^2}\le L_1 -F_{\mathcal{T}\left( \varepsilon \right)+1}\le L_1-\underline{L},
			\end{align*}
			which, in view of the definition of $\mathcal{T}(\varepsilon)$, implies that $\varepsilon ^2\le (L_1-\underline{L})/(\mathcal{T}( \varepsilon )\cdot r_1)$
			or equivalently, $\mathcal{T}\left( \varepsilon \right) \le (L_1-\underline{L})/(r_1\varepsilon ^2)$.
		\end{proof}
		
		\section{Proof of Theorem 5.2}	
		We prove the following two lemmas before giving the proof of Theorem 5.2.
		\begin{lem}\label{lem:5.5}
			Suppose that Assumptions 2.1 and 5.1 hold. Let $\{(x_k,y_k)\}$ be a sequence generated by Algorithm 2 with parameter settings in (5.11). If $\forall k,\vartheta_k > L_{11}$ and $\xi \le \frac{2}{L_{22}^{'}+\hat{c}_1}$, then
			\begin{align}\label{lem:5.5:1}
				&l(x_{k+1},y_{k+1})-l(x_k,y_k)\nonumber \\
				\leq& -\left( \iota +\frac{\vartheta_k}{2}-\frac{\xi L_{12}^{2}}{2} \right) \| x_{k+1}-x_k \| ^2+\frac{1}{\xi} \|y_{k+1}-y_k \| ^2+\frac{1}{2\xi}\| y_k-y_{k-1} \|^2\nonumber\\
				&+\frac{\hat{c}_{k-1}}{2}(\| y_{k+1}\| ^2-\| y_k\|^2).
			\end{align}
		\end{lem}
		\begin{proof}
			The optimality condition for $y_k$ in (5.5) implies that, $\forall y\in \mathcal{Y}$, $\forall k\geq 1$,
			\begin{equation}\label{lem:5.5:2}
				\langle \nabla _y\bar{f}_k(x_{k+1},y_k)-\frac{1}{\xi}(y_{k+1}-y_k) -\partial g_1(y_{k+1}) ,y-y_{k+1} \rangle \le 0.
			\end{equation}
			By choosing $y=y_k$ in \eqref{lem:5.5:2}, we have
			\begin{equation}\label{lem:5.5:3}
				\langle \nabla _y\bar{f}_k(x_{k+1},y_k)-\frac{1}{\xi}(y_{k+1}-y_k)-\partial g_1(y_{k+1}),y_k-y_{k+1} \rangle \le 0.
			\end{equation}
			On the other hand, by replacing $k$ with $k-1$ and choosing $y=y_{k+1}$ in \eqref{lem:5.5:2}, we obtain
			\begin{equation}\label{lem:5.5:4}
				\langle \nabla _y\bar{f}_{k-1}(x_{k},y_{k-1})-\frac{1}{\xi}(y_k-y_{k-1})-\partial g_1(y_{k}) ,y_{k+1}-y_k \rangle \le 0.
			\end{equation}
			The concavity of  $\bar{f}_k(x_{k+1},y)$ w.r.t. $y$ together with  \eqref{lem:5.5:4} and $g_1(y)$ is convex, then imply that
			\begin{align}\label{lem:5.5:5}
				&\bar{l}_k(x_{k+1},y_{k+1})-\bar{l}_k(x_{k+1},y_k)\nonumber\\
				\le& \langle \nabla _y\bar{f}_k(x_{k+1},y_k)-\nabla _x\bar{f}_{k-1}(x_{k},y_{k-1}),y_{k+1}-y_k \rangle+\frac{1}{\xi}\langle y_k-y_{k-1},y_{k+1}-y_k \rangle\nonumber\\
				=&\langle \nabla _y\bar{f}_k(x_{k+1},y_k)-\nabla _y\bar{f}_{k-1}(x_k,y_{k}),y_{k+1}-y_k \rangle+\langle \nabla _y\bar{f}_{k-1}(x_k,y_{k})-\nabla _y\bar{f}_{k-1}(x_{k},y_{k-1}),v_{k+1} \rangle\nonumber\\
				&+\langle \nabla _y\bar{f}_{k-1}(x_k,y_{k})-\nabla _y\bar{f}_{k-1}(x_{k},y_{k-1}),y_k-y_{k-1} \rangle+\frac{1}{\xi}\langle y_k-y_{k-1},y_{k+1}-y_k \rangle,
			\end{align}
			where $v_{k+1}=y_{k+1}-y_k-(y_k-y_{k-1})$. We now provide bounds on the inner product terms of \eqref{lem:5.5:5}. Firstly, by the definition of $\bar{f}_{k}(x_{k+1},y_k)$ and $\bar{f}_{k-1}(x_{k},y_k)$, Assumptions 2.1 and 5.1, and the Cauchy-Schwarz inequality, we have
			\begin{align}\label{lem:5.5:6}
				&\langle  \nabla _y\bar{f}_{k}( x_{k+1},y_k) -\nabla _y\bar{f}_{k-1}( x_{k},y_{k}) ,y_{k+1}-y_k \rangle  \nonumber\\
				=&\langle \nabla_{y} f(x_{k+1},y_k)-\nabla_{y} f(x_k,y_{k}), y_{k+1}-y_k\rangle-\left(\hat{c}_k-\hat{c}_{k-1}\right)\langle y_k, y_{k+1}-y_k\rangle \nonumber\\
				\le&\frac{\xi L_{12}^{2}}{2}\| x_{k+1}-x_{k} \|^2+\frac{1}{2\xi}\| y_{k+1}-y_k \|^2-\frac{\hat{c}_k-\hat{c}_{k-1}}{2}(\| y_{k+1}\| ^2-\| y_k\|^2)\nonumber\\
				&+\frac{\hat{c}_k-\hat{c}_{k-1}}{2}\| y_{k+1}-y_k \|^2\nonumber\\
				\leq& \frac{\xi L_{12}^{2}}{2}\| x_{k+1}-x_{k} \|^2+\frac{1}{2\xi}\| y_{k+1}-y_k \|^2-\frac{\hat{c}_k-\hat{c}_{k-1}}{2}(\| y_{k+1}\| ^2-\| y_k\|^2).	
			\end{align}
			Secondly, by the Cauchy-Schwarz inequality,
			\begin{align}\label{lem:5.5:7}
				\langle  \nabla _y\bar{f}_{k-1}( x_k,y_{k}) -\nabla _y\bar{f}_{k-1}( x_{k},y_{k-1}) ,v_{k+1} \rangle\leq& \frac{\xi}{2}\|\nabla _y\bar{f}_{k-1}( x_k,y_{k}) -\nabla _y\bar{f}_{k-1}( x_{k},y_{k-1}) \|^2 +\frac{1}{2\xi}\|v_{k+1} \|^2.
			\end{align}
			Thirdly, it follows from the concavity of  $\bar{f}_{k-1}(x_{k},y)$ w.r.t. $y$ that
			\begin{align}\label{lem:5.5:8}
				&\langle \nabla _y\bar{f}_{k-1}\left( x_k,y_{k} \right) -\nabla _y\bar{f}_{k-1}\left( x_{k},y_{k-1} \right) ,y_k-y_{k-1} \rangle\nonumber\\
				\le&-\frac{1}{L_{22}^{'}+\hat{c}_{k-1}}\| \nabla _y\bar{f}_{k-1}\left( x_k,y_{k} \right)-\nabla _y\bar{f}_{k-1}\left( x_{k},y_{k-1} \right) \|^2-\frac{\hat{c}_{k-1}L_{22}^{'}}{L_{22}^{'}+\hat{c}_{k-1}}\|y_k-y_{k-1} \|^2\nonumber\\
				\leq& -\frac{1}{L_{22}^{'}+\hat{c}_{k-1}}\| \nabla _y\bar{f}_{k-1}\left( x_k,y_{k} \right)-\nabla _y\bar{f}_{k-1}\left( x_{k},y_{k-1} \right)\|^2.
			\end{align}
			Also observe that
			\begin{equation}\label{lem:5.5:9}
				\frac{1}{\xi}\langle  y_{k+1}-y_k,y_k-y_{k-1} \rangle =\frac{1}{2\xi}\|y_{k+1}-y_k \|^2+\frac{1}{2\xi}\|y_k-y_{k-1} \|^2-\frac{1}{2\xi}\|v_{k+1} \|^2.
			\end{equation}
			Plugging \eqref{lem:5.5:6}-\eqref{lem:5.5:9} into \eqref{lem:5.5:5}, and using the definition of $\bar{l}_k(x_{k+1},y_{k+1})$ and $\bar{l}_k(x_{k+1},y_k)$ and the assumption $\frac{\xi}{2} \le \frac{1}{L_{22}^{'}+\hat{c}_1}$, we obtain
			\begin{align}\label{lem:5.5:10}
				l(x_{k+1},y_{k+1})-l(x_{k+1},y_k)
				&\leq \frac{\xi L_{12}^{2}}{2}\| x_{k+1}-x_{k}\|^2+\frac{1}{\xi} \|y_{k+1}-y_k \| ^2\nonumber\\
				&\quad+\frac{1}{2\xi}\| y_k-y_{k-1} \|^2+\frac{\hat{c}_{k-1}}{2}(\| y_{k+1}\| ^2-\| y_k\|^2).
			\end{align}
			The result in \eqref{lem:5.5:1} the follows by adding \eqref{lem:5.5:10} and (5.12) in Lemma 5.2.
		\end{proof}

		\begin{lem}\label{lem:5.6}
			Suppose that Assumptions 2.1 and 5.1 hold. Let $\{(x_k,y_k)\}$ be a sequence generated by Algorithm 2 with parameter settings in (5.11). Also let us denote
			\begin{align*}
				\mathcal{S}_{k+1}&:=\frac{8}{\xi^2\hat{c}_{k+1}}\| y_{k+1}-y_k \|^2+\frac{8}{\xi}\left(1-\frac{\hat{c}_k}{\hat{c}_{k+1}}\right)\|y_{k+1}\|^2,\\
				\mathcal{F}_{k+1}&:=l(x_{k+1},y_{k+1})+\mathcal{S}_{k+1}- \frac{15}{2\xi}\| y_{k+1}-y_k \|^2-\frac{\hat{c}_k}{2}\|y_{k+1}\|^2.
			\end{align*}
			If
			\begin{equation}\label{lem:5.6:1}
				\vartheta_k >L_{11},~\frac{1}{\hat{c}_{k+1}}-\frac{1}{\hat{c}_k}\leq \frac{\xi}{5},~\xi \leq \frac{2}{L_{22}^{'}+\hat{c}_1},
			\end{equation}
			then $\forall k \geq 2$,
			\begin{align}\label{lem:5.6:2}
				\mathcal{F}_{k+1}-\mathcal{F}_{k}
				&\leq -\left( \iota+\frac{\vartheta_k}{2}-\frac{\xi L_{12}^2}{2}-\frac{16L_{12}^2}{\xi \hat{c}_{k}^2}\right)\| x_{k+1}-x_k \| ^2+\frac{\hat{c}_{k-1}-\hat{c}_{k}}{2}\|y_{k+1}\|^2\\
				&\quad-\frac{9}{10\xi}\| y_{k+1}-y_k \|^2+\frac{8}{\xi}\left(\frac{\hat{c}_{k-1}}{\hat{c}_k}-\frac{\hat{c}_k}{\hat{c}_{k+1}} \right)\|y_{k+1}\|^2.
			\end{align}
		\end{lem}
		
		\begin{proof}
			By \eqref{lem:5.5:3} and \eqref{lem:5.5:4}, we have
			\begin{align}\label{lem:5.6:3}
				&\frac{1}{\xi}\langle v_{k+1},y_{k+1} -y_k \rangle \le \langle \nabla _y\bar{f}_{k}( x_{k+1},y_k ) -\nabla _y\bar{f}_{k-1}\left( x_{k},y_{k-1} \right)-\partial g_1(y_{k+1})+\partial g_1(y_k),y_{k+1}-y_k \rangle.
			\end{align}
			By the $\langle \partial g_1(y_{k+1})-\partial g_1(y_k),y_{k+1}-y_k \rangle \geq0$, similar to \eqref{lem:5.5:5} in Lemma \ref{lem:5.5}, \eqref{lem:5.6:3} can be rewritten as
			\begin{align*}
				&\frac{1}{\xi}\langle v_{k+1}, y_{k+1} -y_k \rangle\nonumber\\
				\leq&\langle \nabla _y\bar{f}_{k}(x_{k+1},y_k)-\nabla _y\bar{f}_{k-1}(x_k,y_{k}),y_{k+1}-y_k \rangle+\langle \nabla _y\bar{f}_{k-1}(x_k,y_{k})-\nabla _y\bar{f}_{k-1}(x_{k},y_{k-1}),v_{k+1} \rangle\nonumber\\
				&+\langle \nabla _y\bar{f}_{k-1}(x_k,y_{k})-\nabla _y\bar{f}_{k-1}(x_{k},y_{k-1}),y_k-y_{k-1} \rangle.
			\end{align*}
			Using an argument similar to the proof of \eqref{lem:5.5:6}-\eqref{lem:5.5:9},
			the concavity of  $\bar{f}_{k-1}(x_{k},y)$ w.r.t. $y$, and the Cauchy-Schwarz inequality, we conclude from the above inequality that
			\begin{align}\label{lem:5.6:4}	&\frac{1}{2\xi}\|y_{k+1}-y_k\|^2+\frac{1}{2\xi}\|v_{k+1}\|^2-\frac{1}{2\xi}\|y_k-y_{k-1}\|^2\nonumber\\
				\leq& \frac{L_{12}^{2}}{2a_k}\|x_{k+1}-x_{k} \|^2+\frac{a_k}{2}\| y_{k+1}-y_k \|^2-\frac{\hat{c}_k-\hat{c}_{k-1}}{2}(\| y_{k+1}\|^2-\| y_k\| ^2)\nonumber\\
				&+\frac{\xi}{2}\|  \nabla _y\bar{f}_{k-1}(x_k,y_{k})-\nabla _y\bar{f}_{k-1}(x_{k},y_{k-1})\|^2-\frac{1}{L_{22}^{'}+\hat{c}_{k-1}}\|\nabla _y\bar{f}_{k-1}(x_k,y_{k})-\nabla _y\bar{f}_{k-1}(x_{k},y_{k-1}) \|^2\nonumber\\
				&+\frac{1}{2\xi}\|v_{k+1} \|^2-\frac{\hat{c}_{k-1}L_{22}^{'}}{L_{22}^{'}+\hat{c}_{k-1}}\| y_k-y_{k-1} \|^2+ \frac{\hat{c}_k-\hat{c}_{k-1}}{2}\|y_{k+1}-y_k\|^2,
			\end{align}
			for any $a_k>0$. Observing by $\hat{c}_1\leq L_{22}^{'}$, and Assumption 5.1, we have
			\begin{equation*}
				-\frac{\hat{c}_{k-1}L_{22}^{'}}{\hat{c}_{k-1}+L_{22}^{'}}\leq-\frac{\hat{c}_{k-1}L_{22}^{'}}{2L_{22}^{'}}
				=-\frac{\hat{c}_{k-1}}{2}\leq-\frac{\hat{c}_k}{2}.
			\end{equation*}
			Combining $\xi \le \frac{2}{L_{22}^{'}+\hat{c}_1}$ and rearranging the terms in \eqref{lem:5.6:4}, we obtain
			\begin{align*}
				& \frac{1}{2\xi}\| y_{k+1}-y_k \|^2+\frac{\hat{c}_k-\hat{c}_{k-1}}{2}\|y_{k+1}\|^2 \nonumber\\			\leq&\frac{1}{2\xi}\|y_k-y_{k-1}\|^2+\frac{\hat{c}_k-\hat{c}_{k-1}}{2}\|y_k\|^2+\frac{L_{12}^{2}}{2a_k}\|x_{k+1}-x_{k} \|^2+\frac{a_k}{2}\| y_{k+1}-y_k \|^2\nonumber\\
				&-\frac{\hat{c}_k}{2}\|y_k-y_{k-1} \|^2.
			\end{align*}
			By multiplying $\frac{16}{\xi \hat{c}_k}$ on both sides of the above inequality, we then obtain
			\begin{align}\label{lem:5.6:5}
				&\frac{8}{\xi^2\hat{c}_k}\| y_{k+1}-y_k \|^2+\frac{8}{\xi}\left(1-\frac{\hat{c}_{k-1}}{\hat{c}_k}\right)\|y_{k+1}\|^2\nonumber\\
				\leq& \frac{8}{\xi^2\hat{c}_k}\|y_k-y_{k-1}\|^2 +\frac{8}{\xi}(1-\frac{\hat{c}_{k-1}}{\hat{c}_k}
				)\|y_k\|^2+\frac{8L_{12}^2}{\xi \hat{c}_ka_k}\|x_{k+1}-x_{k} \|^2\nonumber\\
				& +\frac{8a_k}{\xi \hat{c}_k}\| y_{k+1}-y_k \|^2 -\frac{8}{\xi}\|y_k-y_{k-1}\|^2.
			\end{align}
			Setting $a_k=\frac{\hat{c}_k}{2}$ in the above inequality, and using the definition of $\mathcal{S}_{k+1}$ and \eqref{lem:5.6:1}, we have
			\begin{align}\label{lem:5.6:6}
				\mathcal{S}_{k+1}-\mathcal{S}_{k}
				&\leq \frac{8}{\xi}\left(\frac{\hat{c}_{k-1}}{\hat{c}_k}-\frac{\hat{c}_k}{\hat{c}_{k+1}} \right)\|y_{k+1}\|^2+\frac{16L_{12}^2}{\xi \hat{c}_k^2}\|x_{k+1}-x_{k} \|^2+\frac{28}{5\xi}\| y_{k+1}-y_k \|^2\nonumber\\
				&\quad-\frac{8}{\xi}\|y_k-y_{k-1}\|^2.
			\end{align}
			Combining \eqref{lem:5.6:6} and \eqref{lem:5.5:1} in Lemma \ref{lem:5.5}, and using the definition of $\mathcal{F}_{k+1}$, we conclude
			\begin{align*}
				\mathcal{F}_{k+1}-\mathcal{F}_{k}
				&\leq -\left( \iota+\frac{\vartheta_k}{2}-\frac{\xi L_{12}^2}{2}-\frac{16L_{12}^2}{\xi \hat{c}_{k}^2}\right)\| x_{k+1}-x_k \| ^2+\frac{\hat{c}_{k-1}-\hat{c}_{k}}{2}\|y_{k+1}\|^2\\
				&\quad-\frac{9}{10\xi}\| y_{k+1}-y_k \|^2+\frac{8}{\xi}\left(\frac{\hat{c}_{k-1}}{\hat{c}_k}-\frac{\hat{c}_k}{\hat{c}_{k+1}} \right)\|y_{k+1}\|^2.
			\end{align*}
		\end{proof}
		
		{\bf The proof of Theorem 5.2}
		\begin{proof}
			By $\xi\le \frac{1}{10L_{22}}$ and $\hat{c}_k=\frac{19}{20\xi k^{\text{1/}4}},\forall k\geq 1$, let us denote $\vartheta_k=\xi L_{12}^2+\frac{16\tau L_{12}^2}{\xi \hat{c}_{k}^2}-2\iota$, $\alpha_k=\frac{8(\tau - 2)L_{12}^2}{\xi \hat{c}_{k}^2}$, we can easily see that the relations in \eqref{lem:5.6:1} are satisfied. It follows from the selection of $\vartheta_k$ and $\alpha_k$ that
			$$\iota+\frac{\vartheta_k}{2}-\frac{\xi L_{12}^2}{2}-\frac{16L_{12}^2}{\xi \hat{c}_{k}^2}=\alpha_k.$$
			This observation, in view of Lemma \ref{lem:5.6}, then immediately implies that
			\begin{align}\label{thm5.2:1}
				\alpha_k\|x_{k+1}-x_k \|^2+\frac{9}{10\xi}\| y_{k+1}-y_k \|^2
				\leq& \mathcal{F}_{k}-\mathcal{F}_{k+1}+\frac{8}{\xi}\left(\frac{\hat{c}_{k-1}}{\hat{c}_k} -\frac{\hat{c}_k}{\hat{c}_{k+1}} \right)\|y_{k+1}\|^2\nonumber\\
				&+\frac{\hat{c}_{k-1}-\hat{c}_k}{2}\|y_{k+1}\|^2.
			\end{align}
			We can easily check from the definition of $\bar{f}_{k} (x_k, y_k)$ that  $$\|\nabla {\mathcal{G}}_k\| -\|\nabla \bar{\mathcal{G}}_k \| \leq \hat{c}_{k}\| y_k \|.$$
			By replacing $f$ with $\bar{f}_{k}$, $\hat{\eta}$ with $\vartheta_k+\iota$, similar to \eqref{thm:blockstr:gx} and \eqref{thm:blockstr:gy}, we immediately obtain that
			\begin{align}\label{thm5.2:2}
				\|(\nabla \bar{\mathcal{G}}_k)_{x^{(i)}} \|
				\le\left(\vartheta_k+\iota+L_{11}\right)\| x_{k+1}-x_k \|,
			\end{align}
			and
			\begin{align}\label{thm5.2:3}
					\| (\nabla \bar{\mathcal{G}}_k)_y \|
					\le&\frac{1}{\xi} \| y_{k+1}-y_k \|+L_{12}\|x_{k+1}-x_k\|.
			\end{align}
			Combining \eqref{thm5.2:2} and \eqref{thm5.2:3}, and using the Cauchy-Schwarz inequality, we have
			\begin{align}\label{thm5.2:4}
					\| \nabla \bar{\mathcal{G}}_k \|^2
					\leq& \left(K_1\left(\vartheta_k+\iota+L_{11}\right)^2+2L_{12}^2\right)\| x_{k+1}-x_k \|^2+\frac{2}{\xi ^2}\|y_{k+1}-y_k \|^2.
			\end{align}
			Since both $\alpha_k$ and $\vartheta_k$ are in the same order when $k$ becomes large enough, it then follows from the definition of $r_1$ that $\forall k\ge 1$,
			\begin{align}\label{thm5.2:5}
				\frac{K_1\left(\vartheta_k+\iota+L_{11}\right)^2+2L_{12}^{2}}{\alpha _k^2} &=  \frac{K_1\left( \xi L_{12}^2+\frac{16\tau L_{12}^2}{\xi (\hat{c}_{k})^2}-\iota+L_{11} \right) ^2+2L_{12}^{2}}{\alpha _k^2}\nonumber\\
				&\leq \frac{2K_1\left(\frac{16\tau L_{12}^2}{\xi (\hat{c}_{k})^2}\right)^2+ 2K_1\left(\xi L_{12}^2-\iota+L_{11}\right)^2 + 2L_{12}^2}{\alpha _k^2}\nonumber\\
				&= \frac{8K_1\tau^2}{(\tau-2)^2} + \frac{19^4\left(2K_1\left(\xi L_{12}^2-\iota+L_{11}\right)^2 + 2L_{12}^2\right)}{64\cdot20^4 \xi^2 (\tau-2)^2 L_{12}^4k} \nonumber\\
				&\leq r_1
			\end{align}
			Combining the previous two inequalities in \eqref{thm5.2:5} and \eqref{thm5.2:4}, we obtain
			\begin{equation}\label{thm5.2:6}
					\| \nabla \bar{\mathcal{G}}_k \|^2
					\leq r_1(\alpha_k)^2\| x_{k+1}-x_k \|^2+ \frac{2}{\xi ^2}\|y_{k+1}-y_k \|^2.
			\end{equation}
			Denote $r_k^{(2)}=\frac{1}{\max \left\{ r_1\alpha_k,\frac{20}{9\xi} \right\}}$.
			By multiplying $r_k^{(2)}$ on the both sides of \eqref{thm5.2:6}, and using \eqref{thm5.2:1}, we have
			\begin{align}\label{thm5.2:7}
					r_k^{(2)}\| \nabla \bar{\mathcal{G}}_k \|^2
					\leq &\mathcal{F}_{k}-\mathcal{F}_{k+1}+\frac{8}{\xi}\left(\frac{\hat{c}_{k-1}}{\hat{c}_k} -\frac{\hat{c}_k}{\hat{c}_{k+1}} \right)\|y_{k+1}\|^2+\frac{\hat{c}_{k-1}-\hat{c}_k}{2}\|y_{k+1}\|^2,
			\end{align}
			where the last inequality follows since $r_{k}^{(2)}$ is a decreasing sequence. Denoting
			$$\bar{\mathcal{T}}(\varepsilon):=\min\{k \mid \|\nabla \bar{\mathcal{G}}(x_k,y_k) \|\leq \frac{\varepsilon}{2}, k\geq 2\},$$
			Summing both sides of \eqref{thm5.2:7} from $k=2$ to $k=\bar{\mathcal{T}}(\varepsilon)$, we then obtain
			\begin{align}\label{thm5.2:8}
					&\sum_{k=2}^{\bar{\mathcal{T}}(\varepsilon)}{r_k^{(2)}\lVert \nabla \bar{\mathcal{G}}_k \rVert ^2}\nonumber\\
					\leq& \mathcal{F}_2-\mathcal{F}_{\bar{\mathcal{T}}(\varepsilon)}+\frac{8}{\xi}\left( \frac{\hat{c}_1}{\hat{c}_2}-\frac{\hat{c}_{\bar{\mathcal{T}}(\varepsilon)}}{\hat{c}_{\bar{\mathcal{T}}(\varepsilon)+1} } \right)\hat{\delta}_y^2+\frac{\hat{c}_1-\hat{c}_{\bar{\mathcal{T}}(\varepsilon)}}{2}\hat{\delta}_y^2\nonumber\\
					\leq&\mathcal{F}_2-\mathcal{F}_{\bar{\mathcal{T}}(\varepsilon)}+\frac{8\hat{c}_1}{\xi \hat{c}_2}\hat{\delta}_y^2+\frac{\hat{c}_1}{2}\hat{\delta}_y^2\nonumber\\
					=& \mathcal{F}_2-\mathcal{F}_{\bar{\mathcal{T}}(\varepsilon)}+ \left(8\cdot 2^{1/4} + \frac{19}{40} \right)\frac{\hat{\delta}_y^2}{\xi}.
			\end{align}
			Note that by the definition of $\mathcal{F}_{k+1}$ in Lemma \ref{lem:5.6}, we have
			\begin{align*}
				\mathcal{F}_{\bar{\mathcal{T}}(\varepsilon)}&\ge\underline{l}-\frac{8\hat{\delta}_y^2}{\xi}\left(1+\frac{1}{\bar{\mathcal{T}}(\varepsilon)}\right)^{1/4}
				-\frac{15\hat{\delta}_y^2}{\xi}-\frac{\hat{c}_1\hat{\delta}_y^2}{2}\\
				&=\underline{l}-\left(2^{13/4}+15+\frac{19}{40}\right)\frac{\hat{\delta}_y^2}{\xi}=\underline{\mathcal{L}},
			\end{align*}
			where $\underline{l} := \min_{(x,y) \in \mathcal{X} \times \mathcal{Y}} l(x,y)$. By the definition of $\mathcal{L}_2$, we then conclude from \eqref{thm5.2:8} that
			\begin{align}\label{thm5.2:8n}
					\sum_{k=2}^{\bar{\mathcal{T}}(\varepsilon)}{r_k^{(2)}\lVert \nabla \bar{\mathcal{G}}_k \rVert ^2}&\le \mathcal{F}_2-\underline{\mathcal{L}}+ \left(8\cdot 2^{1/4} + \frac{19}{40} \right)\frac{\hat{\delta}_y^2}{\xi}\nonumber \\
					&\le \mathcal{L}_2-\underline{\mathcal{L}}+ \left(8\cdot 2^{1/4} + \frac{19}{40} \right)\frac{\hat{\delta}_y^2}{\xi}
					=r_2.
			\end{align}
			
			We can see
			from the selection of $r_3$ that $r_3=\max\{r_1,\frac{20}{9\xi\alpha_2}\}$. Observe that $\alpha_k$ is an increasing sequence, when $k\geq 2$, we have that $r_3\ge\max\{r_1,\frac{20}{9\xi\alpha_k}\}$, which implies that $r_k^{(2)}\geq\frac{1}{r_3\alpha_k}$, by multiplying $r_3$ on the both sides of \eqref{thm5.2:8n}, and combining the definition of $r_2$, we have $	\sum_{k=2}^{\bar{\mathcal{T}}(\varepsilon)} \frac{1}{\alpha_k} \|\nabla \bar{\mathcal{G}}_k \|^2\leq r_2r_3$,
			which, by the definition of $\bar{\mathcal{T}}(\varepsilon)$, implies that
			\begin{equation}\label{thm5.2:9}
				\frac{\varepsilon^2}{4} \leq \frac{r_2r_3}{\sum_{k=2}^{\bar{\mathcal{T}}(\varepsilon)}{\frac{1}{\alpha_k}}}.
			\end{equation}
			Note that when $\hat{c}_k=\frac{19}{20\xi k^{1/4}}$, $\alpha_k = \frac{2\cdot40^2\xi (\tau -2) L_{12}^2 \sqrt{k}}{19^2}$. By using the fact $\sum_{k=2}^{\bar{\mathcal{T}}(\varepsilon)}1/\sqrt{k}\geq \sqrt{\bar{\mathcal{T}}(\varepsilon)}-1$ and \eqref{thm5.2:9}, we conclude $\frac{\varepsilon^2}{4} \leq \frac{2\cdot40^2\xi(\tau -2) L_{12}^2r_2r_3}{19^2\left(\sqrt{\bar{\mathcal{T}}(\varepsilon)}-1\right)}$ or equivalently,
			\begin{equation*}
				\bar{T}(\varepsilon) \le \left( \frac{2\cdot80^2\xi (\tau -2) L_{12}^2r_2r_3 }{19^2\varepsilon ^2}+1 \right) ^2.
			\end{equation*}
			On the other hand, if $k \ge \frac{19^4\hat{\delta}_y^4}{10^4\xi ^4\varepsilon ^4}$, then $c_{k}=\frac{19}{20\xi k^{1/4}} \leq \frac{\varepsilon}{2\hat{\delta}_y}$. This inequality together with the definition of $\hat{\delta}_y$ then imply that $\hat{c}_{k}\| y_k \|\leq \frac{\varepsilon}{2}$. Therefore, there exists a
			\begin{align*}
				\mathcal{T}(\varepsilon) &\leq \max \left(	\bar{\mathcal{T}}(\varepsilon), \frac{19^4\hat{\delta} _y^4}{10^4\xi^4\varepsilon ^4}\right)\\
				&\leq \max \left(	\left( \frac{2\cdot80^2\xi (\tau -2) L_{12}^2r_2r_3 }{19^2\varepsilon ^2}+1 \right) ^2, \frac{19^4\hat{\delta}_y^4}{10^4\xi^4\varepsilon ^4}\right),
			\end{align*}
			such that
			$\|\nabla \mathcal{G}_k\|\leq \| \nabla \bar{\mathcal{G}}_k \|+ \hat{c}_{k}\| y_k \| \leq \frac{\varepsilon}{2}+\frac{\varepsilon}{2} =\varepsilon$.
		\end{proof}
		
		\section{Proof of Theorem 5.3}	
		We prove the following two lemmas before giving the proof of Theorem 5.3.
		
		\begin{lem}\label{lem:6.2}
			Suppose that Assumption 2.1 holds. Let $\{\left(x_{k},y_{k}\right)\}$ be a sequence generated by Algorithm 2
			with parameter settings in (5.13). If $\bar{\xi} > L_{22}$, then we have
			\begin{align}\label{lem:6.2:1}
				&l( x_{k+1},y_{k+1}) -l( x_{k},y_{k})  \nonumber \\
				\geq& \frac{\bar{\xi}}{2}\| y_{k+1}-y_{k} \| ^2-\frac{L_{21}^{2}\varrho}{2}\| y_{k}-y_{k-1} \|^2+(\hat{\theta}-\frac{1}{2\varrho}-\frac{\varrho L_{11}^2}{2})\| x_{k}-x_{k-1}\|^2\nonumber\\
				&+ (\frac{\hat{\theta}}{2}-\frac{1}{\varrho})\|x_{k+1}-x_{k} \|^2.
			\end{align}
		\end{lem}
		
		\begin{proof}
			Similar to \eqref{lem:5.2:2}-\eqref{lem:5.2:4} in the proof of Lemma \ref{lem:5.2}, by the optimality condition for $x_k$ in (5.4) implies that $\forall x\in \mathcal{X}$ and $\forall k\geq 1$,
			\begin{align}
				\langle \nabla _xf(x_{k},y_{k})+\partial h_1(x_{k+1})+\frac{1}{\varrho}(x_{k+1}-x_{k}),x-x_{k+1} \rangle &\ge 0,\label{lem:6.2:2}\\
				\langle \nabla _xf(x_{k},y_{k})+\partial h_1(x_{k+1})+\frac{1}{\varrho}(x_{k+1}-x_{k}) ,x_{k}-x_{k+1} \rangle &\ge 0,\label{lem:6.2:3}\\
				\langle \nabla _xf(x_{k-1},y_{k-1})+\partial h_1(x_{k})+\frac{1}{\varrho}(x_{k}-x_{k-1}),x_{k+1}-x_{k} \rangle &\ge 0,\label{lem:6.2:4}
			\end{align}
			which, in view of the fact that $f\left( x,y  \right)$ is $\hat{\theta} $-strongly convex w.r.t. $x$ for any given $y\in \mathcal{Y}$ and $h_1(x)$ is convex, then implies that
			\begin{align}\label{lem:6.2:5}
				l( x_{k+1},y_{k}) -l( x_{k},y_{k} )\geq &  \langle  \nabla _xf( x_{k},y_{k} )+\partial h_1(x_{k}),x_{k+1}-x_{k}\rangle  +\frac{\hat{\theta}}{2}\| x_{k+1}-x_{k} \|^2\  \nonumber\\
				\geq &\langle  \nabla _xf( x_{k},y_{k} ) -\nabla _xf( x_{k-1},y_{k-1}) ,x_{k+1}-x_{k} \rangle  \nonumber\\
				& -\frac{1}{\varrho}\langle  x_{k}-x_{k-1},x_{k+1}-x_{k} \rangle +\frac{\hat{\theta}}{2}\| x_{k+1}-x_{k} \|^2.
			\end{align}
			The rest proof is the same with that of Lemma 4.1 except replacing $\zeta$ by $\varrho$ and $\theta$ by $\hat{\theta}$. We omit the details here.
		\end{proof}

		\begin{lem}\label{lem:6.3}
			Suppose that Assumption 2.1 holds. Let $\{\left(x_k,y_k\right)\}$ be a sequence generated by Algorithm 2
			with parameter settings in (5.13). Denote
			$$l_{k+1}:=l\left( x_{k+1},y_{k+1} \right),\ \hat{S}_{k+1}:=-\frac{2}{\varrho^2\hat{\theta}}\lVert x_{k+1}-x_k \rVert^2,$$
			$$\hat{F}_{k+1}:=l_{k+1}+\hat{S}_{k+1}+ (\hat{\theta}+\frac{7}{2\varrho}-\frac{\varrho L_{11}^2}{2}-\frac{2L_{11}^2}{\hat{\theta}})\lVert x_{k+1}-x_{k} \rVert^2-(\frac{L_{21}^2\varrho}{2}+\frac{2L_{21}^2}{\hat{\theta}^2\varrho})\|y_{k+1}-y_k\|^2.$$
			If $\bar{\xi}>L_{22}$,
			then $\forall k \geq 1$,
			\begin{align}\label{lem:6.3:2}
				\hat{F}_{k+1}-\hat{F}_k\geq& \left(\frac{\bar{\xi}}{2}-\frac{\varrho L_{21}^2}{2}-\frac{2L_{21}^2}{\varrho \hat{\theta}^2}\right)\lVert y_{k+1}-y_k \rVert ^2+\left(\frac{3\hat{\theta}-\varrho L_{11}^2}{2}+\frac{\hat{\theta}-4\varrho L_{11}^2}{2\varrho\hat{\theta}}\right)\lVert x_{k+1}-x_k \rVert ^2.
			\end{align}
		\end{lem}
		
		\begin{proof}
			First by \eqref{lem:6.2:3} and \eqref{lem:6.2:4}, we have
			\begin{align}\label{lem:6.3:3}
				\frac{1}{\varrho}\langle m_{k+1}, x_k -x_{k+1} \rangle
				\ge \langle \nabla _xf( x_k,y_k ) -\nabla _xf\left( x_{k-1},y_{k-1} \right)+\partial h_1(x_{k+1})-\partial h_1(x_{k}),x_{k+1}-x_k \rangle,
			\end{align}
			where $m_{k+1}:=\left(x_{k+1}-x_k\right)-\left(x_k-x_{k-1}\right)$, together with $\langle \partial h_1(x_{k+1})-\partial h_1(x_k),x_{k+1}-x_k \rangle \geq0$ then imply that
			\begin{align}
				\frac{1}{\varrho}\langle m_{k+1},x_k-x_{k+1}  \rangle \ge \langle \nabla _xf( x_k,y_k ) -\nabla _xf\left( x_{k-1},y_{k-1} \right),x_{k+1}-x_k \rangle.\label{lem:6.3:5}
			\end{align}
			The rest proof is the same with that of Lemma 4.2 except replacing $\zeta$ by $\varrho$ and $\theta$ by $\hat{\theta}$. We omit the details here.
		\end{proof}
		
		{\bf The proof of Theorem 5.3}
		\begin{proof}
			Noting that (5.5) is equivalent to $ y_{k+1}^{(j)}= \operatorname{ Prox}_{g_j,\mathcal{Y}_j}^{\bar{\xi}}\left( y^{(j)}_k+\frac{1}{\bar{\xi}} \nabla _{y^{(j)}}f\left(x_{k+1}, w_{k+1}^{(j)} \right) \right)$, we immediately obtain
			\begin{align}\label{thm:5.3:1}
				&\|(\nabla \mathcal{G}_k)_{x^{(j)}} \|\nonumber\\
				\leq&\bar{\xi}\|\operatorname{Prox}_{g_j,\mathcal{Y}_j}^{\bar{\xi}}\left( y^{(j)}_k+\frac{1}{\bar{\xi}} \nabla _{y^{(j)}}f\left(x_{k+1}, w_{k+1}^{(j)} \right) \right)-\operatorname{ Prox}_{g_j,\mathcal{Y}_j}^{\bar{\xi}}\left( y^{(j)}_k+\frac{1}{\bar{\xi}} \nabla _{y^{(j)}}f\left(x_{k}, y_{k} \right) \right) \|\nonumber\\
				&+\bar{\xi}\| y^{(j)}_{k+1}-y^{(j)}_k \|\nonumber\\
				\leq& \bar{\xi}\| y^{(j)}_{k+1}-y^{(j)}_k \| +\|\nabla _{y^{(j)}}f\left(x_{k+1}, w_{k+1}^{(j)} \right)-\nabla _{y^{(j)}}f\left(x_{k}, y_{k} \right)\|\nonumber\\
				\le&\bar{\xi}\| y^{(j)}_{k+1}-y^{(j)}_k \|+L_{12}\|x_{k+1}-x_k \|+L_{22}\|w_{k+1}^{(j)}-y_k \|\nonumber\\
				\le& \left(\bar{\xi}+L_{22}\right)\| y_{k+1}-y_k \|+L_{12}\|x_{k+1}-x_k \|,
			\end{align}
			where the second inequality holds by the nonexpansiveness of the proximal operator $\operatorname{Prox}_{g_j,\mathcal{Y}_j}^{\bar{\xi}}$, the last inequality hold since $\|y^{(j)}_{k+1}-y^{(j)}_k \|\le\| y_{k+1}-y_k \|$, $\|w^{(j)}_{k+1}-y_k \|\le\| y_{k+1}-y_k \|$ by the definition of $y^{(j)}$ and $w^{(j)}_{k+1}$.
			On the other hand, since (5.4) is equivalent to $ x_{k+1}= \operatorname{ Prox}_{\mathcal{X},h_1}^{1/\varrho}\left( x_k-\varrho \nabla _{x}f\left( x_{k} ,y_{k} \right) \right)$, we conclude from the triangle inequality that
			\begin{align}\label{thm:5.3:2}
				\| (\nabla \mathcal{G}_k)_x \|=\frac{1}{\varrho} \| x_{k+1}-x_k \|.
			\end{align}
			By combining \eqref{thm:5.3:1}-\eqref{thm:5.3:2}, and using Cauchy-Schwarz inequality, we obtain
			\begin{align}\label{thm:5.3:3}
				\| \nabla \mathcal{G}_k \|^2
				\leq& \left(2K_2\left(\bar{\xi}+L_{22}\right)^2\right)\| y_{k+1}-y_k \|^2+\left(\frac{1}{\varrho ^2}+2K_2L_{12}^2\right)\|x_{k+1}-x_k \|^2,
			\end{align}
			Observing that $\hat{r}_1 > 0$. Multiplying both sides of \eqref{thm:5.3:3} by $\hat{r}_1$, and using \eqref{lem:5.3:2} in Lemma \ref{lem:6.3},
			we have
			\begin{align}\label{thm:5.3:4}
				\hat{r}_1\|\nabla \mathcal{G}_k \|^2\le \hat{F}_{k+1}-\hat{F}_k.
			\end{align}
			Summing up the above inequalities from $k=1$ to $k=\mathcal{\mathcal{T}}(\varepsilon)$ and using the definition of $\hat{L}_1$, we obtain
			\begin{align}\label{thm:5.3:5}
				\sum_{k=1}^{\mathcal{T}\left( \varepsilon \right)}{\hat{r}_1\|\nabla \mathcal{G}_k \| ^2}\le \hat{F}_{\mathcal{T}\left( \varepsilon \right)+1}-\hat{F}_1 \le \hat{F}_{\mathcal{T}\left( \varepsilon \right)+1}-\hat{L}_1.
			\end{align}
			Note that by the definition of $\hat{F}_{k+1}$ in Lemma \ref{lem:6.3}, we have
			\begin{align*}
				\hat{F}_{\mathcal{T}\left( \varepsilon \right)+1}
				&\le \hat{l}+(\hat{\theta}+\frac{7}{2\varrho}-\frac{\varrho L_{11}^2}{2}-\frac{2L_{11}^2}{\hat{\theta}})\delta_x^2= \bar{L}.
			\end{align*}
			We then conclude from \eqref{thm:5.3:5} that
			\begin{align*}
				\sum_{k=1}^{\mathcal{T}\left( \varepsilon \right)}{\hat{r}_1\| \nabla \mathcal{G}_k \|^2}\le \bar{L}-\hat{L}_1 ,
			\end{align*}
			which, in view of the definition of $\mathcal{T}(\varepsilon)$, implies that $\varepsilon ^2\le (\bar{L}-\hat{L}_1)/(\mathcal{T}( \varepsilon )\cdot \hat{r}_1)$
			or equivalently, $\mathcal{T}\left( \varepsilon \right) \le (\bar{L}-\hat{L}_1)/(\hat{r}_1\varepsilon ^2)$.
		\end{proof}
		
		\section{Proof of Theorem 5.4}	
		We prove the following three lemmas before giving the proof of Theorem 5.4.
		\begin{lem}\label{lem:6.4}
			Suppose that Assumption 2.1 holds. Let $\{(x_k,y_k)\}$ be a sequence generated by Algorithm 2 with parameter settings in (5.17), If $\forall k,\bar{\vartheta}_k>L_{22}$,  we have
			\begin{equation}\label{lem:6.4:1}
				l(x_{k+1},y_{k+1})-l(x_{k+1},y_k)\ge \left( \bar{\iota} +\frac{\bar{\vartheta}_k}{2} \right) \| y_{k+1}-y_k \|^2.
			\end{equation}
		\end{lem}
		\begin{proof}
			The proof is the same with that of Lemma 5.3 except replacing $\xi^{(j)}_k$ by $\bar{\iota} +\bar{\vartheta}_k$. We omit the details here.
		\end{proof}

		\begin{lem}\label{lem:6.5}
			Suppose that Assumption 2.1 and 5.2 hold. Let $\{\left(x_k,y_k\right)\}$ be a sequence generated by Algorithm 2
			with parameter settings in (5.17). If $\forall k,\bar{\vartheta}_k >L_{22}$ and $\bar{\varrho} \le  \frac{2}{L_{11}^{'}+\hat{b}_1}$, then
			\begin{align}\label{lem:6.5:1}
				&l(x_{k+1},y_{k+1})-l(x_k,y_k)\nonumber \\
				\geq& \left( \bar{\iota} +\frac{\bar{\vartheta}_k}{2} \right) \| y_{k+1}-y_k \| ^2-\frac{L_{21}^{2}\bar{\varrho}}{2}\| y_k-y_{k-1}\|^2-\frac{1}{\bar{\varrho}} \|x_{k+1}-x_k \| ^2-\frac{1}{2\bar{\varrho}}\| x_k-x_{k-1} \|^2\nonumber\\
				&-\frac{\hat{b}_{k-1}}{2}(\| x_{k+1}\| ^2-\| x_k\|^2).
			\end{align}
		\end{lem}
		
		\begin{proof}
			Similar to \eqref{lem:6.2:2}-\eqref{lem:6.2:4} in the proof of Lemma \ref{lem:6.2}, by replacing $l$ with $\bar{l}_{k}$, $f(x_k,y_k)$ with $\bar{f}_{k}$, $\varrho$ with $\bar{\varrho}$ respectively, and setting $\hat{\theta}=0$, we obtain that
			\begin{align}\label{lem:6.5:2}
				\bar{l}_{k}( x_{k+1},y_{k}) -\bar{l}_{k}( x_{k},y_{k} )\geq &  \langle  \nabla _x\bar{f}_{k}( x_{k},y_{k} )+\partial h_1(x_{k}),x_{k+1}-x_{k}\rangle   \nonumber\\
				\geq &\langle  \nabla _x\bar{f}_{k}( x_{k},y_{k} ) -\nabla _x\bar{f}_{k-1}( x_{k-1},y_{k-1}) ,x_{k+1}-x_{k} \rangle  \nonumber\\
				& -\frac{1}{\bar{\varrho}}\langle  x_{k}-x_{k-1},x_{k+1}-x_{k} \rangle .
			\end{align}
			The rest proof is the same with Lemma 4.3 except replacing $\bar{\zeta}$ by $\bar{\varrho}$ and $q_k$ by $\hat{b}_k$. We omit the details here.
			
		\end{proof}

		\begin{lem}\label{lem:6.6}
			Suppose that Assumptions 2.1 and 5.2 hold. Let $\{\left(x_k,y_k\right)\}$ be a sequence generated by Algorithm 2
			with parameter settings in (5.17). Denote
			\begin{align*}
				\hat{\mathcal{S}}_{k+1}&:=-\frac{8}{\bar{\varrho}^2\hat{b}_{k+1}}\| x_{k+1}-x_k \|^2-\frac{8}{\bar{\varrho}}\left(1-\frac{\hat{b}_k}{\hat{b}_{k+1}}\right)\|x_{k+1}\|^2,\\
				\hat{\mathcal{F}}_{k+1}&:=l(x_{k+1},y_{k+1})+\hat{\mathcal{S}}_{k+1}+ \frac{15}{2\bar{\varrho}}\| x_{k+1}-x_k \|^2+\frac{\hat{b}_k}{2}\|x_{k+1}\|^2\\
				&\quad-\left( \frac{\bar{\varrho} L_{21}^2}{2}+\frac{16L_{21}^2}{\bar{\varrho} (\hat{b}_{k+1})^2} \right)\| y_{k+1}-y_k \| ^2.
			\end{align*}
			If
			\begin{equation}\label{lem:6.6:1}
				\bar{\vartheta}_k >L_{22}, ~\frac{1}{\hat{b}_{k+1}}-\frac{1}{\hat{b}_k}\leq \frac{\bar{\varrho}}{5},~\bar{\varrho} \leq \frac{2}{L_{11}^{'}+\hat{b}_1},
			\end{equation}
			then $\forall k \geq 1$,
			\begin{align*}
				&\hat{\mathcal{F}}_{k+1}-\hat{\mathcal{F}}_{k} \nonumber\\
				\geq 	&\left( \bar{\iota}+\frac{\bar{\vartheta}_k}{2}-\frac{\bar{\varrho} L_{21}^2}{2}-\frac{16L_{21}^2}{\bar{\varrho} (\hat{b}_{k+1})^2}\right)\| y_{k+1}-y_k \| ^2+\frac{\hat{b}_k-\hat{b}_{k-1}}{2}\|x_{k+1}\|^2\nonumber\\
				&+\frac{9}{10\bar{\varrho}}\| x_{k+1}-x_k \|^2 +\frac{8}{\bar{\varrho}}\left(\frac{\hat{b}_k}{\hat{b}_{k+1}}- \frac{\hat{b}_{k-1}}{\hat{b}_k} \right)\|x_{k+1}\|^2.
			\end{align*}
		\end{lem}
		
		\begin{proof}
			Similar to \eqref{lem:6.3:3} in the proof of Lemma \ref{lem:6.3}, by replacing $f(x_k,y_k)$ with $\bar{f}_{k}$, $\varrho$ with $\bar{\varrho}$ respectively, we have
			\begin{align*}
				\frac{1}{\bar{\varrho}}\langle m_{k+1}, x_k -x_{k+1} \rangle
				\ge \langle \nabla _x\bar{f}_{k}( x_k,y_k ) -\nabla _x\bar{f}_{k-1}\left( x_{k-1},y_{k-1} \right),x_{k+1}-x_k \rangle,
			\end{align*}
			The rest proof is the same with Lemma 4.4 except replacing $\bar{\zeta}$ by $\bar{\varrho}$ and $q_k$ by $\hat{b}_k$. We omit the details here.
		\end{proof}

		{\bf The proof of Theorem 5.4}
		\begin{proof}
			By $\bar{\varrho}\le \frac{1}{10L_{11}}$,$\hat{b}_k=\frac{19}{20\bar{\varrho} k^{\text{1/}4}},\forall k\geq 1$, let us denote $\bar{\vartheta}_k=\bar{\varrho} L_{21}^2+\frac{16\tau L_{21}^2}{\bar{\varrho} (\hat{b}_{k+1})^2}-2\bar{\iota}$, $p_k=\frac{8(\tau - 2)L_{21}^2}{\bar{\varrho}( \hat{b}_{k+1})^2}$, it can be easily checked that the relations in \eqref{lem:6.6:1} are satisfied. It follows from the selection of $\bar{\vartheta}_k$ and $p_k$ that
			$$\bar{\iota}+\frac{\bar{\vartheta}_k}{2}-\frac{\bar{\varrho} L_{21}^2}{2}-\frac{16L_{21}^2}{\bar{\varrho}(\hat{b}_{k+1})^2}=p_k.$$
			This observation, in view of Lemma \ref{lem:6.6}, then immediately implies that
			\begin{align}\label{thm:5.4:1}
				p_k\|y_{k+1}-y_k \|^2+\frac{9}{10\bar{\varrho}}\| x_{k+1}-x_k \|^2
				\leq& \hat{\mathcal{F}}_{k+1}-\hat{\mathcal{F}}_{k}+\frac{8}{\bar{\varrho}}\left(\frac{\hat{b}_{k-1}}{\hat{b}_k} -\frac{\hat{b}_k}{\hat{b}_{k+1}} \right)\|x_{k+1}\|^2\nonumber\\
				&+\frac{\hat{b}_{k-1}-\hat{b}_k}{2}\|x_{k+1}\|^2.
			\end{align}
			We can easily check from the definition of $\bar{f}_{k}(x_k, y_k)$ that  $$\|\nabla {\mathcal{G}}_k\| -\|\nabla \bar{\mathcal{G}}_k \| \leq \hat{b}_{k}\| x_k \|.$$
			Similar to \eqref{thm:5.3:1} and \eqref{thm:5.3:2} in the proof of Theorem 5.3, by replacing $f$ with $\bar{f}_{k}$, $\nabla {\mathcal{G}}_k$ with $\nabla \bar{\mathcal{G}}_k$, $\varrho$ with $\bar{\varrho}$, $\bar{\xi}$ with $\bar{\vartheta}_k+\bar{\iota}$ respectively, we conclude that
			\begin{align}\label{thm:5.4:2}
				\| \nabla \bar{\mathcal{G}}_k \|^2
				\leq& \left(2K_2\left(\bar{\vartheta}_k+\bar{\iota}+L_{22}\right)^2\right)\| y_{k+1}-y_k \|^2+\left(\frac{1}{\bar{\varrho} ^2}+2K_2L_{12}^2\right)\|x_{k+1}-x_k \|^2.
			\end{align}
			Since both $p_k$ and $\bar{\vartheta}_k$ are in the same order when $k$ becomes large enough, it then follows from the definition of $\hat{d}_1$ that $\forall k\ge 1$,
			\begin{align}\label{thm:5.4:3}
				\frac{2K_2\left(\bar{\vartheta}_k+\bar{\iota}+L_{22}\right)^2}{p_k^2} &=  \frac{2K_1\left(\bar{\varrho} L_{21}^2+\frac{16\tau L_{21}^2}{\bar{\varrho} (\hat{b}_{k+1})^2}-\bar{\iota}+L_{22} \right) ^2}{p_k^2}\nonumber\\
				&\leq \frac{4K_2\left(\frac{16\tau L_{21}^2}{\bar{\varrho} (\hat{b}_{k+1})^2}\right)^2+ 4K_2\left(\bar{\varrho} L_{21}^2-\bar{\iota}+L_{22}\right)^2}{p_k^2}\nonumber\\
				&= \frac{16K_2\tau^2}{(\tau-2)^2} + \frac{19^4\left(4K_2\left(\bar{\varrho} L_{21}^2-\bar{\iota}+L_{22}\right)^2 \right)}{64\cdot20^4 \bar{\varrho}^2 (\tau-2)^2 L_{21}^4k} \nonumber\\
				&\leq \hat{d}_1.
			\end{align}
			Combining the previous two inequalities in \eqref{thm:5.4:2} and \eqref{thm:5.4:3}, we obtain
			\begin{equation}\label{thm:5.4:4}
				\| \nabla \bar{\mathcal{G}}_k \|^2
				\leq \hat{d}_1(p_k)^2\| y_{k+1}-y_k \|^2+\left( \frac{1}{\bar{\varrho}^2}+2K_2L_{12}^{2} \right)\|x_{k+1}-x_k \|^2.
			\end{equation}
			Denote $\hat{d}_k^{(2)}=\frac{1}{\max \left\{ \hat{d}_1p_k,\frac{10+20K_2\bar{\varrho}^2 L_{12}^2}{9\bar{\varrho}} \right\}}$. By multiplying $\hat{d}_k^{(2)}$ on the both sides of \eqref{thm:5.4:4}, and using \eqref{thm:5.4:1}, we get
			\begin{align}\label{thm:5.4:5}
				\hat{d}_k^{(2)}\| \nabla \bar{\mathcal{G}}_k \|^2
				\leq \hat{\mathcal{F}}_{k+1}-\hat{\mathcal{F}}_{k}+\frac{8}{\bar{\varrho}}\left(\frac{\hat{b}_{k-1}}{\hat{b}_k} -\frac{\hat{b}_k}{\hat{b}_{k+1}} \right)\|x_{k+1}\|^2+\frac{\hat{b}_{k-1}-\hat{b}_k}{2}\|x_{k+1}\|^2.
			\end{align}
			Denoting
			$$\bar{\mathcal{T}}(\varepsilon):=\min\{k \mid \|\nabla \bar{\mathcal{G}}(x_k,y_k) \|\leq \frac{\varepsilon}{2}, k\geq 2\},$$
			Summing both sides of \eqref{thm:5.4:5} from $k=2$ to $k=\bar{\mathcal{T}}(\varepsilon)$, we then obtain
			\begin{align}\label{thm:5.4:6}
				&\sum_{k=2}^{\bar{\mathcal{T}}(\varepsilon)}{\hat{d}_k^{(2)}\lVert \nabla \bar{\mathcal{G}}_k \rVert ^2}\nonumber\\
				\leq & \hat{\mathcal{F}}_{\bar{\mathcal{T}}(\varepsilon)+1}-\hat{\mathcal{F}}_2+\frac{8}{\bar{\varrho}}\left( \frac{\hat{b}_1}{\hat{b}_2}-\frac{\hat{b}_{\bar{\mathcal{T}}(\varepsilon)}}{\hat{b}_{\bar{\mathcal{T}}(\varepsilon)+1} } \right)\hat{\delta}_x^2+\frac{\hat{b}_1-\hat{b}_{\bar{\mathcal{T}}(\varepsilon)}}{2}\hat{\delta}_x^2\nonumber\\
				\leq & \hat{\mathcal{F}}_{\bar{\mathcal{T}}(\varepsilon)+1}-\hat{\mathcal{F}}_2+\frac{8\hat{b}_1}{\bar{\varrho} \hat{b}_2}\hat{\delta}_x^2+\frac{\hat{b}_1}{2}\hat{\delta}_x^2\nonumber\\
				= & \hat{\mathcal{F}}_{\bar{\mathcal{T}}(\varepsilon)+1}-\hat{\mathcal{F}}_2+ \frac{8\cdot 2^{1/4} }{\bar{\varrho}}\hat{\delta}_x^2 + \frac{19}{40\bar{\varrho}} \hat{\delta}_x^2.
			\end{align}
			Note that by the definition of $\hat{\mathcal{F}}_{k+1}$ in Lemma \ref{lem:6.6}, we have
			\begin{align*}
				\hat{\mathcal{F}}_{\bar{\mathcal{T}}(\varepsilon)+1}&\le
				\hat{l}+\frac{8\hat{\delta}_x^2}{\bar{\varrho}}\left(1+\frac{1}{\bar{\mathcal{T}}(\varepsilon)}\right)^{1/4}
				+\frac{15\hat{\delta}_x^2}{\bar{\varrho}}+\frac{q_1\hat{\delta}_x^2}{2}\\
				&=\hat{l}+\left(2^{13/4}+15+\frac{19}{40}\right)\frac{\hat{\delta}_x^2}{\bar{\varrho}}= \hat{\mathcal{L}}_0,
			\end{align*}
			where $\hat{l} := \max_{(x,y) \in \mathcal{X} \times \mathcal{Y}} l(x,y)$.  By the definition of $\hat{\mathcal{L}}_2$, we then conclude from \eqref{thm:5.4:6} that
			\begin{align}\label{thm:5.4:6n}
				\sum_{k=2}^{\bar{\mathcal{T}}(\varepsilon)}{\hat{d}_k^{(2)}\lVert \nabla \bar{\mathcal{G}}_k \rVert ^2}&\le\hat{\mathcal{L}}_0-\hat{\mathcal{F}}_2+ \frac{8\cdot 2^{1/4} }{\bar{\varrho}}\hat{\delta}_x^2 + \frac{19}{40\bar{\varrho}} \hat{\delta}_x^2\nonumber\\
				&\le\hat{\mathcal{L}}_0-\hat{\mathcal{L}}_2+ \frac{8\cdot 2^{1/4} }{\bar{\varrho}}\hat{\delta}_x^2 + \frac{19}{40\bar{\varrho}} \hat{\delta}_x^2=\hat{r}_2.
			\end{align}
			Note that $\hat{r}_3=\max\{\hat{d}_1,\frac{10+20K_2\bar{\varrho}^2 L_{12}^2}{9\bar{\varrho}p_2}\}$. Since $p_k$ is an increasing sequence when $k\geq 2$, we have that $\hat{r}_3\ge\max\{\hat{d}_1,\frac{10+20\bar{\varrho}^2 L_{12}^2}{9\bar{\varrho}p_k}\}$, which implies that $\hat{d}_k^{(2)}\geq\frac{1}{\hat{r}_3p_k}$. By multiplying $\hat{r}_3$ on the both sides of \eqref{thm:5.4:6n}, and using the definition of $\hat{r}_2$, we have $	\sum_{k=2}^{\bar{\mathcal{T}}(\varepsilon)} \frac{1}{p_k} \lVert \nabla \bar{\mathcal{G}}_k \rVert ^2\leq \hat{r}_2\hat{r}_3$,
			which, by the definition of $\bar{\mathcal{T}}(\varepsilon)$, implies that
			\begin{equation}\label{thm:5.4:7}
				\frac{\varepsilon^2}{4} \leq \frac{\hat{r}_2\hat{r}_3}{\sum_{k=2}^{\bar{\mathcal{T}}(\varepsilon)}{\frac{1}{p_k}}}.
			\end{equation}
			Using the assumptions $\hat{b}_k=\frac{19}{20\bar{\varrho} k^{1/4}}$ and $p_k = \frac{2\cdot40^2\bar{\varrho} (\tau -2) L_{21}^2 \sqrt{k+1}}{19^2}$, \eqref{thm:5.4:7} and the fact $\sum_{k=2}^{\bar{\mathcal{T}}(\varepsilon)}1/\sqrt{k+1}\geq \sqrt{\bar{\mathcal{T}}(\varepsilon)}-2$ , we conclude that $\frac{\varepsilon^2}{4} \leq \frac{2\cdot40^2\bar{\varrho}(\tau -2) L_{21}^2\hat{r}_2\hat{r}_3}{19^2\left(\sqrt{\bar{\mathcal{T}}(\varepsilon)}-2\right)}$ or equivalently,
			\begin{equation*}
				\bar{\mathcal{T}}(\varepsilon) \le \left( \frac{2\cdot80^2\bar{\varrho} (\tau -2) L_{21}^2\hat{r}_2\hat{r}_3 }{19^2\varepsilon ^2}+2 \right) ^2.
			\end{equation*}
			On the other hand, if $k> \frac{19^4\hat{\delta}_x^4}{10^4\bar{\varrho} ^4\varepsilon ^4}$, then $q_k=\frac{19}{20\bar{\varrho} k^{1/4}} \leq \frac{\varepsilon}{2\hat{\delta}_x}$, this inequality together with the definition of $\hat{\delta}_x$ then imply that $\hat{b}_k\lVert x_k \rVert \leq \frac{\varepsilon}{2}$. Therefore, there exists a
			\begin{align*}
				\mathcal{T}(\varepsilon) \leq \max \left(	\bar{\mathcal{T}}(\varepsilon), \frac{19^4\hat{\delta} _x^4}{10^4\bar{\varrho}^4\varepsilon ^4}\right)
				\leq \max \left(	\left( \frac{2\cdot80^2 \bar{\varrho} (\tau -2) L_{21}^2\hat{r}_2\hat{r}_3}{19^2\varepsilon ^2}+2 \right) ^2, \frac{19^4\hat{\delta} _x^4}{10^4\bar{\varrho}^4\varepsilon ^4}\right),
			\end{align*}
			such that
			$\lVert \nabla \mathcal{G}_k\rVert \leq \lVert \nabla \bar{\mathcal{G}}_k \rVert + \hat{b}_k\lVert x_k \rVert \leq \frac{\varepsilon}{2}+\frac{\varepsilon}{2} =\varepsilon$.
	\end{proof}
\end{appendix}
	
\end{document}